\documentclass{article}
\usepackage[margin=3.5 cm]{geometry}
\usepackage{graphicx, stmaryrd, amssymb}
\usepackage{color}

\graphicspath{ {images/} }
\usepackage{amsmath}
\usepackage{grffile}
\usepackage{physics}
\usepackage{float}
\usepackage{amsthm}

\usepackage{caption}
\usepackage{subcaption}
\usepackage{textcomp}
\usepackage{indentfirst}
\usepackage[T1]{fontenc}
\usepackage{titling}
\usepackage{tikz}
\usepackage{pgf}
\usepackage{amsfonts}
\newtheorem{theorem}{Theorem}

\usepackage{cite}
\theoremstyle{definition}

\theoremstyle{remark}

\title{\Large{From Convex Ideal Polyhedra to Fundamental Domains in $\mathbb{H}^3$}}
\author{Laurel Heck\\Oberlin College}
\date{August 2018}

\begin{document}

\maketitle

\section{Abstract}

Our research goal is to better understand the relationship between the polyhedron and the group associated with a fundamental domain in $\mathbb{H}^3$. In this paper, we will study torsion-free groups and determine a formula for how many edge classes a given abstract polyhedron must have. We will use that result to classify all fundamental domains on the cube with torsion-free groups, including a discussion of the explicit groups associated to those domains. We will then turn to more general fundamental domains and prove a series of results about how properties of the group place restrictions on the edge classes in the quotient manifold. These results give insight into how the polyhedron and the group associated to a fundamental domain interact, as well as offering concrete tools to find fundamental domains.

\section{Introduction}

Fundamental domains are an important tool for studying and understanding hyperbolic manifolds. A fundamental domain is a region that disjointly tiles hyperbolic space under the action of a Kleinian group. The tiling induces a set of identifications between the faces, edges, and vertices of the domain. We can study the manifold resulting from these identifications by considering the fundamental domain. We are particularly interested in fundamental domains that are polyhedral; here we note a few useful properties of fundamental polyhedra arising from the Poincare construction [1] (adapted from Proposition 3.5.1).

\begin{theorem}
A fundamental polyhedron P has the following properties:

\begin{itemize}
    \item[(i)] The faces are arranged in pairs ($\sigma,\sigma')$; to each pair corresponds an element $g\in G$ (called a face identification) such that $g(\sigma)=\sigma'$ and $g(P)\bigcap P=\sigma$.
    \item[(ii)] If a face pairing transformation is elliptic, there is an edge contained in its rotation axis.
    \item[(iii)] To each edge $e$ corresponds an edge relation: $g_1g_2...g_n=g_e$ where either $g_e$ is an elliptic of finite order with rotation axis containing $e$, or $g_e=id$.
   
    *We will refer to the set of edges in a given edge relation as an edge class.
   
    \item[(iv)] The orbit of P under G fills $\mathbb{H}^3$ without overlap on interiors.
    \item[(v)] The face pairing transformations generate G; the edge relations generate the relations in G.
\end{itemize}

\end{theorem}

In this paper, we will study fundamental domains by considering both the abstract polyhedron P and the Kleinian group G. Specifically, we will investigate how the structure of the polyhedron and the properties of the associated group interact.

This question is both inspired and guided by the work of Igor Rivin with regard to convex ideal polyhedra [2]. It should be noted that we are interested in fundamental domains arising from the Poincare construction referenced in Theorem 1, which yields convex domains. Thus understanding convex polyhedra in $\mathbb{H}^3$ is a step towards understanding fundamental polyhedra. In particular, Rivin's strategy determines whether or not a given abstract polyhedron with certain prescribed dihedral angles can be realized as a convex ideal polyhedron in $\mathbb{H}^3$. We will use Rivin's result as a tool and also generalize his method of imposing conditions on an abstract polyhedron to study fundamental domains.

Rivin's work is an extension of a result by E.M. Andre'ev that characterizes finite-volume convex polyhedra with non-obtuse angles in $\mathbb{H}^3$ [3]. Rivin generalizes Andre'ev's result by describing a method of constructing convex ideal polyhedra with any dihedral angles. Before stating the result, we will briefly touch on the idea of dual polyhedra, a concept that Rivin uses heavily.

Given some polyhedron $P$, the Poincare dual $P^*$ is another polyhedron whose faces correspond to the vertices of the original polyhedron, and vice versa. For example, an edge between two faces A and B in $P$ becomes an edge connecting two vertices $A^*$ and $B^*$ in $P^*$.
\begin{figure}[H]
    \centering

\tikzset{every picture/.style={line width=0.75pt}} 

\begin{tikzpicture}[x=0.55pt,y=0.55pt,yscale=-1,xscale=1]

\draw [color={rgb, 255:red, 245; green, 166; blue, 35 }  ,draw opacity=1 ][line width=2.25]    (443.7,33) -- (535.5,133.8) ;

\draw [color={rgb, 255:red, 74; green, 144; blue, 226 }  ,draw opacity=1 ][line width=2.25]    (443.7,33) -- (414.9,151.48) ;

\draw [line width=2.25]    (443.7,33) -- (355.5,109.04) ;

\draw [line width=2.25]  [dash pattern={on 2.53pt off 3.02pt}]  (443.7,33) -- (470.7,89.59) ;

\draw [color={rgb, 255:red, 80; green, 227; blue, 194 }  ,draw opacity=1 ][line width=2.25]    (414.9,151.48) -- (535.5,133.8) ;

\draw [line width=2.25]    (355.5,109.04) -- (414.9,151.48) ;

\draw [line width=2.25]    (535.5,133.8) -- (432.9,201) ;

\draw [line width=2.25]    (414.9,151.48) -- (432.9,201) ;

\draw [line width=2.25]    (470.7,89.59) -- (535.5,133.8) ;

\draw [line width=2.25]    (355.5,109.04) -- (470.7,89.59) ;

\draw [line width=2.25]  [dash pattern={on 2.53pt off 3.02pt}]  (470.7,89.59) -- (432.9,201) ;

\draw [line width=2.25]    (355.5,109.04) -- (432.9,201) ;

\draw  [line width=2.25]  (82,89.1) -- (126.1,45) -- (232.5,45) -- (232.5,147.9) -- (188.4,192) -- (82,192) -- cycle ; \draw  [line width=2.25]  (232.5,45) -- (188.4,89.1) -- (82,89.1) ; \draw  [line width=2.25]  (188.4,89.1) -- (188.4,192) ;
\draw [line width=2.25]  [dash pattern={on 2.53pt off 3.02pt}]  (82,192) -- (124.75,147.9) ;

\draw  [dash pattern={on 2.53pt off 3.02pt}][line width=2.25]   (124.75, 45) rectangle (232.5, 147.9)   ;
\draw [color={rgb, 255:red, 74; green, 144; blue, 226 }  ,draw opacity=1 ][line width=2.25]    (82,89.1) -- (188.4,89.1) ;

\draw [color={rgb, 255:red, 80; green, 227; blue, 194 }  ,draw opacity=1 ][line width=2.25]    (188.4,89.1) -- (188.4,192) ;

\draw [color={rgb, 255:red, 245; green, 166; blue, 35 }  ,draw opacity=1 ][line width=2.25]    (188.4,89.1) -- (232.5,45) ;

\draw    (259,116) -- (311.5,116) ;
\draw [shift={(313.5,116)}, rotate = 180] [color={rgb, 255:red, 0; green, 0; blue, 0 }  ][line width=0.75]    (10.93,-3.29) .. controls (6.95,-1.4) and (3.31,-0.3) .. (0,0) .. controls (3.31,0.3) and (6.95,1.4) .. (10.93,3.29)   ;

\draw (478,108) node  [align=left] {B};
\draw (443,17) node  [align=left] {T};
\draw (434,218) node  [align=left] {Bo};
\draw (339,107) node  [align=left] {L};
\draw (557,133) node  [align=left] {R};
\draw (432,136) node  [align=left] {F};
\draw (184,74) node  [align=left] {V};
\draw (68,38) node  [align=left] {P};
\draw (525,45) node  [align=left] {P*};

\end{tikzpicture}

    \caption{The cube and octahedron are dual polyhedra. Edges go to edges, faces to vertices, and vertices to faces. T refers to the top face in the cube, Bo refers to bottom, etc.}
    \label{fig:my_label}
\end{figure}
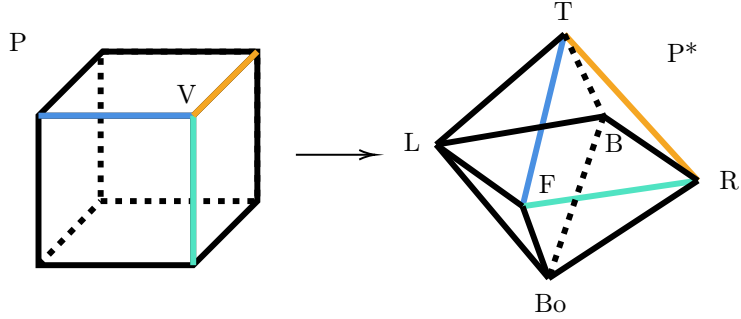

Figure 1 shows another example that will be useful for our purposes: the blue, teal, and orange edges meet at vertex $V$ in $P$. Vertex $V$ is dual to the face in $P^*$ bounded by the duals of those edges. This idea will be useful when we discuss Rivin's equations for ideal convex polyhedra.

Rivin's result is as follows [2, Theorem 0.1]. For $P$ a convex ideal polyhedron in $\mathbb{H}^3$, consider the dual polyhedron $P^*$. For each edge $e^*$ of $P^*$, assign a weight $w(e^*)$, where $w(e^*)$ is the exterior dihedral angle at the dual edge $e\in P$. Then the following conditions must be true:

\begin{itemize}
    \item[1)] $0<w(e^*)<\pi$ for all $e^*\in P^*$
    \item[2)] If $e_1^*,...,e_k^*$ bound a face in $P^*$, then $w(e_1^*)+...+w(e_k^*)=2\pi$.
    \item[3)] If $e_1^*,...,e_k^*$ form a simple circuit that does not bound a face of $P^*$, then $w(e_1^*)+...+w(e_k^*)>2\pi$.
\end{itemize}

His result further states that given some abstract polyhedron $P^*$ with weighted edges $e_i^*$ that satisfy the above conditions, the dual of $P^*$ can be realized as a convex ideal polyhedron P with exterior dihedral angles determined by the weights. In this paper, we will refer to this result as \textit{the Rivin equations and inequalities}. The most useful part comes from condition 2: we need the sum of edge weights for edges bounding a face in the dual to be $2\pi$. Edges bounding a face in the dual become edges incident to a vertex in the original polyhedron. Thus, we will require that the sum of exterior dihedral angles for all edges incident to a given vertex is $2\pi$. An example is shown in Figure 1: the blue, teal, and orange edges meet at vertex $V$ in $P$. Vertex $V$ is dual to the face in $P^*$ bounded by the duals of those edges.

We will combine Rivin's results with several restrictions on the types of groups and polyhedra that we consider in order to obtain results about fundamental domains. First, we will assume that we are working with ideal polyhedra so that we can invoke Rivin's result. Second, we will assume that the groups have no elliptic generators. Finally, we will assume that the interior angle sum of edges in a class will always be $2\pi$. This restriction will give rise to what we will call the fundamental domain equations.

\vspace{1em}
\noindent \textbf{Def}: For an abstract polyhedron with $\overline{E}$ edges and $k$ edge classes $E_1,...,E_k$, where $E_i$ has $n_i$ elements, we obtain $k$ equations of the form
\begin{equation}
    \sum_{j=1}^{n_i} (\pi-x_j)=2\pi,
\end{equation} where $x_j$ is the exterior dihedral angle at edge $j$, and edge $j$ is in class $E_i$. These are the fundamental domain equations.

From these assumptions, we will proceed as follows. In Section 3, we will make a general statement about the role of elliptic elements in general fundamental domains. We will then justify that our assumptions restrict our consideration to torsion-free groups. In Section 4, we will explore how the Rivin equations and inequalities combine to give certain conditions on the generators and relators of a group associated with a general abstract polyhedron. In section 5 we will consider the case that the abstract polyhedron is a cube and classify all possible fundamental domains on the cube with torsion free groups. Finally, in Section 6, we will explore the question of how properties of the group restrict the combinatorics of the associated abstract polyhedron.

\section{Role of Elliptics in Fundamental Domains}

Imagine that we are given a fundamental polyhedron with associated abstract polyhedron $P$ and group $G$. We will consider the possibility that the group has elliptic elements and determine what role they play. Elliptic elements must fix an axis in the interior of $\mathbb{H}^3$, and we will use this fact to characterize them. First, we will show that if the group does have elliptic elements, the axis of each elliptic must only go along edges of the polyhedron and its copies. We will then fully classify the ways in which elliptic elements can appear in the group $G$. Finally, we will argue that by assuming that there are no elliptic generators and that the interior angles in an edge class sum to $2\pi$, we obtain a torsion free group.

\subsection{Elliptic elements must fix an edge}

In order to understand the role of elliptic elements, we make the following claim:

\textbf{Claim 3.1}: The axis of an elliptic element of G must only intersect P and its copies along edges.

\begin{proof} We first quickly note what is meant by the copies of P. Each element of G is an isometry of $\mathbb{H}^3$ that moves the original polyhedron P to another location in $\mathbb{H}^3$ while maintaining its geometry. These copies together with the original polyhedron P tile $\mathbb{H}^3$. Now to prove the claim: we know that the axis is in the interior of $\mathbb{H}^3$. Thus, it must intersect some set of interiors, faces, edges, and vertices of P and its copies. We will prove the statement by ruling out all cases except for edges. We will call the elliptic element $g_e$ and denote its order by $|g_e|.$

Case 1: The axis intersects the interior of P or one of its copies.

Under the action of G, P disjointly tiles all of $\mathbb{H}^3$, and its interior must move around $\mathbb{H}^3$.
Thus, the axis cannot be in the interior of P or any of its copies because the axis is fixed by $g_e$ and would not be able to move as the polyhedron tiles the space.

Case 2: The axis intersects a face of one of the polyhedra.

If $|g_e|>2$, then we have a contradiction: $g_e$ rotates by some angle $\theta\neq\pi$, a rational divisor of $2\pi$, and the axis in face $F$ is fixed. Thus, $g_e(F)\bigcap F\neq0$; namely, the axis is the intersection. However, because $\theta\neq\pi$, the faces are not lined up. Thus, we know that $g_e(F)$ intersects the interior of two copies of $P$, namely the two copies that share face $F$ (see Figure 2(a)). The tiling should be disjoint, so this is a contradiction.

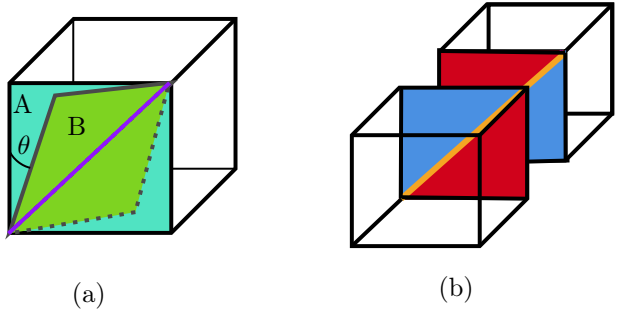
\begin{figure}[H]
    \centering

\tikzset{every picture/.style={line width=0.75pt}} 

\begin{tikzpicture}[x=0.5pt,y=0.5pt,yscale=-1,xscale=1]

\draw    (130.26, 42) rectangle (253, 155.4)   ;
\draw [rotate around= { 359.79: (143.37, 147.3)
    }] [fill={rgb, 255:red, 80; green, 227; blue, 194 }  ,fill opacity=1 ]  (82, 90.6) rectangle (204.74, 204)   ;
\draw  [line width=1.5]  (82,90.6) -- (130.6,42) -- (253,42) -- (253,155.4) -- (204.4,204) -- (82,204) -- cycle ; \draw  [line width=1.5]  (253,42) -- (204.4,90.6) -- (82,90.6) ; \draw  [line width=1.5]  (204.4,90.6) -- (204.4,204) ;
\draw [color={rgb, 255:red, 144; green, 19; blue, 254 }  ,draw opacity=1 ][line width=1.5]    (204.74,90.6) -- (82,204) ;

\draw  [color={rgb, 255:red, 74; green, 74; blue, 74 }  ,draw opacity=1 ][fill={rgb, 255:red, 126; green, 211; blue, 33 }  ,fill opacity=1 ][line width=1.5]  (116.44,99.97) -- (205.15,90.16) -- (82,204) -- cycle ;
\draw  [color={rgb, 255:red, 74; green, 74; blue, 74 }  ,draw opacity=1 ][fill={rgb, 255:red, 126; green, 211; blue, 33 }  ,fill opacity=1 ][dash pattern={on 1.69pt off 2.76pt}][line width=1.5]  (177.26,188.25) -- (81.25,204.44) -- (202.68,90.93) -- cycle ;
\draw [color={rgb, 255:red, 144; green, 19; blue, 254 }  ,draw opacity=1 ][line width=1.5]    (203.08,90.49) -- (81.25,204.44) ;

\draw  [draw opacity=0] (97.73,154.74) .. controls (94.69,154.08) and (91.49,152.57) .. (88.58,150.25) .. controls (85.38,147.71) and (83.05,144.62) .. (81.77,141.53) -- (95.66,141.34) -- cycle ; \draw   (97.73,154.74) .. controls (94.69,154.08) and (91.49,152.57) .. (88.58,150.25) .. controls (85.38,147.71) and (83.05,144.62) .. (81.77,141.53) ;
\draw  [line width=1.5]   (444.45, 31) rectangle (538, 113.73)   ;
\draw  [line width=1.5]  (407.67,66.45) -- (443.12,31) -- (538,31) -- (538,113.73) -- (502.55,149.18) -- (407.67,149.18) -- cycle ; \draw  [line width=1.5]  (538,31) -- (502.55,66.45) -- (407.67,66.45) ; \draw  [line width=1.5]  (502.55,66.45) -- (502.55,149.18) ;
\draw [color={rgb, 255:red, 144; green, 19; blue, 254 }  ,draw opacity=1 ][line width=1.5]    (501.22,66.45) -- (407.67,149.18) ;

\draw [color={rgb, 255:red, 144; green, 19; blue, 254 }  ,draw opacity=1 ][line width=1.5]    (499.96,66.37) -- (407.1,149.5) ;

\draw  [color={rgb, 255:red, 208; green, 2; blue, 27 }  ,draw opacity=1 ][fill={rgb, 255:red, 208; green, 2; blue, 27 }  ,fill opacity=1 ] (500.7,66.72) -- (407.43,149.18) -- (406.8,67.45) -- cycle ;
\draw  [color={rgb, 255:red, 74; green, 144; blue, 226 }  ,draw opacity=1 ][fill={rgb, 255:red, 74; green, 144; blue, 226 }  ,fill opacity=1 ] (409.04,149.52) -- (502.33,67.09) -- (502.63,149.18) -- cycle ;
\draw [color={rgb, 255:red, 245; green, 166; blue, 35 }  ,draw opacity=1 ][line width=2.25]    (502.43,67.6) -- (409.27,149.52) ;

\draw  [color={rgb, 255:red, 208; green, 2; blue, 27 }  ,draw opacity=1 ][fill={rgb, 255:red, 208; green, 2; blue, 27 }  ,fill opacity=1 ] (378,178) -- (473.32,94.83) -- (473.62,177.65) -- cycle ;
\draw  [color={rgb, 255:red, 74; green, 144; blue, 226 }  ,draw opacity=1 ][fill={rgb, 255:red, 74; green, 144; blue, 226 }  ,fill opacity=1 ] (473.33,94.63) -- (379,177) -- (378.37,95.36) -- cycle ;
\draw    (378.8, 94.03) rectangle (472.96, 178.1)   ;
\draw [line width=1.5]    (378.8,178.1) -- (341,214) ;

\draw [color={rgb, 255:red, 245; green, 166; blue, 35 }  ,draw opacity=1 ][line width=2.25]    (472.96,94.03) -- (378.35,178.19) ;

\draw  [line width=1.5]  (341,130.02) -- (376.99,94.03) -- (473.74,94.03) -- (473.74,178.01) -- (437.75,214) -- (341,214) -- cycle ; \draw  [line width=1.5]  (473.74,94.03) -- (437.75,130.02) -- (341,130.02) ; \draw  [line width=1.5]  (437.75,130.02) -- (437.75,214) ;
\draw [line width=1.5]    (407.03,66.27) -- (502.43,67.6) ;

\draw [line width=1.5]    (502.43,67.6) -- (502.73,149.18) ;

\draw [line width=1.5]    (475,150) -- (505.65,149.17) ;

\draw [line width=1.5]    (407.03,64.16) -- (407,95) ;

\draw [line width=1.5]    (379,177) -- (472.96,178.1) ;

\draw [line width=1.5]    (378.37,95.36) -- (378,178) ;

\draw (92,106) node  [align=left] {A};
\draw (133,123) node  [align=left] {B};
\draw (93,137) node   {$\theta $};
\draw (142,251) node  [align=left] {(a)};
\draw (420,247) node  [align=left] {(b)};

\end{tikzpicture}

    \caption{(a) If $|g_e|>2$, then $g_e$ takes $A$ to $B$ and $B$ intersects the original polyhedron in the interior. (b) If $|g_e|=2$, the colored faces are glued together, but with a reflection over the axis. The colored face is identified with two faces, a contradiction.}
    \label{fig:my_label}
\end{figure}

If $|g_e|=2$, we cannot use the same argument because $F$ and $g_e(F)$ will be completely lined up, but with a reflection over the axis (see Figure 2(b)).

The axis breaks the face $F$ into two pieces, and we can consider them as two new faces that have an exterior dihedral angle of $\pi$. Because $P$ is a fundamental domain, all of the elements in the group $G$ send faces to faces. We know that in the $\pi$ rotation over the axis, these two new faces are sent to each other; we have just reflected the face $F$ over the axis. However, in doing so, we have identified face F with itself. From Theorem 1, face identifications come in pairs for a fundamental polyhedron, so now $F$ is identified with itself and the original paired face, a contradiction.

Case 3: The axis is tangent to a vertex of one of the polyhedra.

The axis must continue through $\mathbb{H}^3$ outside of the polyhedron containing the tangent vertex. It cannot go in the interior or along a face of another copy of $P$, so it must only go along edges.

Thus, we conclude that if $G$ has an elliptic element, its axis must only intersect $P$ and its copies along edges.
\end{proof}

\subsection{Identifying All Elliptic Elements in G}

We will first consider the case for a general polyhedron and then discuss what happens when we impose the assumptions stated in Section 2.
\subsubsection{In General}
Suppose that $G$ has an elliptic element $g_e$. By Claim 3.1, the axis of $g_e$ only goes along edges of $P$ and its copies. First consider the case that the axis contains an edge of $P$. There are two options: either $g_e$ is a face identification, or it is not. In the second case, say that $g_e$ takes face $F$ of $P$ to face $F'$ of $P'$. Face $F$ is \textit{not} identified with $F'$ by a generator in this case (see Figure 3).

\begin{figure}[H]
    \centering

\tikzset{every picture/.style={line width=0.75pt}} 

\begin{tikzpicture}[x=0.5pt,y=0.5pt,yscale=-1,xscale=1]

\draw    (245,74) -- (152,209) ;

\draw    (245,74) -- (229,223) ;

\draw    (245,74) -- (364,173) ;

\draw    (245,74) -- (399,43) ;

\draw    (397,117) -- (245,74) ;

\draw    (245,74) -- (300,217) ;

\draw  [color={rgb, 255:red, 208; green, 2; blue, 27 }  ,draw opacity=1 ][fill={rgb, 255:red, 208; green, 2; blue, 27 }  ,fill opacity=1 ]  (245, 74) circle [x radius= 5.5, y radius= 5.5]  ;
\draw [line width=1.5]    (512.5,115.5) -- (512.03,160) ;
\draw [shift={(512,163)}, rotate = 270.6] [color={rgb, 255:red, 0; green, 0; blue, 0 }  ][line width=1.5]    (14.21,-4.28) .. controls (9.04,-1.82) and (4.3,-0.39) .. (0,0) .. controls (4.3,0.39) and (9.04,1.82) .. (14.21,4.28)   ;

\draw (261,180) node   {$P_{1}$};
\draw (311,169) node   {$P_{2}$};
\draw (351,134) node   {$P_{3}$};
\draw (201,178) node   {$P$};
\draw (378,80) node   {$P^{'}$};
\draw (511,187) node   {$g_{0} g_{1} g_{2} g_{3} =g_{e}$};
\draw (526,68) node  [align=left] {$\displaystyle \exists $ an elliptic taking\\P to P' $\displaystyle \Longrightarrow $We can\\construct it via a \\series of face IDs};
\draw (229,234) node  [align=left] {F};
\draw (414,43) node  [align=left] {F'};

\end{tikzpicture}

    \caption{Constructing $g_e$ as a sequence of face identifications. $g_i$} is the face identification that ends at polyhedron copy $P_i$.
    \label{fig:my_label}
\end{figure}

We are looking at the figure from the side; the red dot represents the edge $e$ fixed by $g_e$, coming out of the page. All of the space between $P$ and $P'$ is filled in by other copies of $P$ that also share edge $e$. Thus, we can construct $g_e$ as a sequence of face identifications, as shown in Figure 3.

We begin at face $F$ and construct the relation corresponding to the edge class of $e$ using face identifications.

There are three cases to consider. Starting at face F, the relation for the class of $e$ terminates before we get to $F'$, at $F'$, or after $F'$. Note that the relation terminates when we return to the original edge.

The relation cannot end after $F'$ because we return to the original edge at $F'$, and thus the relation has to end at or before $F'$.

If the edge class finishes at $F'$, then the relator is either $g_e$ itself or the identity.

Finally, if the edge class finishes before $F'$, then by Theorem 1, $g_e$ is a power of the relation, which itself is a smaller elliptic. To see this, note that by discreteness, there must be some smallest rotation element in $G$. Call it $g_{\alpha}$ with rotation by an angle $\alpha$. We claim that $g_e$ must be a power of $g_{\alpha}$. If it was not, then $\exists m\in\mathbb{Z}$ such that $g_{\alpha}^m$ is within $\alpha$ of $F'$. But then $g_e\circ g_{\alpha}^{-m}$ is a smaller rotation element than $g_{\alpha}$, a contradiction.

We conclude that any elliptic that fixes an edge in $P$ is either a generator, a relator, or a finite power of a relator.

To get an elliptic that fixes an edge in one of the copies of $P$, we simply conjugate one of the elliptic elements described above. The generators and relators in any copy of $P$ are conjugate to the generators and relators of P. Conjugacy preserves type of Mobius transformation, so these elements are also elliptic elements. Conversely, any elliptic fixing an edge in one of the copies of P is conjugate to a generator, relator, or power of a relator in P, and we could have found it by conjugating one of those elements.

We conclude that any elliptic in $G$ is either a generator, a relator, a finite power of a relator, or something conjugate to one of the above.

\subsubsection{Under our Assumptions}

We have shown that an elliptic can be constructed as a generator, a relator, or a finite power of a relator, depending on where the relator terminates relative to face $F'$. When we impose that there are no elliptic generators and that the fundamental domain equations are satisfied, we eliminate both of these cases. Clearly if there are no elliptic generators, then $g_e$ is not a generator. If the relator ends exactly at $F'$, then either $g_e$ is the relation or the identity. From the fundamental domain equations, we  know that the interior angle sum of the edges in the class is $2\pi$. Thus, $g_e$ must be the identity. This eliminates the case that one of the relators is an elliptic.

Now we claim that the relator cannot end before we reach $F'$; for if it did, then the elliptic that takes $F$ to $F'$ must have a rotation angle $>2\pi$, which is not possible. This eliminates the case where the elliptic is a finite power of a relator.

With our two assumptions, we have eliminated the cases that the elliptic is a generator, a relator, or a finite power of a relator. Since there are no elliptics with axes containing edges of P, there are no elliptics in the group: suppose that there was an elliptic fixing an edge in another copy of P. Then that elliptic is conjugate to an elliptic fixing an edge in P. But there are no elliptics fixing an edge in P. Thus, we conclude that under these assumptions, there are no elliptics in $G$.

\noindent\textbf{Corollary 3.1:} A group corresponding to a fundamental domain with no elliptic generators and where all edge classes sum to $2\pi$ is a torsion-free group.

\section{Number of Edge Classes for a Given Combinatorial Polyhedron}

Recall that we have two sets of equations: the fundamental domain equations and the Rivin equations. We can combine them in an interesting way to determine several facts about how many edge classes, faces and vertices should be in the quotient manifold resulting from the identifications.

First, we will say a few words about the Euler characteristic equation. For a given spherical polyhedron, we have the equation
\begin{equation} \overline{V}-\overline{E}+\overline{F}=2,
\end{equation} where $\overline{V}$, $\overline{E}$, and $\overline{F}$ are the number of vertices, edges, and faces, respectively. Once we take the quotient, we have the following equation for the resulting manifold:
\begin{equation} V-E+F-P=q,
\end{equation} where again V, E, and F are the number of vertices, edges, and faces, respectively. P is the number of interiors and is equal to 1 in our case; $q$ is the Euler characteristic of the quotient manifold.

To make a few connections to our terminology, we have the following:
\begin{itemize}
    \item $F=\dfrac{\overline{F}}{2}$ because we identify faces pairwise (Theorem 1)
    \item $\overline{V}$ is the number of Rivin equations that we obtain: each vertex is dual to a face, and the incident edges bound the dual face.
    \item $\overline{E}$ is the number of variables in our equations.
    \item $E$ is the number of edge classes and consequently the number of fundamental domain equations.
\end{itemize}

Thus, in general we have a system of up to $\overline{V}+E$ equations in $\overline{E}$ variables (there may be fewer equations if there are edge classes completely at one vertex). We will derive an equation for the number of edge classes given the number of edges and vertices in the original polyhedron.

\textbf{Proposition 4.1}: A polyhedron with $\overline{E}$ edges and $\overline{V}$ vertices must have $E=\dfrac{\overline{E}-\overline{V}}{2}$ edge classes.
\begin{proof}

We will proceed by summing the exterior dihedral angles of all of the edges using each of the sets of equations.

\vspace{1em}

\noindent \textbf{Method 1}: Rivin Equations\vspace{1em}

Each edge is adjacent to two vertices. Thus, if we sum over the angle sums at all of the vertices, we will double-count each edge. Defining $v_i$ to be the sum of the exterior dihedral angles at a given vertex, we have
\begin{equation}
    \sum_{i=0}^{\overline{V}}v_i=2(x_1+...+x_{\overline{E}})=2\pi\overline{V},
\end{equation} from which we obtain that \begin{equation}
    \sum_{j=1}^{\overline{E}}x_j=\pi\overline{V}.
\end{equation}
\vspace{1em}
\noindent \textbf{Method 2}: Fundamental Domain Equations

The edge classes disjointly partition the set of edges in the original polyhedron. Thus, if we sum the exterior dihedral angles in each edge class and then sum over the edge classes, we will obtain the sum of all exterior dihedral angles in the polyhedron. In general, for an edge class consisting of $n_i$ edges, the sum of the exterior angles will be $(n_i-2)\pi$. Summing over all $E$ classes, we obtain
\begin{equation}
    \sum_{i=1}^{E}(n_i-2)\pi=x_1+...+x_{\overline{E}}.
\end{equation}
Rearranging the sum, we find
\begin{equation}
    \pi\sum_{i=1}^{E}n_i-2\pi E=\sum_{j=1}^{\overline{E}}x_j.
\end{equation} Note that $\sum_{i=1}^{E}n_i=\overline{E}$.

\vspace{1em}
Because any fundamental domain must satisfy the Rivin equations and the fundamental domain equations, these two quantities are equal. Cancelling a factor of $\pi$ on each side, we obtain \begin{equation} \overline{E}-2E=\overline{V},
\end{equation}which is what we set out to prove.
\end{proof}

From this result, we conclude the following for the Platonic solids:
\begin{table}[H]
\centering
\begin{tabular}{|c|c|c|c|}
\hline
\textbf{Polyhedron} & \textbf{Vertices} & \textbf{Edges} & \textbf{Edge Classes in FD} \\ \hline
Tetrahedron         & 4                 & 6              & 1                           \\ \hline
Cube                & 8                 & 12             & 2                           \\ \hline
Octahedron          & 6                 & 12             & 3                           \\ \hline
Dodecahedron        & 20                & 30             & 5                           \\ \hline
Icosahedron         & 12                & 30             & 9                           \\ \hline
\end{tabular}
\end{table}

We can push this result a bit farther by combining it with the Euler characteristic equations.

\begin{figure}[h]
\centering

\tikzset{every picture/.style={line width=0.75pt}} 

\begin{tikzpicture}[x=0.75pt,y=0.75pt,yscale=-1,xscale=1]

\draw    (221.5,130) -- (408.5,131) ;

\draw (312,95) node  [align=left] {$\displaystyle 2V-2E+2F-2P=2q$};
\draw (312,114) node  [align=left] {$\displaystyle \overline{V} \ -\ \overline{E} +\ \ \overline{F} \ \ \ \ \ \ \ \ \ =2$};
\draw (316,143) node   {$2V-2E-2P-\overline{V} +\overline{E} =2q-2$};
\draw (215,114) node   {$-$};

\end{tikzpicture}

\end{figure}To simplify, note that $P=1$, so we have
\begin{equation}
    2V-2E-\overline{V}+\overline{E}=2q.
\end{equation}
Then we can rearrange to obtain
\begin{equation}
    \overline{E}-\overline{V}=2E-(2V+2q).
\end{equation} The part of the equation outside of parentheses is our equation from Proposition 4.1; we conclude that $V=q$ and thus that the number of vertex classes in a fundamental domain is equal to the Euler characteristic of the quotient space.

\section{Fundamental Domains on the Cube with Torsion-Free Groups}
Using several of the results from the previous section, we will classify all fundamental domains with the cube as an abstract polyhedron.
\subsection{Finding Candidates for Fundamental Domains on the Cube with Torsion-Free Groups}
\begin{figure}[H]
    \centering

\tikzset{every picture/.style={line width=0.75pt}} 

\begin{tikzpicture}[x=0.75pt,y=0.75pt,yscale=-1,xscale=1]

\draw [color={rgb, 255:red, 245; green, 166; blue, 35 }  ,draw opacity=1 ][line width=2.25]    (454.5,104.42) -- (522.5,104.42) ;

\draw [color={rgb, 255:red, 208; green, 2; blue, 27 }  ,draw opacity=1 ][line width=2.25]    (328.4,99.61) -- (328.4,194.06) ;

\draw [color={rgb, 255:red, 208; green, 2; blue, 27 }  ,draw opacity=1 ][line width=2.25]    (229.41,99.61) -- (328.4,99.61) ;

\draw [color={rgb, 255:red, 189; green, 16; blue, 224 }  ,draw opacity=1 ][line width=2.25]    (228.5,195) -- (328.4,194.06) ;

\draw [color={rgb, 255:red, 208; green, 2; blue, 27 }  ,draw opacity=1 ][line width=2.25]    (273,59) -- (229.41,99.61) ;

\draw [color={rgb, 255:red, 189; green, 16; blue, 224 }  ,draw opacity=1 ][line width=2.25]    (328.4,99.61) -- (362,59) ;

\draw [color={rgb, 255:red, 189; green, 16; blue, 224 }  ,draw opacity=1 ][line width=2.25]    (273,59) -- (362,59) ;

\draw [color={rgb, 255:red, 189; green, 16; blue, 224 }  ,draw opacity=1 ][line width=2.25]    (362,59) -- (362,161.94) ;

\draw [color={rgb, 255:red, 208; green, 2; blue, 27 }  ,draw opacity=1 ][line width=2.25]    (362,161.94) -- (328.4,194.06) ;

\draw [color={rgb, 255:red, 208; green, 2; blue, 27 }  ,draw opacity=1 ][line width=2.25]    (273,59) -- (273,161.94) ;

\draw [color={rgb, 255:red, 189; green, 16; blue, 224 }  ,draw opacity=1 ][line width=2.25]    (273,161.94) -- (228.5,195) ;

\draw [color={rgb, 255:red, 208; green, 2; blue, 27 }  ,draw opacity=1 ][line width=2.25]    (273,161.94) -- (362,161.94) ;

\draw [color={rgb, 255:red, 189; green, 16; blue, 224 }  ,draw opacity=1 ][line width=2.25]    (229.41,99.61) -- (228.5,195) ;

\draw [color={rgb, 255:red, 74; green, 144; blue, 226 }  ,draw opacity=1 ][line width=2.25]    (143.33,61.07) -- (143.33,126.35) ;

\draw [color={rgb, 255:red, 74; green, 144; blue, 226 }  ,draw opacity=1 ][line width=2.25]    (69.18,61.07) -- (143.33,61.07) ;

\draw [color={rgb, 255:red, 144; green, 19; blue, 254 }  ,draw opacity=1 ][line width=2.25]    (68.5,127) -- (143.33,126.35) ;

\draw [color={rgb, 255:red, 74; green, 144; blue, 226 }  ,draw opacity=1 ][line width=2.25]    (101.83,33) -- (69.18,61.07) ;

\draw [color={rgb, 255:red, 144; green, 19; blue, 254 }  ,draw opacity=1 ][line width=2.25]    (143.33,61.07) -- (168.5,33) ;

\draw [color={rgb, 255:red, 144; green, 19; blue, 254 }  ,draw opacity=1 ][line width=2.25]    (101.83,33) -- (168.5,33) ;

\draw [color={rgb, 255:red, 74; green, 144; blue, 226 }  ,draw opacity=1 ][line width=2.25]    (168.5,33) -- (168.5,104.15) ;

\draw [color={rgb, 255:red, 144; green, 19; blue, 254 }  ,draw opacity=1 ][line width=2.25]    (168.5,104.15) -- (143.33,126.35) ;

\draw [color={rgb, 255:red, 144; green, 19; blue, 254 }  ,draw opacity=1 ][line width=2.25]    (101.83,33) -- (101.83,104.15) ;

\draw [color={rgb, 255:red, 74; green, 144; blue, 226 }  ,draw opacity=1 ][line width=2.25]    (101.83,104.15) -- (68.5,127) ;

\draw [color={rgb, 255:red, 74; green, 144; blue, 226 }  ,draw opacity=1 ][line width=2.25]    (101.83,104.15) -- (168.5,104.15) ;

\draw [color={rgb, 255:red, 144; green, 19; blue, 254 }  ,draw opacity=1 ][line width=2.25]    (69.18,61.07) -- (68.5,127) ;

\draw [color={rgb, 255:red, 65; green, 117; blue, 5 }  ,draw opacity=1 ][line width=2.25]    (496.83,59.97) -- (496.83,127.33) ;

\draw [color={rgb, 255:red, 245; green, 166; blue, 35 }  ,draw opacity=1 ][line width=2.25]    (420.5,128) -- (496.83,127.33) ;

\draw [color={rgb, 255:red, 65; green, 117; blue, 5 }  ,draw opacity=1 ][line width=2.25]    (454.5,31) -- (421.19,59.97) ;

\draw [color={rgb, 255:red, 245; green, 166; blue, 35 }  ,draw opacity=1 ][line width=2.25]    (496.83,59.97) -- (522.5,31) ;

\draw [color={rgb, 255:red, 245; green, 166; blue, 35 }  ,draw opacity=1 ][line width=2.25]    (454.5,31) -- (522.5,31) ;

\draw [color={rgb, 255:red, 65; green, 117; blue, 5 }  ,draw opacity=1 ][line width=2.25]    (522.5,31) -- (522.5,104.42) ;

\draw [color={rgb, 255:red, 245; green, 166; blue, 35 }  ,draw opacity=1 ][line width=2.25]    (522.5,104.42) -- (496.83,127.33) ;

\draw [color={rgb, 255:red, 245; green, 166; blue, 35 }  ,draw opacity=1 ][line width=2.25]    (454.5,31) -- (454.5,104.42) ;

\draw [color={rgb, 255:red, 65; green, 117; blue, 5 }  ,draw opacity=1 ][line width=2.25]    (454.5,104.42) -- (420.5,128) ;

\draw [color={rgb, 255:red, 245; green, 166; blue, 35 }  ,draw opacity=1 ][line width=2.25]    (421.19,59.97) -- (420.5,128) ;

\draw [color={rgb, 255:red, 65; green, 117; blue, 5 }  ,draw opacity=1 ][line width=2.25]    (421.19,59.97) -- (496.83,59.97) ;

\draw [color={rgb, 255:red, 144; green, 19; blue, 254 }  ,draw opacity=1 ][line width=2.25]    (137.58,166.37) -- (137.58,232.34) ;

\draw [color={rgb, 255:red, 74; green, 144; blue, 226 }  ,draw opacity=1 ][line width=2.25]    (64.17,166.37) -- (137.58,166.37) ;

\draw [color={rgb, 255:red, 144; green, 19; blue, 254 }  ,draw opacity=1 ][line width=2.25]    (63.5,233) -- (137.58,232.34) ;

\draw [color={rgb, 255:red, 144; green, 19; blue, 254 }  ,draw opacity=1 ][line width=2.25]    (96.5,138) -- (64.17,166.37) ;

\draw [color={rgb, 255:red, 74; green, 144; blue, 226 }  ,draw opacity=1 ][line width=2.25]    (137.58,166.37) -- (162.5,138) ;

\draw [color={rgb, 255:red, 144; green, 19; blue, 254 }  ,draw opacity=1 ][line width=2.25]    (96.5,138) -- (162.5,138) ;

\draw [color={rgb, 255:red, 144; green, 19; blue, 254 }  ,draw opacity=1 ][line width=2.25]    (162.5,138) -- (162.5,209.91) ;

\draw [color={rgb, 255:red, 74; green, 144; blue, 226 }  ,draw opacity=1 ][line width=2.25]    (162.5,209.91) -- (137.58,232.34) ;

\draw [color={rgb, 255:red, 74; green, 144; blue, 226 }  ,draw opacity=1 ][line width=2.25]    (96.5,138) -- (96.5,209.91) ;

\draw [color={rgb, 255:red, 144; green, 19; blue, 254 }  ,draw opacity=1 ][line width=2.25]    (96.5,209.91) -- (63.5,233) ;

\draw [color={rgb, 255:red, 74; green, 144; blue, 226 }  ,draw opacity=1 ][line width=2.25]    (96.5,209.91) -- (162.5,209.91) ;

\draw [color={rgb, 255:red, 74; green, 144; blue, 226 }  ,draw opacity=1 ][line width=2.25]    (64.17,166.37) -- (63.5,233) ;

\draw [color={rgb, 255:red, 245; green, 166; blue, 35 }  ,draw opacity=1 ][line width=2.25]    (450.83,213.42) -- (517.5,213.42) ;

\draw [color={rgb, 255:red, 245; green, 166; blue, 35 }  ,draw opacity=1 ][line width=2.25]    (492.33,168.97) -- (492.33,236.33) ;

\draw [color={rgb, 255:red, 245; green, 166; blue, 35 }  ,draw opacity=1 ][line width=2.25]    (417.5,237) -- (492.33,236.33) ;

\draw [color={rgb, 255:red, 245; green, 166; blue, 35 }  ,draw opacity=1 ][line width=2.25]    (450.83,140) -- (418.18,168.97) ;

\draw [color={rgb, 255:red, 65; green, 117; blue, 5 }  ,draw opacity=1 ][line width=2.25]    (492.33,168.97) -- (517.5,140) ;

\draw [color={rgb, 255:red, 245; green, 166; blue, 35 }  ,draw opacity=1 ][line width=2.25]    (450.83,140) -- (517.5,140) ;

\draw [color={rgb, 255:red, 245; green, 166; blue, 35 }  ,draw opacity=1 ][line width=2.25]    (517.5,140) -- (517.5,213.42) ;

\draw [color={rgb, 255:red, 65; green, 117; blue, 5 }  ,draw opacity=1 ][line width=2.25]    (517.5,213.42) -- (492.33,236.33) ;

\draw [color={rgb, 255:red, 65; green, 117; blue, 5 }  ,draw opacity=1 ][line width=2.25]    (450.83,140) -- (450.83,213.42) ;

\draw [color={rgb, 255:red, 245; green, 166; blue, 35 }  ,draw opacity=1 ][line width=2.25]    (450.83,213.42) -- (417.5,237) ;

\draw [color={rgb, 255:red, 65; green, 117; blue, 5 }  ,draw opacity=1 ][line width=2.25]    (418.18,168.97) -- (417.5,237) ;

\draw [color={rgb, 255:red, 65; green, 117; blue, 5 }  ,draw opacity=1 ][line width=2.25]    (418.18,168.97) -- (492.33,168.97) ;

\draw (114,267) node  [align=left] {FD (1)};
\draw (289,220) node  [align=left] {FD (2)};
\draw (466,264) node  [align=left] {FD (3)};

\end{tikzpicture}

    \caption{Possible fundamental domains on the cube.}
   
    \label{fig:my_label}
\end{figure}
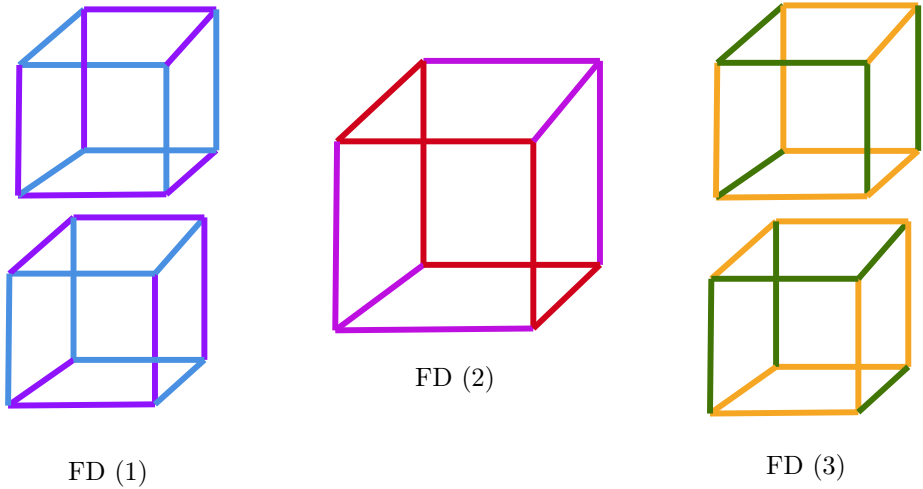

Our goal in this section is to classify all possible fundamental domains on the cube with torsion-free groups. We will first prove that the only candidates to be fundamental domains are the three shown in Figure 3. By candidates, we mean that there are 2 edge classes, a set of face identifications producing them, and both the Rivin and Fundamental domain equations are satisfied. The last aspect to check will be that all of the relations in the group are equal to the identity and not a loxodromic element. It is possible for all of the above to be true and have a relation that is a loxodromic element with axis fixing the original edge in the relator. In this case, we will not have a fundamental domain. Note that there are two polyhedrons associated with FD (1) and FD (3). The two versions are not exactly the same but are analogous in a way that will be made precise when we calculate the associated group in Section 5.2.

The full descriptions are as follows.
\begin{itemize}
    \item FD (1) (top version)

    \subitem A: Top$\longrightarrow$Bottom with $\dfrac{\pi}{2}$ clockwise twist on bottom interior face
    \subitem B: Front$\longrightarrow$Back with $\dfrac{\pi}{2}$ clockwise twist on back interior face
    \subitem C: Left$\longrightarrow$Right with $\dfrac{\pi}{2}$ clockwise twist on right interior face

All of the exterior dihedral angles are $\dfrac{2\pi}{3}$. It is straightforward to check that these angles satisfy the Rivin equations and inequalities, as well as the fundamental domain equations: each vertex has three incident edges, each with angle equal to $\dfrac{2\pi}{3}$, and $3\cdot\dfrac{2\pi}{3}=2\pi$. Any closed curve not bounding a face in the dual (an octahedron) has at least four edges, and $4\cdot\dfrac{2\pi}{3}>2\pi$, as desired. Finally, each edge class has six edges, and the sum of any six interior dihedral angles is $6\cdot\dfrac{\pi}{3}=2\pi$.
    \item FD (2)
   
    Face Identifications:
    \subitem Top$\longrightarrow$Left with $\dfrac{\pi}{2}$ clockwise twist on left interior face
    \subitem Right$\longrightarrow$Bottom with $\dfrac{\pi}{2}$ counterclockwise twist on bottom interior face
    \subitem Front$\longrightarrow$Back with $\pi$ twist

All external dihedral angles are $\dfrac{2\pi}{3}$; note that as with the first example, this satisfies all of Rivin's equations and inequalities, and it satisfies the fundamental domain equations.

\item FD (3) (top version)

Face Identifications:
\subitem Top$\longrightarrow$Front with $\dfrac{\pi}{2}$ clockwise twist on front exterior face
\subitem Left$\longrightarrow$Right with $\dfrac{\pi}{2}$ counterclockwise twist on right interior face
\subitem Back$\longrightarrow$Bottom with $\pi$ twist on bottom interior face
\end{itemize}
\begin{figure}[H]
    \centering

\tikzset{every picture/.style={line width=0.75pt}} 

\begin{tikzpicture}[x=0.75pt,y=0.75pt,yscale=-1,xscale=1]

\draw [color={rgb, 255:red, 245; green, 166; blue, 35 }  ,draw opacity=1 ][line width=2.25]    (142.5,171) -- (240.5,171) ;

\draw [color={rgb, 255:red, 245; green, 166; blue, 35 }  ,draw opacity=1 ][line width=2.25]    (142.5,62) -- (142.5,171) ;

\draw [color={rgb, 255:red, 245; green, 166; blue, 35 }  ,draw opacity=1 ][line width=2.25]    (94.5,105) -- (93.5,206) ;

\draw [color={rgb, 255:red, 65; green, 117; blue, 5 }  ,draw opacity=1 ][line width=2.25]    (203.5,105) -- (203.5,205) ;

\draw [color={rgb, 255:red, 65; green, 117; blue, 5 }  ,draw opacity=1 ][line width=2.25]    (94.5,105) -- (203.5,105) ;

\draw [color={rgb, 255:red, 245; green, 166; blue, 35 }  ,draw opacity=1 ][line width=2.25]    (93.5,206) -- (203.5,205) ;

\draw [color={rgb, 255:red, 65; green, 117; blue, 5 }  ,draw opacity=1 ][line width=2.25]    (142.5,62) -- (94.5,105) ;

\draw [color={rgb, 255:red, 245; green, 166; blue, 35 }  ,draw opacity=1 ][line width=2.25]    (203.5,105) -- (240.5,62) ;

\draw [color={rgb, 255:red, 245; green, 166; blue, 35 }  ,draw opacity=1 ][line width=2.25]    (142.5,62) -- (240.5,62) ;

\draw [color={rgb, 255:red, 65; green, 117; blue, 5 }  ,draw opacity=1 ][line width=2.25]    (240.5,62) -- (240.5,171) ;

\draw [color={rgb, 255:red, 245; green, 166; blue, 35 }  ,draw opacity=1 ][line width=2.25]    (240.5,171) -- (203.5,205) ;

\draw [color={rgb, 255:red, 65; green, 117; blue, 5 }  ,draw opacity=1 ][line width=2.25]    (142.5,171) -- (93.5,206) ;

\draw (546.1,94.16) node   {$\textcolor[rgb]{0.25,0.46,0.02}{x}\textcolor[rgb]{0.25,0.46,0.02}{_{10}}\textcolor[rgb]{0.25,0.46,0.02}{=}\textcolor[rgb]{0.25,0.46,0.02}{\dfrac{3\pi }{5}}$};
\draw (458.26,145.46) node   {$\textcolor[rgb]{0.96,0.65,0.14}{x}\textcolor[rgb]{0.96,0.65,0.14}{_{7}}\textcolor[rgb]{0.96,0.65,0.14}{=}\textcolor[rgb]{0.96,0.65,0.14}{\dfrac{4\pi }{5}}$};
\draw (546.1,203.46) node   {$\textcolor[rgb]{0.96,0.65,0.14}{x}\textcolor[rgb]{0.96,0.65,0.14}{_{12}}\textcolor[rgb]{0.96,0.65,0.14}{=}\textcolor[rgb]{0.96,0.65,0.14}{\dfrac{4\pi }{5}}$};
\draw (546.1,43.97) node   {$\textcolor[rgb]{0.96,0.65,0.14}{x}\textcolor[rgb]{0.96,0.65,0.14}{_{9}}\textcolor[rgb]{0.96,0.65,0.14}{=}\textcolor[rgb]{0.96,0.65,0.14}{\dfrac{4\pi }{5}}$};
\draw (457.26,47.31) node   {$\textcolor[rgb]{0.25,0.46,0.02}{x}\textcolor[rgb]{0.25,0.46,0.02}{_{5}}\textcolor[rgb]{0.25,0.46,0.02}{=}\textcolor[rgb]{0.25,0.46,0.02}{\dfrac{3\pi }{5}}$};
\draw (371.42,204.58) node   {$\textcolor[rgb]{0.96,0.65,0.14}{x}\textcolor[rgb]{0.96,0.65,0.14}{_{4}}\textcolor[rgb]{0.96,0.65,0.14}{=}\textcolor[rgb]{0.96,0.65,0.14}{\dfrac{3\pi }{5}}$};
\draw (372.42,147.69) node   {$\textcolor[rgb]{0.96,0.65,0.14}{x}\textcolor[rgb]{0.96,0.65,0.14}{_{3}}\textcolor[rgb]{0.96,0.65,0.14}{=}\textcolor[rgb]{0.96,0.65,0.14}{\dfrac{4\pi }{5}}$};
\draw (371.42,99.73) node   {$\textcolor[rgb]{0.25,0.46,0.02}{x}\textcolor[rgb]{0.25,0.46,0.02}{_{2}}\textcolor[rgb]{0.25,0.46,0.02}{=}\textcolor[rgb]{0.25,0.46,0.02}{\dfrac{3\pi }{5}}$};
\draw (545.1,144.35) node   {$\textcolor[rgb]{0.25,0.46,0.02}{x}\textcolor[rgb]{0.25,0.46,0.02}{_{11}}\textcolor[rgb]{0.25,0.46,0.02}{=}\textcolor[rgb]{0.25,0.46,0.02}{\dfrac{3\pi }{5}}$};
\draw (460.26,203.46) node   {$\textcolor[rgb]{0.96,0.65,0.14}{x}\textcolor[rgb]{0.96,0.65,0.14}{_{8}}\textcolor[rgb]{0.96,0.65,0.14}{=}\textcolor[rgb]{0.96,0.65,0.14}{\dfrac{3\pi }{5}}$};
\draw (457.26,97.5) node   {$\textcolor[rgb]{0.25,0.46,0.02}{x}\textcolor[rgb]{0.25,0.46,0.02}{_{6}}\textcolor[rgb]{0.25,0.46,0.02}{=}\textcolor[rgb]{0.25,0.46,0.02}{\dfrac{3\pi }{5}}$};
\draw (369.43,48.43) node   {$\textcolor[rgb]{0.96,0.65,0.14}{x}\textcolor[rgb]{0.96,0.65,0.14}{_{1}}\textcolor[rgb]{0.96,0.65,0.14}{=}\textcolor[rgb]{0.96,0.65,0.14}{\dfrac{3\pi }{5}}$};
\draw (116,171) node  [align=left] {\textcolor[rgb]{0.25,0.46,0.02}{2}};
\draw (151,217) node  [align=left] {\textcolor[rgb]{0.96,0.65,0.14}{1}};
\draw (234,193) node  [align=left] {\textcolor[rgb]{0.96,0.65,0.14}{3}};
\draw (183,161) node  [align=left] {\textcolor[rgb]{0.96,0.65,0.14}{4}};
\draw (165,117) node  [align=left] {\textcolor[rgb]{0.25,0.46,0.02}{5}};
\draw (107,80) node  [align=left] {\textcolor[rgb]{0.25,0.46,0.02}{6}};
\draw (227,96) node  [align=left] {\textcolor[rgb]{0.96,0.65,0.14}{7}};
\draw (191,50) node  [align=left] {\textcolor[rgb]{0.96,0.65,0.14}{8}};
\draw (82,147) node  [align=left] {\textcolor[rgb]{0.96,0.65,0.14}{9}};
\draw (216,145) node  [align=left] {\textcolor[rgb]{0.25,0.46,0.02}{10}};
\draw (254,125) node  [align=left] {\textcolor[rgb]{0.25,0.46,0.02}{11}};
\draw (129,129) node  [align=left] {\textcolor[rgb]{0.96,0.65,0.14}{12}};
\draw (161,251) node  [align=left] {(a)};
\draw (462,251) node  [align=left] {(b)};

\end{tikzpicture}

    \caption{(a) Edge labeling for FD (3) (b) Exterior dihedral angle assignments}
    \label{fig:my_label}
\end{figure}
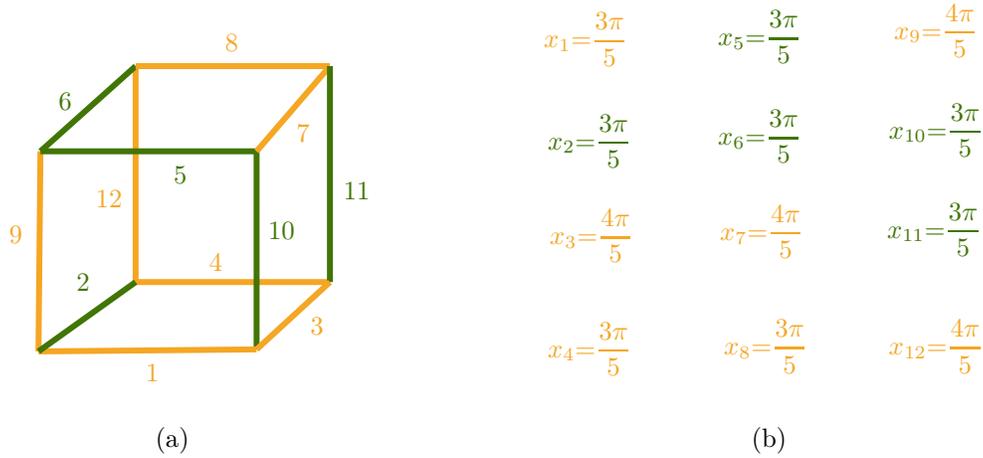

The angle assignments in Figure 5(b) satisfy the Rivin equations; each vertex has two angles of $\dfrac{3\pi}{5}$ and one angle of $\dfrac{4\pi}{5}$, which sum to $2\pi$. Any closed cycle not bounding a face in the dual has to consist of at least 4 edges and has a minimum length of $\dfrac{12\pi}{5}>2\pi$. Finally, it is straightforward to check that the sum of each edge class is what it should be: $3\pi$ for the green class and $5\pi$ for the yellow class.

\textbf{Proposition 5.1}: FDs (1),(2), and (3) are the only possible candidates for fundamental domains on the cube with torsion-free groups, up to rotation.

First note that by Proposition 4.1, we know that a fundamental domain on the cube without elliptic generators must have 2 edge classes. Thus, the first goal is to figure out what size edge classes are allowed. For notation purposes, we will denote a given edge breakdown by $x$-$y$, where $x+y=12$.

\textbf{Claim 1}: The only allowed edge breakdowns for the cube with non-elliptic generators are 5-7 and 6-6.
\begin{proof}
Certainly we cannot have an edge class of less than or equal to 2 elements: for if we did, then by the fundamental domain equations we have \begin{equation}
    \pi-x_i+\pi-x_j=2\pi\  \longrightarrow \ x_i+x_j=0,
\end{equation} which is not possible. Now consider the case that there is an edge class of size 3. Because 3 is odd, we will need to have adjacent faces identified and sharing an edge, so we consider all of those cases, shown in Figure 6. WLOG we will consider the front and left faces the adjacent faces. Either the edges on the face are opposite each other or adjacent. If they are opposite, there is only one case, as shown in Figure 6(i). If they are adjacent, then there are in principle 4 cases corresponding to the twists by 0,$\frac{\pi}{2}$, $\pi$, and $\frac{3\pi}{2}$; however, $\pi$ and $\frac{3\pi}{2}$ result in the edge class having four edges. The remaining two cases are shown in Figure 6(ii,iii).

\begin{figure}[H]
    \centering

\tikzset{every picture/.style={line width=0.75pt}} 

\begin{tikzpicture}[x=0.6pt,y=0.6pt,yscale=-1,xscale=1]

\draw [color={rgb, 255:red, 74; green, 144; blue, 226 }  ,draw opacity=1 ][line width=2.25]    (205.72,108.16) -- (205.72,187.59) ;

\draw [color={rgb, 255:red, 0; green, 0; blue, 0 }  ,draw opacity=1 ][line width=2.25]  [dash pattern={on 2.53pt off 3.02pt}]  (123.87,108.16) -- (205.72,108.16) ;

\draw [color={rgb, 255:red, 0; green, 0; blue, 0 }  ,draw opacity=1 ][line width=2.25]  [dash pattern={on 2.53pt off 3.02pt}]  (123.12,188.38) -- (205.72,187.59) ;

\draw [color={rgb, 255:red, 0; green, 0; blue, 0 }  ,draw opacity=1 ][line width=2.25]  [dash pattern={on 2.53pt off 3.02pt}]  (159.91,74) -- (123.87,108.16) ;

\draw [color={rgb, 255:red, 0; green, 0; blue, 0 }  ,draw opacity=1 ][line width=2.25]  [dash pattern={on 2.53pt off 3.02pt}]  (205.72,108.16) -- (233.5,74) ;

\draw [color={rgb, 255:red, 0; green, 0; blue, 0 }  ,draw opacity=1 ][line width=2.25]  [dash pattern={on 2.53pt off 3.02pt}]  (159.91,74) -- (233.5,74) ;

\draw [color={rgb, 255:red, 0; green, 0; blue, 0 }  ,draw opacity=1 ][line width=2.25]  [dash pattern={on 2.53pt off 3.02pt}]  (233.5,74) -- (233.5,160.58) ;

\draw [color={rgb, 255:red, 0; green, 0; blue, 0 }  ,draw opacity=1 ][line width=2.25]  [dash pattern={on 2.53pt off 3.02pt}]  (233.5,160.58) -- (205.72,187.59) ;

\draw [color={rgb, 255:red, 74; green, 144; blue, 226 }  ,draw opacity=1 ][line width=2.25]    (159.91,74) -- (159.91,160.58) ;

\draw [color={rgb, 255:red, 0; green, 0; blue, 0 }  ,draw opacity=1 ][line width=2.25]  [dash pattern={on 2.53pt off 3.02pt}]  (159.91,160.58) -- (123.12,188.38) ;

\draw [color={rgb, 255:red, 0; green, 0; blue, 0 }  ,draw opacity=1 ][line width=2.25]  [dash pattern={on 2.53pt off 3.02pt}]  (159.91,160.58) -- (233.5,160.58) ;

\draw [color={rgb, 255:red, 74; green, 144; blue, 226 }  ,draw opacity=1 ][line width=2.25]    (123.87,108.16) -- (123.12,188.38) ;

\draw [rotate around= { 270: (123.49, 148.27)
    }] [color={rgb, 255:red, 208; green, 2; blue, 27 }  ,draw opacity=1 ][line width=1.5]   (123.49, 148.27) circle [x radius= 47.73, y radius= 16.82]  ;
\draw [color={rgb, 255:red, 0; green, 0; blue, 0 }  ,draw opacity=1 ][line width=2.25]  [dash pattern={on 2.53pt off 3.02pt}]  (485.59,108.56) -- (485.59,183.25) ;

\draw [color={rgb, 255:red, 74; green, 144; blue, 226 }  ,draw opacity=1 ][line width=2.25]    (400.28,108.56) -- (485.59,108.56) ;

\draw [color={rgb, 255:red, 0; green, 0; blue, 0 }  ,draw opacity=1 ][line width=2.25]  [dash pattern={on 2.53pt off 3.02pt}]  (399.5,184) -- (485.59,183.25) ;

\draw [color={rgb, 255:red, 208; green, 2; blue, 27 }  ,draw opacity=1 ][line width=2.25]  [dash pattern={on 2.53pt off 3.02pt}]  (437.85,76.44) -- (400.28,108.56) ;

\draw [color={rgb, 255:red, 0; green, 0; blue, 0 }  ,draw opacity=1 ][line width=2.25]  [dash pattern={on 2.53pt off 3.02pt}]  (485.59,108.56) -- (514.55,76.44) ;

\draw [color={rgb, 255:red, 208; green, 2; blue, 27 }  ,draw opacity=1 ][line width=2.25]  [dash pattern={on 2.53pt off 3.02pt}]  (437.85,76.44) -- (514.55,76.44) ;

\draw [color={rgb, 255:red, 208; green, 2; blue, 27 }  ,draw opacity=1 ][line width=2.25]  [dash pattern={on 2.53pt off 3.02pt}]  (514.55,76.44) -- (514.55,157.86) ;

\draw [color={rgb, 255:red, 208; green, 2; blue, 27 }  ,draw opacity=1 ][line width=2.25]  [dash pattern={on 2.53pt off 3.02pt}]  (514.55,157.86) -- (485.59,183.25) ;

\draw [color={rgb, 255:red, 208; green, 2; blue, 27 }  ,draw opacity=1 ][line width=2.25]  [dash pattern={on 2.53pt off 3.02pt}]  (437.85,76.44) -- (437.85,157.86) ;

\draw [color={rgb, 255:red, 74; green, 144; blue, 226 }  ,draw opacity=1 ][line width=2.25]    (437.85,157.86) -- (399.5,184) ;

\draw [color={rgb, 255:red, 208; green, 2; blue, 27 }  ,draw opacity=1 ][line width=2.25]  [dash pattern={on 2.53pt off 3.02pt}]  (437.85,157.86) -- (514.55,157.86) ;

\draw [color={rgb, 255:red, 74; green, 144; blue, 226 }  ,draw opacity=1 ][line width=2.25]    (400.28,108.56) -- (399.5,184) ;

\draw [color={rgb, 255:red, 208; green, 2; blue, 27 }  ,draw opacity=1 ][line width=2.25]  [dash pattern={on 2.53pt off 3.02pt}]  (341.26,107.29) -- (341.26,186.59) ;

\draw [color={rgb, 255:red, 74; green, 144; blue, 226 }  ,draw opacity=1 ][line width=2.25]    (259.25,107.29) -- (341.26,107.29) ;

\draw [color={rgb, 255:red, 208; green, 2; blue, 27 }  ,draw opacity=1 ][line width=2.25]  [dash pattern={on 2.53pt off 3.02pt}]  (258.5,187.38) -- (341.26,186.59) ;

\draw [color={rgb, 255:red, 74; green, 144; blue, 226 }  ,draw opacity=1 ][line width=2.25]    (295.37,73.19) -- (259.25,107.29) ;

\draw [color={rgb, 255:red, 208; green, 2; blue, 27 }  ,draw opacity=1 ][line width=2.25]  [dash pattern={on 2.53pt off 3.02pt}]  (341.26,107.29) -- (369.1,73.19) ;

\draw [color={rgb, 255:red, 208; green, 2; blue, 27 }  ,draw opacity=1 ][line width=2.25]  [dash pattern={on 2.53pt off 3.02pt}]  (295.37,73.19) -- (369.1,73.19) ;

\draw [color={rgb, 255:red, 208; green, 2; blue, 27 }  ,draw opacity=1 ][line width=2.25]  [dash pattern={on 2.53pt off 3.02pt}]  (369.1,73.19) -- (369.1,159.63) ;

\draw [color={rgb, 255:red, 208; green, 2; blue, 27 }  ,draw opacity=1 ][line width=2.25]  [dash pattern={on 2.53pt off 3.02pt}]  (369.1,159.63) -- (341.26,186.59) ;

\draw [color={rgb, 255:red, 0; green, 0; blue, 0 }  ,draw opacity=1 ][line width=2.25]  [dash pattern={on 2.53pt off 3.02pt}]  (295.37,73.19) -- (295.37,159.63) ;

\draw [color={rgb, 255:red, 0; green, 0; blue, 0 }  ,draw opacity=1 ][line width=2.25]  [dash pattern={on 2.53pt off 3.02pt}]  (295.37,159.63) -- (258.5,187.38) ;

\draw [color={rgb, 255:red, 0; green, 0; blue, 0 }  ,draw opacity=1 ][line width=2.25]  [dash pattern={on 2.53pt off 3.02pt}]  (295.37,159.63) -- (369.1,159.63) ;

\draw [color={rgb, 255:red, 74; green, 144; blue, 226 }  ,draw opacity=1 ][line width=2.25]    (259.25,107.29) -- (258.5,187.38) ;

\draw  [color={rgb, 255:red, 208; green, 2; blue, 27 }  ,draw opacity=1 ][line width=1.5]  (520.96,165.64) .. controls (516.64,169.73) and (509.69,169.4) .. (505.44,164.91) .. controls (501.18,160.41) and (501.24,153.46) .. (505.56,149.37) .. controls (509.88,145.28) and (516.83,145.61) .. (521.08,150.1) .. controls (525.34,154.59) and (525.28,161.55) .. (520.96,165.64) -- cycle ;
\draw  [color={rgb, 255:red, 208; green, 2; blue, 27 }  ,draw opacity=1 ][line width=1.5]  (445.98,84.97) .. controls (441.66,89.06) and (434.71,88.73) .. (430.46,84.24) .. controls (426.2,79.74) and (426.26,72.78) .. (430.58,68.69) .. controls (434.9,64.6) and (441.85,64.93) .. (446.11,69.43) .. controls (450.36,73.92) and (450.3,80.88) .. (445.98,84.97) -- cycle ;
\draw  [color={rgb, 255:red, 208; green, 2; blue, 27 }  ,draw opacity=1 ][line width=1.5]  (348.85,193.21) .. controls (344.64,197.19) and (337.87,196.87) .. (333.72,192.49) .. controls (329.58,188.12) and (329.63,181.34) .. (333.84,177.35) .. controls (338.05,173.37) and (344.83,173.69) .. (348.97,178.07) .. controls (353.11,182.45) and (353.06,189.23) .. (348.85,193.21) -- cycle ;
\draw  [color={rgb, 255:red, 208; green, 2; blue, 27 }  ,draw opacity=1 ][line width=1.5]  (376.19,82.38) .. controls (371.98,86.36) and (365.21,86.04) .. (361.06,81.67) .. controls (356.92,77.29) and (356.97,70.51) .. (361.18,66.52) .. controls (365.39,62.54) and (372.16,62.86) .. (376.31,67.24) .. controls (380.45,71.62) and (380.4,78.4) .. (376.19,82.38) -- cycle ;
\draw    (358.5,217) -- (325.21,197.03) ;
\draw [shift={(323.5,196)}, rotate = 390.96000000000004] [color={rgb, 255:red, 0; green, 0; blue, 0 }  ][line width=0.75]    (10.93,-3.29) .. controls (6.95,-1.4) and (3.31,-0.3) .. (0,0) .. controls (3.31,0.3) and (6.95,1.4) .. (10.93,3.29)   ;

\draw    (358.5,217) -- (392.81,195.08) ;
\draw [shift={(394.5,194)}, rotate = 507.43] [color={rgb, 255:red, 0; green, 0; blue, 0 }  ][line width=0.75]    (10.93,-3.29) .. controls (6.95,-1.4) and (3.31,-0.3) .. (0,0) .. controls (3.31,0.3) and (6.95,1.4) .. (10.93,3.29)   ;

\draw (173,218) node  [align=left] {Elliptic Generator};
\draw (383,226) node  [align=left] {2 opposite vertices};
\draw (181,42) node  [align=left] {(i) Opposite\\ \ \ \ \ \ edges};
\draw (328,41) node  [align=left] {(ii) Adjacent edges, \\ \ \ \ \ \ \ \ \ \ \ no twist};
\draw (490,43) node  [align=left] {{\small (iii) Adjacent edges, }\\{\small  \ \ \ \ \ }{\scriptsize  $\displaystyle \frac{3\pi }{2}$}{\small  twist CW}};

\end{tikzpicture}

    \caption{Ways to have an edge class of size 3, as described above. Each poses a problem.}
    \label{fig:my_label}
    \end{figure}
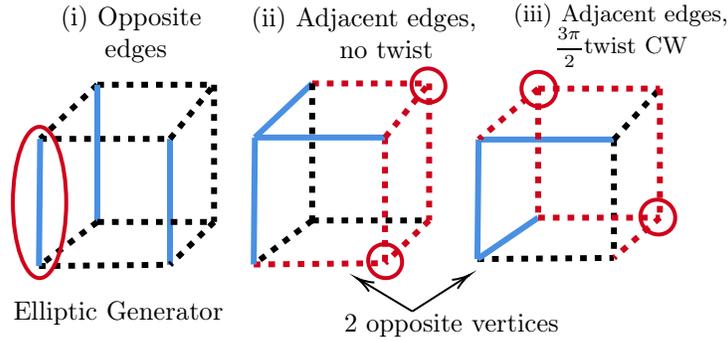

We are not considering the case with elliptic generators, so case (i) can be discarded. For cases (ii) and (iii) where we have 2 opposite vertices with all edges in one class, we know that six of the elements have to sum to $4\pi$ (each set of three edges at a vertex sums to $2\pi$ by the Rivin equations). The other class has to have 9 elements, so we have that \begin{equation}
    x_1+...+x_9=7\pi \ \mbox{and}\ x_1+...+x_6=4\pi.
\end{equation} Thus, we arrive at the fact that
\begin{equation}
    x_7+x_8+x_9=7\pi-4\pi=3\pi,
\end{equation}so at least one of the angles must be greater than or equal to $\pi$, a contradiction to convexity. The last case to discard (not pictured) is that there are three separate edges, each shared between two faces. However, it is clear that all three face pairs would be identified by an elliptic fixing the shared edge, so we cannot have this case. These are all possible ways to put three edges on the cube, up to rotation, and thus we conclude that we cannot have an edge class of three elements.

Suppose we had an edge class of four elements. Note that no face can have three edges, nor can any face have one edge, as that face will not have a face identification. This leaves us with two possibilities: we can have two edges on each of two non-adjacent faces, or we can have two adjacent faces sharing an edge with the fourth edge elsewhere on the cube. However, when we try to place the fourth edge, the only possible choices lead to an elliptic generator. We can thus consider only the case with two edges on opposite faces.  There are only 4 ways to achieve this, as shown in Figure 7: we can have the edges on the faces either opposite (Figure 7 (i) and (ii)) or adjacent (Figure 7 (iii) and (iv)), and in each case there are two possible twists. For adjacent edges, there is a $\frac{\pi}{2}$ twist (iv) and a $\pi$ twist (iii). The $\frac{3\pi}{2}$ case is analogous to the $\frac{\pi}{2}$ case, and the $0$ case requires an elliptic generator.
\begin{figure}[H]
    \centering

\tikzset{every picture/.style={line width=0.75pt}} 

\begin{tikzpicture}[x=0.6pt,y=0.6pt,yscale=-1,xscale=1]

\draw [color={rgb, 255:red, 74; green, 144; blue, 226 }  ,draw opacity=1 ][line width=2.25]    (260.32,67.56) -- (260.32,136.31) ;

\draw [color={rgb, 255:red, 0; green, 0; blue, 0 }  ,draw opacity=1 ][line width=2.25]  [dash pattern={on 2.53pt off 3.02pt}]  (183.21,67.56) -- (260.32,67.56) ;

\draw [color={rgb, 255:red, 0; green, 0; blue, 0 }  ,draw opacity=1 ][line width=2.25]  [dash pattern={on 2.53pt off 3.02pt}]  (182.5,137) -- (260.32,136.31) ;

\draw [color={rgb, 255:red, 0; green, 0; blue, 0 }  ,draw opacity=1 ][line width=2.25]  [dash pattern={on 2.53pt off 3.02pt}]  (217.17,38) -- (183.21,67.56) ;

\draw [color={rgb, 255:red, 0; green, 0; blue, 0 }  ,draw opacity=1 ][line width=2.25]  [dash pattern={on 2.53pt off 3.02pt}]  (260.32,67.56) -- (286.5,38) ;

\draw [color={rgb, 255:red, 0; green, 0; blue, 0 }  ,draw opacity=1 ][line width=2.25]  [dash pattern={on 2.53pt off 3.02pt}]  (217.17,38) -- (286.5,38) ;

\draw [color={rgb, 255:red, 74; green, 144; blue, 226 }  ,draw opacity=1 ][line width=2.25]    (286.5,38) -- (286.5,112.94) ;

\draw [color={rgb, 255:red, 0; green, 0; blue, 0 }  ,draw opacity=1 ][line width=2.25]  [dash pattern={on 2.53pt off 3.02pt}]  (286.5,112.94) -- (260.32,136.31) ;

\draw [color={rgb, 255:red, 74; green, 144; blue, 226 }  ,draw opacity=1 ][line width=2.25]    (217.17,38) -- (217.17,112.94) ;

\draw [color={rgb, 255:red, 0; green, 0; blue, 0 }  ,draw opacity=1 ][line width=2.25]  [dash pattern={on 2.53pt off 3.02pt}]  (217.17,112.94) -- (182.5,137) ;

\draw [color={rgb, 255:red, 0; green, 0; blue, 0 }  ,draw opacity=1 ][line width=2.25]  [dash pattern={on 2.53pt off 3.02pt}]  (217.17,112.94) -- (286.5,112.94) ;

\draw [color={rgb, 255:red, 74; green, 144; blue, 226 }  ,draw opacity=1 ][line width=2.25]    (183.21,67.56) -- (182.5,137) ;

\draw [color={rgb, 255:red, 0; green, 0; blue, 0 }  ,draw opacity=1 ][line width=2.25]  [dash pattern={on 2.53pt off 3.02pt}]  (415.32,196.56) -- (415.32,265.31) ;

\draw [color={rgb, 255:red, 0; green, 0; blue, 0 }  ,draw opacity=1 ][line width=2.25]  [dash pattern={on 2.53pt off 3.02pt}]  (338.21,196.56) -- (415.32,196.56) ;

\draw [color={rgb, 255:red, 74; green, 144; blue, 226 }  ,draw opacity=1 ][line width=2.25]    (337.5,266) -- (415.32,265.31) ;

\draw [color={rgb, 255:red, 0; green, 0; blue, 0 }  ,draw opacity=1 ][line width=2.25]  [dash pattern={on 2.53pt off 3.02pt}]  (372.17,167) -- (338.21,196.56) ;

\draw [color={rgb, 255:red, 74; green, 144; blue, 226 }  ,draw opacity=1 ][line width=2.25]    (415.32,196.56) -- (441.5,167) ;

\draw [color={rgb, 255:red, 74; green, 144; blue, 226 }  ,draw opacity=1 ][line width=2.25]    (372.17,167) -- (441.5,167) ;

\draw [color={rgb, 255:red, 0; green, 0; blue, 0 }  ,draw opacity=1 ][line width=2.25]  [dash pattern={on 2.53pt off 3.02pt}]  (441.5,167) -- (441.5,241.94) ;

\draw [color={rgb, 255:red, 74; green, 144; blue, 226 }  ,draw opacity=1 ][line width=2.25]    (441.5,241.94) -- (415.32,265.31) ;

\draw [color={rgb, 255:red, 0; green, 0; blue, 0 }  ,draw opacity=1 ][line width=2.25]  [dash pattern={on 2.53pt off 3.02pt}]  (372.17,167) -- (372.17,241.94) ;

\draw [color={rgb, 255:red, 0; green, 0; blue, 0 }  ,draw opacity=1 ][line width=2.25]  [dash pattern={on 2.53pt off 3.02pt}]  (372.17,241.94) -- (337.5,266) ;

\draw [color={rgb, 255:red, 0; green, 0; blue, 0 }  ,draw opacity=1 ][line width=2.25]  [dash pattern={on 2.53pt off 3.02pt}]  (372.17,241.94) -- (441.5,241.94) ;

\draw [color={rgb, 255:red, 0; green, 0; blue, 0 }  ,draw opacity=1 ][line width=2.25]  [dash pattern={on 2.53pt off 3.02pt}]  (338.21,196.56) -- (337.5,266) ;

\draw [color={rgb, 255:red, 0; green, 0; blue, 0 }  ,draw opacity=1 ][line width=2.25]  [dash pattern={on 2.53pt off 3.02pt}]  (252.32,197.56) -- (252.32,266.31) ;

\draw [color={rgb, 255:red, 0; green, 0; blue, 0 }  ,draw opacity=1 ][line width=2.25]  [dash pattern={on 2.53pt off 3.02pt}]  (175.21,197.56) -- (252.32,197.56) ;

\draw [color={rgb, 255:red, 74; green, 144; blue, 226 }  ,draw opacity=1 ][line width=2.25]    (174.5,267) -- (252.32,266.31) ;

\draw [color={rgb, 255:red, 0; green, 0; blue, 0 }  ,draw opacity=1 ][line width=2.25]  [dash pattern={on 2.53pt off 3.02pt}]  (209.17,168) -- (175.21,197.56) ;

\draw [color={rgb, 255:red, 74; green, 144; blue, 226 }  ,draw opacity=1 ][line width=2.25]    (252.32,197.56) -- (278.5,168) ;

\draw [color={rgb, 255:red, 74; green, 144; blue, 226 }  ,draw opacity=1 ][line width=2.25]    (209.17,168) -- (278.5,168) ;

\draw [color={rgb, 255:red, 0; green, 0; blue, 0 }  ,draw opacity=1 ][line width=2.25]  [dash pattern={on 2.53pt off 3.02pt}]  (278.5,168) -- (278.5,242.94) ;

\draw [color={rgb, 255:red, 0; green, 0; blue, 0 }  ,draw opacity=1 ][line width=2.25]  [dash pattern={on 2.53pt off 3.02pt}]  (278.5,242.94) -- (252.32,266.31) ;

\draw [color={rgb, 255:red, 0; green, 0; blue, 0 }  ,draw opacity=1 ][line width=2.25]  [dash pattern={on 2.53pt off 3.02pt}]  (209.17,168) -- (209.17,242.94) ;

\draw [color={rgb, 255:red, 74; green, 144; blue, 226 }  ,draw opacity=1 ][line width=2.25]    (209.17,242.94) -- (174.5,267) ;

\draw [color={rgb, 255:red, 0; green, 0; blue, 0 }  ,draw opacity=1 ][line width=2.25]  [dash pattern={on 2.53pt off 3.02pt}]  (209.17,242.94) -- (278.5,242.94) ;

\draw [color={rgb, 255:red, 0; green, 0; blue, 0 }  ,draw opacity=1 ][line width=2.25]  [dash pattern={on 2.53pt off 3.02pt}]  (175.21,197.56) -- (174.5,267) ;

\draw [color={rgb, 255:red, 74; green, 144; blue, 226 }  ,draw opacity=1 ][line width=2.25]    (421.32,64.56) -- (421.32,133.31) ;

\draw [color={rgb, 255:red, 0; green, 0; blue, 0 }  ,draw opacity=1 ][line width=2.25]  [dash pattern={on 2.53pt off 3.02pt}]  (344.21,64.56) -- (421.32,64.56) ;

\draw [color={rgb, 255:red, 0; green, 0; blue, 0 }  ,draw opacity=1 ][line width=2.25]  [dash pattern={on 2.53pt off 3.02pt}]  (343.5,134) -- (421.32,133.31) ;

\draw [color={rgb, 255:red, 0; green, 0; blue, 0 }  ,draw opacity=1 ][line width=2.25]  [dash pattern={on 2.53pt off 3.02pt}]  (378.17,35) -- (344.21,64.56) ;

\draw [color={rgb, 255:red, 0; green, 0; blue, 0 }  ,draw opacity=1 ][line width=2.25]  [dash pattern={on 2.53pt off 3.02pt}]  (421.32,64.56) -- (447.5,35) ;

\draw [color={rgb, 255:red, 74; green, 144; blue, 226 }  ,draw opacity=1 ][line width=2.25]    (378.17,35) -- (447.5,35) ;

\draw [color={rgb, 255:red, 0; green, 0; blue, 0 }  ,draw opacity=1 ][line width=2.25]  [dash pattern={on 2.53pt off 3.02pt}]  (447.5,35) -- (447.5,109.94) ;

\draw [color={rgb, 255:red, 0; green, 0; blue, 0 }  ,draw opacity=1 ][line width=2.25]  [dash pattern={on 2.53pt off 3.02pt}]  (447.5,109.94) -- (421.32,133.31) ;

\draw [color={rgb, 255:red, 0; green, 0; blue, 0 }  ,draw opacity=1 ][line width=2.25]  [dash pattern={on 2.53pt off 3.02pt}]  (378.17,35) -- (378.17,109.94) ;

\draw [color={rgb, 255:red, 0; green, 0; blue, 0 }  ,draw opacity=1 ][line width=2.25]  [dash pattern={on 2.53pt off 3.02pt}]  (378.17,109.94) -- (343.5,134) ;

\draw [color={rgb, 255:red, 74; green, 144; blue, 226 }  ,draw opacity=1 ][line width=2.25]    (378.17,109.94) -- (447.5,109.94) ;

\draw [color={rgb, 255:red, 74; green, 144; blue, 226 }  ,draw opacity=1 ][line width=2.25]    (344.21,64.56) -- (343.5,134) ;

\draw (144,52) node  [align=left] {(i)};
\draw (470,50) node  [align=left] {(ii)};
\draw (144,193) node  [align=left] {(iii)};
\draw (470,194) node  [align=left] {(iv)};

\end{tikzpicture}

    \caption{Possible ways to put 4 edges on the cube. Again, each poses a problem.}
    \label{fig:my_label}
\end{figure}
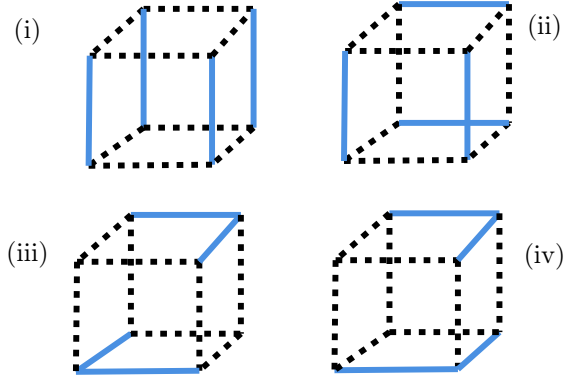

In case (i), the problem arises from one of Rivin's inequalities: the blue edge class must have exterior dihedral angles summing to $4\pi-2\pi=2\pi$, but in the dual, the blue edges form a path around the "waist" of the octahedron. Thus, we have a closed path not bounding a face of length equal to $2\pi$ in the dual, violating Rivin's inequalities. For (ii), the front face has to be identified with the back face because they are the only faces that are combinatorially equivalent. From there, we have a few choices: we can either identify the top face with the bottom, left, or right face. However, no matter which one we choose, we end up with either a $\pi$ twist or a $0$ twist for both of the other face identifications, and the result is that the edge class of 8 is actually broken into at least 2 smaller edge classes. Since it does not have two edge classes, we can discard this case.

For (iii), we have two opposite vertices in the same edge class, and \begin{equation} x_1+...+x_8=6\pi \ \mbox{and}\  x_1+...+x_6=2(2\pi)=4\pi.
\end{equation} Thus, we have that $x_7+x_8=2\pi$, and at least one of $x_7, x_8$ is greater than or equal to $\pi$, a contradiction.

Finally, for case (iv), the left face is all in one class, but there is not another face with four edges in that class, so we cannot assign face identifications.

Thus, we cannot have an edge class of four elements. Since we cannot have edge classes of sizes 1, 2, 3, or 4, we cannot have classes of size 11, 10, 9, or 8. Thus, we are left with the 6-6 and 5-7 cases.

\end{proof}

We have already shown that there are examples of edge classes with 5, 6, and 7 elements. Our goal now will be to show that if we have a fundamental domain with either a 6-6 breakdown or a 5-7 breakdown, it is one of the three cases shown in Figure 4.

\textbf{Claim 2}: FDs (1) and (2) are the only two 6-6 fundamental domains on the cube with non-elliptic generators.

\begin{proof}

We will break the proof into two parts by considering the following two cases: two edges per face and 3 edges per face (on some subset of the faces of the cube). At the end of the proof we will justify that these are indeed the only possible cases for the 6-6 breakdown.

\textbf{Case 1: 2 edges per face}
   
    To prove this claim, it suffices to pick two edges on one face and consider the possible rotations that send it to the opposite (and adjacent) face(s). There are two ways to put 2 edges on a face: the edges can either be opposite one another or adjacent to each other. We will then consider opposite and adjacent face identifications for each edge case, for a total of four possibilities. We will use the following color scheme: purple 'x's correspond to having $\geq$3 edges on a face. WLOG, we will have the opposite face identification pair front and back and the adjacent face identification pair top and front. For a different combination, we would simply rotate the picture.

\textbf{Adjacent edges, opposite faces}: We omit the case where there is no rotation between the front and back face because that case will require and elliptic generator over the shared edge between the top and back faces. The other three cases are shown in Figure 9.

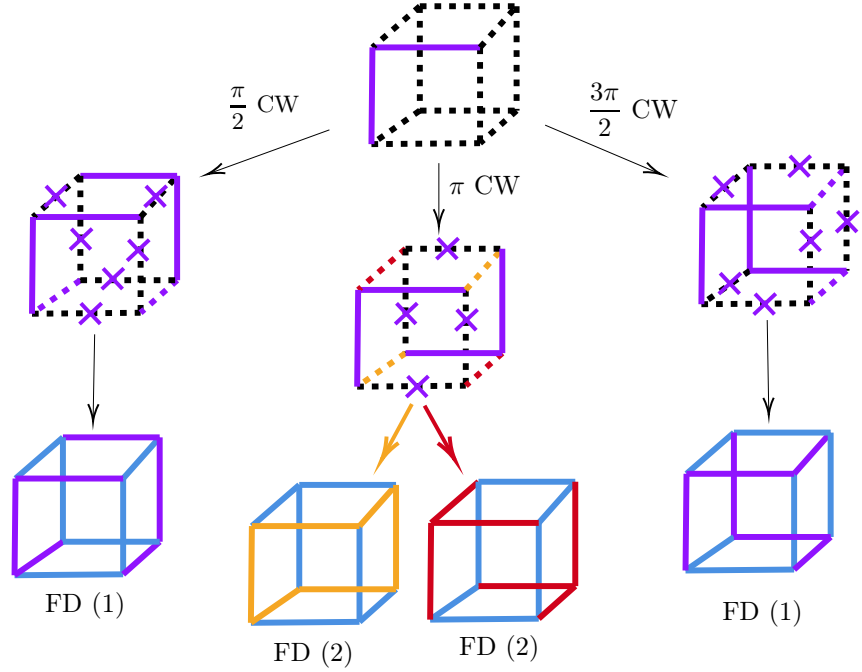
\begin{figure}[H]
    \centering
\tikzset{every picture/.style={line width=0.3pt}} 

\begin{tikzpicture}[x=0.6pt,y=0.6pt,yscale=-1,xscale=1]

\draw [color={rgb, 255:red, 0; green, 0; blue, 0 }  ,draw opacity=1 ][line width=2.25]  [dash pattern={on 2.53pt off 3.02pt}]  (348.97,46.98) -- (348.97,107.4) ;

\draw [color={rgb, 255:red, 144; green, 19; blue, 254 }  ,draw opacity=1 ][line width=2.25]    (281.12,46.98) -- (348.97,46.98) ;

\draw [color={rgb, 255:red, 0; green, 0; blue, 0 }  ,draw opacity=1 ][line width=2.25]  [dash pattern={on 2.53pt off 3.02pt}]  (280.5,108) -- (348.97,107.4) ;

\draw [color={rgb, 255:red, 0; green, 0; blue, 0 }  ,draw opacity=1 ][line width=2.25]  [dash pattern={on 2.53pt off 3.02pt}]  (311,21) -- (281.12,46.98) ;

\draw [color={rgb, 255:red, 0; green, 0; blue, 0 }  ,draw opacity=1 ][line width=2.25]  [dash pattern={on 2.53pt off 3.02pt}]  (348.97,46.98) -- (372,21) ;

\draw [color={rgb, 255:red, 0; green, 0; blue, 0 }  ,draw opacity=1 ][line width=2.25]  [dash pattern={on 2.53pt off 3.02pt}]  (311,21) -- (372,21) ;

\draw [color={rgb, 255:red, 0; green, 0; blue, 0 }  ,draw opacity=1 ][line width=2.25]  [dash pattern={on 2.53pt off 3.02pt}]  (372,21) -- (372,86.85) ;

\draw [color={rgb, 255:red, 0; green, 0; blue, 0 }  ,draw opacity=1 ][line width=2.25]  [dash pattern={on 2.53pt off 3.02pt}]  (372,86.85) -- (348.97,107.4) ;

\draw [color={rgb, 255:red, 0; green, 0; blue, 0 }  ,draw opacity=1 ][line width=2.25]  [dash pattern={on 2.53pt off 3.02pt}]  (311,21) -- (311,86.85) ;

\draw [color={rgb, 255:red, 0; green, 0; blue, 0 }  ,draw opacity=1 ][line width=2.25]  [dash pattern={on 2.53pt off 3.02pt}]  (311,86.85) -- (280.5,108) ;

\draw [color={rgb, 255:red, 0; green, 0; blue, 0 }  ,draw opacity=1 ][line width=2.25]  [dash pattern={on 2.53pt off 3.02pt}]  (311,86.85) -- (372,86.85) ;

\draw [color={rgb, 255:red, 144; green, 19; blue, 254 }  ,draw opacity=1 ][line width=2.25]    (281.12,46.98) -- (280.5,108) ;

\draw    (254,99) -- (178.88,126.32) ;
\draw [shift={(177,127)}, rotate = 340.02] [color={rgb, 255:red, 0; green, 0; blue, 0 }  ][line width=0.75]    (10.93,-3.29) .. controls (6.95,-1.4) and (3.31,-0.3) .. (0,0) .. controls (3.31,0.3) and (6.95,1.4) .. (10.93,3.29)   ;

\draw    (323,120) -- (323,158) ;
\draw [shift={(323,160)}, rotate = 270] [color={rgb, 255:red, 0; green, 0; blue, 0 }  ][line width=0.75]    (10.93,-3.29) .. controls (6.95,-1.4) and (3.31,-0.3) .. (0,0) .. controls (3.31,0.3) and (6.95,1.4) .. (10.93,3.29)   ;

\draw    (390,96) -- (461.16,126.22) ;
\draw [shift={(463,127)}, rotate = 203.01] [color={rgb, 255:red, 0; green, 0; blue, 0 }  ][line width=0.75]    (10.93,-3.29) .. controls (6.95,-1.4) and (3.31,-0.3) .. (0,0) .. controls (3.31,0.3) and (6.95,1.4) .. (10.93,3.29)   ;

\draw [color={rgb, 255:red, 0; green, 0; blue, 0 }  ,draw opacity=1 ][line width=2.25]  [dash pattern={on 2.53pt off 3.02pt}]  (134.97,153.98) -- (134.97,214.4) ;

\draw [color={rgb, 255:red, 144; green, 19; blue, 254 }  ,draw opacity=1 ][line width=2.25]    (67.12,153.98) -- (134.97,153.98) ;

\draw [color={rgb, 255:red, 0; green, 0; blue, 0 }  ,draw opacity=1 ][line width=2.25]  [dash pattern={on 2.53pt off 3.02pt}]  (66.5,215) -- (134.97,214.4) ;

\draw [color={rgb, 255:red, 0; green, 0; blue, 0 }  ,draw opacity=1 ][line width=2.25]  [dash pattern={on 2.53pt off 3.02pt}]  (97,128) -- (67.12,153.98) ;

\draw [color={rgb, 255:red, 0; green, 0; blue, 0 }  ,draw opacity=1 ][line width=2.25]  [dash pattern={on 2.53pt off 3.02pt}]  (134.97,153.98) -- (158,128) ;

\draw [color={rgb, 255:red, 144; green, 19; blue, 254 }  ,draw opacity=1 ][line width=2.25]    (97,128) -- (158,128) ;

\draw [color={rgb, 255:red, 144; green, 19; blue, 254 }  ,draw opacity=1 ][line width=2.25]    (158,128) -- (158,193.85) ;

\draw [color={rgb, 255:red, 144; green, 19; blue, 254 }  ,draw opacity=1 ][line width=2.25]  [dash pattern={on 2.53pt off 3.02pt}]  (158,193.85) -- (134.97,214.4) ;

\draw [color={rgb, 255:red, 0; green, 0; blue, 0 }  ,draw opacity=1 ][line width=2.25]  [dash pattern={on 2.53pt off 3.02pt}]  (97,128) -- (97,193.85) ;

\draw [color={rgb, 255:red, 144; green, 19; blue, 254 }  ,draw opacity=1 ][line width=2.25]  [dash pattern={on 2.53pt off 3.02pt}]  (97,193.85) -- (66.5,215) ;

\draw [color={rgb, 255:red, 0; green, 0; blue, 0 }  ,draw opacity=1 ][line width=2.25]  [dash pattern={on 2.53pt off 3.02pt}]  (97,193.85) -- (158,193.85) ;

\draw [color={rgb, 255:red, 144; green, 19; blue, 254 }  ,draw opacity=1 ][line width=2.25]    (67.12,153.98) -- (66.5,215) ;

\draw [color={rgb, 255:red, 0; green, 0; blue, 0 }  ,draw opacity=1 ][line width=2.25]  [dash pattern={on 2.53pt off 3.02pt}]  (339.97,199.98) -- (339.97,260.4) ;

\draw [color={rgb, 255:red, 144; green, 19; blue, 254 }  ,draw opacity=1 ][line width=2.25]    (272.12,199.98) -- (339.97,199.98) ;

\draw [color={rgb, 255:red, 0; green, 0; blue, 0 }  ,draw opacity=1 ][line width=2.25]  [dash pattern={on 2.53pt off 3.02pt}]  (271.5,261) -- (339.97,260.4) ;

\draw [color={rgb, 255:red, 208; green, 2; blue, 27 }  ,draw opacity=1 ][line width=2.25]  [dash pattern={on 2.53pt off 3.02pt}]  (302,174) -- (272.12,199.98) ;

\draw [color={rgb, 255:red, 245; green, 166; blue, 35 }  ,draw opacity=1 ][line width=2.25]  [dash pattern={on 2.53pt off 3.02pt}]  (339.97,199.98) -- (363,174) ;

\draw [color={rgb, 255:red, 0; green, 0; blue, 0 }  ,draw opacity=1 ][line width=2.25]  [dash pattern={on 2.53pt off 3.02pt}]  (302,174) -- (363,174) ;

\draw [color={rgb, 255:red, 144; green, 19; blue, 254 }  ,draw opacity=1 ][line width=2.25]    (363,174) -- (363,239.85) ;

\draw [color={rgb, 255:red, 208; green, 2; blue, 27 }  ,draw opacity=1 ][line width=2.25]  [dash pattern={on 2.53pt off 3.02pt}]  (363,239.85) -- (339.97,260.4) ;

\draw [color={rgb, 255:red, 0; green, 0; blue, 0 }  ,draw opacity=1 ][line width=2.25]  [dash pattern={on 2.53pt off 3.02pt}]  (302,174) -- (302,239.85) ;

\draw [color={rgb, 255:red, 245; green, 166; blue, 35 }  ,draw opacity=1 ][line width=2.25]  [dash pattern={on 2.53pt off 3.02pt}]  (302,239.85) -- (271.5,261) ;

\draw [color={rgb, 255:red, 144; green, 19; blue, 254 }  ,draw opacity=1 ][line width=2.25]    (302,239.85) -- (363,239.85) ;

\draw [color={rgb, 255:red, 144; green, 19; blue, 254 }  ,draw opacity=1 ][line width=2.25]    (272.12,199.98) -- (271.5,261) ;

\draw [color={rgb, 255:red, 0; green, 0; blue, 0 }  ,draw opacity=1 ][line width=2.25]  [dash pattern={on 2.53pt off 3.02pt}]  (556.97,147.98) -- (556.97,208.4) ;

\draw [color={rgb, 255:red, 144; green, 19; blue, 254 }  ,draw opacity=1 ][line width=2.25]    (489.12,147.98) -- (556.97,147.98) ;

\draw [color={rgb, 255:red, 0; green, 0; blue, 0 }  ,draw opacity=1 ][line width=2.25]  [dash pattern={on 2.53pt off 3.02pt}]  (488.5,209) -- (556.97,208.4) ;

\draw [color={rgb, 255:red, 0; green, 0; blue, 0 }  ,draw opacity=1 ][line width=2.25]  [dash pattern={on 2.53pt off 3.02pt}]  (519,122) -- (489.12,147.98) ;

\draw [color={rgb, 255:red, 144; green, 19; blue, 254 }  ,draw opacity=1 ][line width=2.25]  [dash pattern={on 2.53pt off 3.02pt}]  (556.97,147.98) -- (580,122) ;

\draw [color={rgb, 255:red, 0; green, 0; blue, 0 }  ,draw opacity=1 ][line width=2.25]  [dash pattern={on 2.53pt off 3.02pt}]  (519,122) -- (580,122) ;

\draw [color={rgb, 255:red, 0; green, 0; blue, 0 }  ,draw opacity=1 ][line width=2.25]  [dash pattern={on 2.53pt off 3.02pt}]  (580,122) -- (580,187.85) ;

\draw [color={rgb, 255:red, 144; green, 19; blue, 254 }  ,draw opacity=1 ][line width=2.25]  [dash pattern={on 2.53pt off 3.02pt}]  (580,187.85) -- (556.97,208.4) ;

\draw [color={rgb, 255:red, 144; green, 19; blue, 254 }  ,draw opacity=1 ][line width=2.25]    (519,122) -- (519,187.85) ;

\draw [color={rgb, 255:red, 0; green, 0; blue, 0 }  ,draw opacity=1 ][line width=2.25]  [dash pattern={on 2.53pt off 3.02pt}]  (519,187.85) -- (488.5,209) ;

\draw [color={rgb, 255:red, 144; green, 19; blue, 254 }  ,draw opacity=1 ][line width=2.25]    (519,187.85) -- (580,187.85) ;

\draw [color={rgb, 255:red, 144; green, 19; blue, 254 }  ,draw opacity=1 ][line width=2.25]    (489.12,147.98) -- (488.5,209) ;

\draw    (106,228) -- (105.04,281) ;
\draw [shift={(105,283)}, rotate = 271.04] [color={rgb, 255:red, 0; green, 0; blue, 0 }  ][line width=0.75]    (10.93,-3.29) .. controls (6.95,-1.4) and (3.31,-0.3) .. (0,0) .. controls (3.31,0.3) and (6.95,1.4) .. (10.93,3.29)   ;

\draw [color={rgb, 255:red, 208; green, 2; blue, 27 }  ,draw opacity=1 ][line width=1.5]    (314,273) -- (331.54,304.38) ;
\draw [shift={(333,307)}, rotate = 240.8] [color={rgb, 255:red, 208; green, 2; blue, 27 }  ,draw opacity=1 ][line width=1.5]    (14.21,-4.28) .. controls (9.04,-1.82) and (4.3,-0.39) .. (0,0) .. controls (4.3,0.39) and (9.04,1.82) .. (14.21,4.28)   ;

\draw [color={rgb, 255:red, 245; green, 166; blue, 35 }  ,draw opacity=1 ][line width=1.5]    (306,272) -- (288.43,304.36) ;
\draw [shift={(287,307)}, rotate = 298.5] [color={rgb, 255:red, 245; green, 166; blue, 35 }  ,draw opacity=1 ][line width=1.5]    (14.21,-4.28) .. controls (9.04,-1.82) and (4.3,-0.39) .. (0,0) .. controls (4.3,0.39) and (9.04,1.82) .. (14.21,4.28)   ;

\draw    (530,219) -- (530.97,277) ;
\draw [shift={(531,279)}, rotate = 269.05] [color={rgb, 255:red, 0; green, 0; blue, 0 }  ][line width=0.75]    (10.93,-3.29) .. controls (6.95,-1.4) and (3.31,-0.3) .. (0,0) .. controls (3.31,0.3) and (6.95,1.4) .. (10.93,3.29)   ;

\draw [color={rgb, 255:red, 74; green, 144; blue, 226 }  ,draw opacity=1 ][line width=2.25]    (385.97,346.98) -- (385.97,407.4) ;

\draw [color={rgb, 255:red, 74; green, 144; blue, 226 }  ,draw opacity=1 ][line width=2.25]    (317.5,408) -- (385.97,407.4) ;

\draw [color={rgb, 255:red, 208; green, 2; blue, 27 }  ,draw opacity=1 ][line width=2.25]    (348,321) -- (318.12,346.98) ;

\draw [color={rgb, 255:red, 74; green, 144; blue, 226 }  ,draw opacity=1 ][line width=2.25]    (385.97,346.98) -- (409,321) ;

\draw [color={rgb, 255:red, 74; green, 144; blue, 226 }  ,draw opacity=1 ][line width=2.25]    (348,321) -- (409,321) ;

\draw [color={rgb, 255:red, 208; green, 2; blue, 27 }  ,draw opacity=1 ][line width=2.25]    (409,321) -- (409,386.85) ;

\draw [color={rgb, 255:red, 208; green, 2; blue, 27 }  ,draw opacity=1 ][line width=2.25]    (409,386.85) -- (385.97,407.4) ;

\draw [color={rgb, 255:red, 74; green, 144; blue, 226 }  ,draw opacity=1 ][line width=2.25]    (348,321) -- (348,386.85) ;

\draw [color={rgb, 255:red, 74; green, 144; blue, 226 }  ,draw opacity=1 ][line width=2.25]    (348,386.85) -- (317.5,408) ;

\draw [color={rgb, 255:red, 208; green, 2; blue, 27 }  ,draw opacity=1 ][line width=2.25]    (348,386.85) -- (409,386.85) ;

\draw [color={rgb, 255:red, 208; green, 2; blue, 27 }  ,draw opacity=1 ][line width=2.25]    (318.12,346.98) -- (317.5,408) ;

\draw [color={rgb, 255:red, 74; green, 144; blue, 226 }  ,draw opacity=1 ][line width=2.25]    (272.97,348.98) -- (272.97,409.4) ;

\draw [color={rgb, 255:red, 74; green, 144; blue, 226 }  ,draw opacity=1 ][line width=2.25]    (204.5,410) -- (272.97,409.4) ;

\draw [color={rgb, 255:red, 74; green, 144; blue, 226 }  ,draw opacity=1 ][line width=2.25]    (235,323) -- (205.12,348.98) ;

\draw [color={rgb, 255:red, 245; green, 166; blue, 35 }  ,draw opacity=1 ][line width=2.25]    (272.97,348.98) -- (296,323) ;

\draw [color={rgb, 255:red, 74; green, 144; blue, 226 }  ,draw opacity=1 ][line width=2.25]    (235,323) -- (296,323) ;

\draw [color={rgb, 255:red, 245; green, 166; blue, 35 }  ,draw opacity=1 ][line width=2.25]    (296,323) -- (296,388.85) ;

\draw [color={rgb, 255:red, 74; green, 144; blue, 226 }  ,draw opacity=1 ][line width=2.25]    (296,388.85) -- (272.97,409.4) ;

\draw [color={rgb, 255:red, 74; green, 144; blue, 226 }  ,draw opacity=1 ][line width=2.25]    (235,323) -- (235,388.85) ;

\draw [color={rgb, 255:red, 245; green, 166; blue, 35 }  ,draw opacity=1 ][line width=2.25]    (235,388.85) -- (204.5,410) ;

\draw [color={rgb, 255:red, 245; green, 166; blue, 35 }  ,draw opacity=1 ][line width=2.25]    (235,388.85) -- (296,388.85) ;

\draw [color={rgb, 255:red, 245; green, 166; blue, 35 }  ,draw opacity=1 ][line width=2.25]    (205.12,348.98) -- (204.5,410) ;

\draw [color={rgb, 255:red, 208; green, 2; blue, 27 }  ,draw opacity=1 ][line width=2.25]    (318.12,346.98) -- (385.97,346.98) ;

\draw [color={rgb, 255:red, 245; green, 166; blue, 35 }  ,draw opacity=1 ][line width=2.25]    (205.12,348.98) -- (272.97,348.98) ;

\draw [color={rgb, 255:red, 74; green, 144; blue, 226 }  ,draw opacity=1 ][line width=2.25]    (123.97,318.98) -- (123.97,379.4) ;

\draw [color={rgb, 255:red, 74; green, 144; blue, 226 }  ,draw opacity=1 ][line width=2.25]    (55.5,380) -- (123.97,379.4) ;

\draw [color={rgb, 255:red, 74; green, 144; blue, 226 }  ,draw opacity=1 ][line width=2.25]    (86,293) -- (56.12,318.98) ;

\draw [color={rgb, 255:red, 74; green, 144; blue, 226 }  ,draw opacity=1 ][line width=2.25]    (123.97,318.98) -- (147,293) ;

\draw [color={rgb, 255:red, 144; green, 19; blue, 254 }  ,draw opacity=1 ][line width=2.25]    (86,293) -- (147,293) ;

\draw [color={rgb, 255:red, 144; green, 19; blue, 254 }  ,draw opacity=1 ][line width=2.25]    (147,293) -- (147,358.85) ;

\draw [color={rgb, 255:red, 144; green, 19; blue, 254 }  ,draw opacity=1 ][line width=2.25]    (147,358.85) -- (123.97,379.4) ;

\draw [color={rgb, 255:red, 74; green, 144; blue, 226 }  ,draw opacity=1 ][line width=2.25]    (86,293) -- (86,358.85) ;

\draw [color={rgb, 255:red, 144; green, 19; blue, 254 }  ,draw opacity=1 ][line width=2.25]    (86,358.85) -- (55.5,380) ;

\draw [color={rgb, 255:red, 74; green, 144; blue, 226 }  ,draw opacity=1 ][line width=2.25]    (86,358.85) -- (147,358.85) ;

\draw [color={rgb, 255:red, 144; green, 19; blue, 254 }  ,draw opacity=1 ][line width=2.25]    (56.12,318.98) -- (55.5,380) ;

\draw [color={rgb, 255:red, 144; green, 19; blue, 254 }  ,draw opacity=1 ][line width=2.25]    (56.12,318.98) -- (123.97,318.98) ;

\draw [color={rgb, 255:red, 74; green, 144; blue, 226 }  ,draw opacity=1 ][line width=2.25]    (546.97,315.98) -- (546.97,376.4) ;

\draw [color={rgb, 255:red, 74; green, 144; blue, 226 }  ,draw opacity=1 ][line width=2.25]    (478.5,377) -- (546.97,376.4) ;

\draw [color={rgb, 255:red, 74; green, 144; blue, 226 }  ,draw opacity=1 ][line width=2.25]    (509,290) -- (479.12,315.98) ;

\draw [color={rgb, 255:red, 144; green, 19; blue, 254 }  ,draw opacity=1 ][line width=2.25]    (546.97,315.98) -- (570,290) ;

\draw [color={rgb, 255:red, 74; green, 144; blue, 226 }  ,draw opacity=1 ][line width=2.25]    (509,290) -- (570,290) ;

\draw [color={rgb, 255:red, 74; green, 144; blue, 226 }  ,draw opacity=1 ][line width=2.25]    (570,290) -- (570,355.85) ;

\draw [color={rgb, 255:red, 144; green, 19; blue, 254 }  ,draw opacity=1 ][line width=2.25]    (570,355.85) -- (546.97,376.4) ;

\draw [color={rgb, 255:red, 144; green, 19; blue, 254 }  ,draw opacity=1 ][line width=2.25]    (509,290) -- (509,355.85) ;

\draw [color={rgb, 255:red, 74; green, 144; blue, 226 }  ,draw opacity=1 ][line width=2.25]    (509,355.85) -- (478.5,377) ;

\draw [color={rgb, 255:red, 144; green, 19; blue, 254 }  ,draw opacity=1 ][line width=2.25]    (509,355.85) -- (570,355.85) ;

\draw [color={rgb, 255:red, 144; green, 19; blue, 254 }  ,draw opacity=1 ][line width=2.25]    (479.12,315.98) -- (478.5,377) ;

\draw [color={rgb, 255:red, 144; green, 19; blue, 254 }  ,draw opacity=1 ][line width=2.25]    (479.12,315.98) -- (546.97,315.98) ;

\draw [rotate around= { 46.67: (309.44, 260.92)
    }] [color={rgb, 255:red, 144; green, 19; blue, 254 }  ,draw opacity=1 ][line width=1.5]  (300.52,260.92) -- (318.36,260.92)(309.44,251.34) -- (309.44,270.5) ;
\draw [rotate around= { 46.67: (302.44, 214.92)
    }] [color={rgb, 255:red, 144; green, 19; blue, 254 }  ,draw opacity=1 ][line width=1.5]  (293.52,214.92) -- (311.36,214.92)(302.44,205.34) -- (302.44,224.5) ;
\draw [rotate around= { 46.67: (340.44, 218.92)
    }] [color={rgb, 255:red, 144; green, 19; blue, 254 }  ,draw opacity=1 ][line width=1.5]  (331.52,218.92) -- (349.36,218.92)(340.44,209.34) -- (340.44,228.5) ;
\draw [rotate around= { 46.67: (328.44, 173.92)
    }] [color={rgb, 255:red, 144; green, 19; blue, 254 }  ,draw opacity=1 ][line width=1.5]  (319.52,173.92) -- (337.36,173.92)(328.44,164.34) -- (328.44,183.5) ;
\draw [rotate around= { 46.67: (144.44, 140.92)
    }] [color={rgb, 255:red, 144; green, 19; blue, 254 }  ,draw opacity=1 ][line width=1.5]  (135.52,140.92) -- (153.36,140.92)(144.44,131.34) -- (144.44,150.5) ;
\draw [rotate around= { 46.67: (117.44, 192.92)
    }] [color={rgb, 255:red, 144; green, 19; blue, 254 }  ,draw opacity=1 ][line width=1.5]  (108.52,192.92) -- (126.36,192.92)(117.44,183.34) -- (117.44,202.5) ;
\draw [rotate around= { 46.67: (97.44, 167.92)
    }] [color={rgb, 255:red, 144; green, 19; blue, 254 }  ,draw opacity=1 ][line width=1.5]  (88.52,167.92) -- (106.36,167.92)(97.44,158.34) -- (97.44,177.5) ;
\draw [rotate around= { 46.67: (81.44, 140.92)
    }] [color={rgb, 255:red, 144; green, 19; blue, 254 }  ,draw opacity=1 ][line width=1.5]  (72.52,140.92) -- (90.36,140.92)(81.44,131.34) -- (81.44,150.5) ;
\draw [rotate around= { 46.67: (133.44, 174.92)
    }] [color={rgb, 255:red, 144; green, 19; blue, 254 }  ,draw opacity=1 ][line width=1.5]  (124.52,174.92) -- (142.36,174.92)(133.44,165.34) -- (133.44,184.5) ;
\draw [rotate around= { 46.67: (103.44, 214.92)
    }] [color={rgb, 255:red, 144; green, 19; blue, 254 }  ,draw opacity=1 ][line width=1.5]  (94.52,214.92) -- (112.36,214.92)(103.44,205.34) -- (103.44,224.5) ;
\draw [rotate around= { 46.67: (550.44, 121.92)
    }] [color={rgb, 255:red, 144; green, 19; blue, 254 }  ,draw opacity=1 ][line width=1.5]  (541.52,121.92) -- (559.36,121.92)(550.44,112.34) -- (550.44,131.5) ;
\draw [rotate around= { 46.67: (580.44, 156.92)
    }] [color={rgb, 255:red, 144; green, 19; blue, 254 }  ,draw opacity=1 ][line width=1.5]  (571.52,156.92) -- (589.36,156.92)(580.44,147.34) -- (580.44,166.5) ;
\draw [rotate around= { 46.67: (506.44, 195.92)
    }] [color={rgb, 255:red, 144; green, 19; blue, 254 }  ,draw opacity=1 ][line width=1.5]  (497.52,195.92) -- (515.36,195.92)(506.44,186.34) -- (506.44,205.5) ;
\draw [rotate around= { 46.67: (528.44, 208.92)
    }] [color={rgb, 255:red, 144; green, 19; blue, 254 }  ,draw opacity=1 ][line width=1.5]  (519.52,208.92) -- (537.36,208.92)(528.44,199.34) -- (528.44,218.5) ;
\draw [rotate around= { 46.67: (502.44, 134.92)
    }] [color={rgb, 255:red, 144; green, 19; blue, 254 }  ,draw opacity=1 ][line width=1.5]  (493.52,134.92) -- (511.36,134.92)(502.44,125.34) -- (502.44,144.5) ;
\draw [rotate around= { 46.67: (557.44, 168.92)
    }] [color={rgb, 255:red, 144; green, 19; blue, 254 }  ,draw opacity=1 ][line width=1.5]  (548.52,168.92) -- (566.36,168.92)(557.44,159.34) -- (557.44,178.5) ;

\draw (215,85) node  [align=left] {{\small $\displaystyle \dfrac{\pi }{2}$ CW }};
\draw (527,404) node  [align=left] {FD (1)};
\draw (361,426) node  [align=left] {FD (2)};
\draw (244,430) node  [align=left] {FD (2)};
\draw (100,399) node  [align=left] {FD (1)};
\draw (352,134) node  [align=left] {$\displaystyle \pi $ CW};
\draw (445,86) node  [align=left] {$\displaystyle \dfrac{3\pi }{2}$ CW};

\end{tikzpicture}

    \caption{Ways to put 6 edges in a class with 2 edges per face. A purple X correspond to an edge that would lead to a face with more than 2 edges on it.}
    \label{fig:my_label}
\end{figure}  

Thus, in this case we always end up with either FD (1) or FD (2).

\textbf{Opposite edges, opposite faces}:

\begin{figure}[H]
    \centering
\tikzset{every picture/.style={line width=0.4pt}} 

\begin{tikzpicture}[x=0.5pt,y=0.5pt,yscale=-1,xscale=1]

\draw [color={rgb, 255:red, 144; green, 19; blue, 254 }  ,draw opacity=1 ][line width=2.25]    (163.97,118.98) -- (163.97,179.4) ;

\draw [color={rgb, 255:red, 0; green, 0; blue, 0 }  ,draw opacity=1 ][line width=2.25]  [dash pattern={on 2.53pt off 3.02pt}]  (96.12,118.98) -- (163.97,118.98) ;

\draw [color={rgb, 255:red, 0; green, 0; blue, 0 }  ,draw opacity=1 ][line width=2.25]  [dash pattern={on 2.53pt off 3.02pt}]  (95.5,180) -- (163.97,179.4) ;

\draw [color={rgb, 255:red, 0; green, 0; blue, 0 }  ,draw opacity=1 ][line width=2.25]  [dash pattern={on 2.53pt off 3.02pt}]  (126,93) -- (96.12,118.98) ;

\draw [color={rgb, 255:red, 0; green, 0; blue, 0 }  ,draw opacity=1 ][line width=2.25]  [dash pattern={on 2.53pt off 3.02pt}]  (163.97,118.98) -- (187,93) ;

\draw [color={rgb, 255:red, 0; green, 0; blue, 0 }  ,draw opacity=1 ][line width=2.25]  [dash pattern={on 2.53pt off 3.02pt}]  (126,93) -- (187,93) ;

\draw [color={rgb, 255:red, 144; green, 19; blue, 254 }  ,draw opacity=1 ][line width=2.25]    (187,93) -- (187,158.85) ;

\draw [color={rgb, 255:red, 0; green, 0; blue, 0 }  ,draw opacity=1 ][line width=2.25]  [dash pattern={on 2.53pt off 3.02pt}]  (187,158.85) -- (163.97,179.4) ;

\draw [color={rgb, 255:red, 144; green, 19; blue, 254 }  ,draw opacity=1 ][line width=2.25]    (126,93) -- (126,158.85) ;

\draw [color={rgb, 255:red, 0; green, 0; blue, 0 }  ,draw opacity=1 ][line width=2.25]  [dash pattern={on 2.53pt off 3.02pt}]  (126,158.85) -- (95.5,180) ;

\draw [color={rgb, 255:red, 0; green, 0; blue, 0 }  ,draw opacity=1 ][line width=2.25]  [dash pattern={on 2.53pt off 3.02pt}]  (126,158.85) -- (187,158.85) ;

\draw [color={rgb, 255:red, 144; green, 19; blue, 254 }  ,draw opacity=1 ][line width=2.25]    (96.12,118.98) -- (95.5,180) ;

\draw [color={rgb, 255:red, 144; green, 19; blue, 254 }  ,draw opacity=1 ][line width=2.25]    (379.97,117.98) -- (379.97,178.4) ;

\draw [color={rgb, 255:red, 0; green, 0; blue, 0 }  ,draw opacity=1 ][line width=2.25]  [dash pattern={on 2.53pt off 3.02pt}]  (312.12,117.98) -- (379.97,117.98) ;

\draw [color={rgb, 255:red, 0; green, 0; blue, 0 }  ,draw opacity=1 ][line width=2.25]  [dash pattern={on 2.53pt off 3.02pt}]  (311.5,179) -- (379.97,178.4) ;

\draw [color={rgb, 255:red, 80; green, 227; blue, 194 }  ,draw opacity=1 ][line width=2.25]    (342,92) -- (312.12,117.98) ;

\draw [color={rgb, 255:red, 208; green, 2; blue, 27 }  ,draw opacity=1 ][line width=2.25]    (379.97,117.98) -- (403,92) ;

\draw [color={rgb, 255:red, 144; green, 19; blue, 254 }  ,draw opacity=1 ][line width=2.25]    (342,92) -- (403,92) ;

\draw [color={rgb, 255:red, 0; green, 0; blue, 0 }  ,draw opacity=1 ][line width=2.25]  [dash pattern={on 2.53pt off 3.02pt}]  (403,92) -- (403,157.85) ;

\draw [color={rgb, 255:red, 80; green, 227; blue, 194 }  ,draw opacity=1 ][line width=2.25]    (403,157.85) -- (379.97,178.4) ;

\draw [color={rgb, 255:red, 0; green, 0; blue, 0 }  ,draw opacity=1 ][line width=2.25]  [dash pattern={on 2.53pt off 3.02pt}]  (342,92) -- (342,157.85) ;

\draw [color={rgb, 255:red, 208; green, 2; blue, 27 }  ,draw opacity=1 ][line width=2.25]    (342,157.85) -- (311.5,179) ;

\draw [color={rgb, 255:red, 144; green, 19; blue, 254 }  ,draw opacity=1 ][line width=2.25]    (342,157.85) -- (403,157.85) ;

\draw [color={rgb, 255:red, 144; green, 19; blue, 254 }  ,draw opacity=1 ][line width=2.25]    (312.12,117.98) -- (311.5,179) ;

\draw [color={rgb, 255:red, 80; green, 227; blue, 194 }  ,draw opacity=1 ][line width=1.5]    (419,118) -- (445.77,94) ;
\draw [shift={(448,92)}, rotate = 498.12] [color={rgb, 255:red, 80; green, 227; blue, 194 }  ,draw opacity=1 ][line width=1.5]    (14.21,-4.28) .. controls (9.04,-1.82) and (4.3,-0.39) .. (0,0) .. controls (4.3,0.39) and (9.04,1.82) .. (14.21,4.28)   ;

\draw [color={rgb, 255:red, 208; green, 2; blue, 27 }  ,draw opacity=1 ][line width=1.5]    (416,148) -- (448.43,167.46) ;
\draw [shift={(451,169)}, rotate = 210.96] [color={rgb, 255:red, 208; green, 2; blue, 27 }  ,draw opacity=1 ][line width=1.5]    (14.21,-4.28) .. controls (9.04,-1.82) and (4.3,-0.39) .. (0,0) .. controls (4.3,0.39) and (9.04,1.82) .. (14.21,4.28)   ;

\draw [color={rgb, 255:red, 80; green, 227; blue, 194 }  ,draw opacity=1 ][line width=2.25]    (528.97,48.98) -- (528.97,109.4) ;

\draw [color={rgb, 255:red, 74; green, 144; blue, 226 }  ,draw opacity=1 ][line width=2.25]    (460.5,110) -- (528.97,109.4) ;

\draw [color={rgb, 255:red, 80; green, 227; blue, 194 }  ,draw opacity=1 ][line width=2.25]    (491,23) -- (461.12,48.98) ;

\draw [color={rgb, 255:red, 74; green, 144; blue, 226 }  ,draw opacity=1 ][line width=2.25]    (528.97,48.98) -- (552,23) ;

\draw [color={rgb, 255:red, 80; green, 227; blue, 194 }  ,draw opacity=1 ][line width=2.25]    (491,23) -- (552,23) ;

\draw [color={rgb, 255:red, 74; green, 144; blue, 226 }  ,draw opacity=1 ][line width=2.25]    (552,23) -- (552,88.85) ;

\draw [color={rgb, 255:red, 80; green, 227; blue, 194 }  ,draw opacity=1 ][line width=2.25]    (552,88.85) -- (528.97,109.4) ;

\draw [color={rgb, 255:red, 74; green, 144; blue, 226 }  ,draw opacity=1 ][line width=2.25]    (491,23) -- (491,88.85) ;

\draw [color={rgb, 255:red, 74; green, 144; blue, 226 }  ,draw opacity=1 ][line width=2.25]    (491,88.85) -- (460.5,110) ;

\draw [color={rgb, 255:red, 80; green, 227; blue, 194 }  ,draw opacity=1 ][line width=2.25]    (491,88.85) -- (552,88.85) ;

\draw [color={rgb, 255:red, 80; green, 227; blue, 194 }  ,draw opacity=1 ][line width=2.25]    (461.12,48.98) -- (460.5,110) ;

\draw [color={rgb, 255:red, 74; green, 144; blue, 226 }  ,draw opacity=1 ][line width=2.25]    (461.12,48.98) -- (528.97,48.98) ;

\draw [color={rgb, 255:red, 208; green, 2; blue, 27 }  ,draw opacity=1 ][line width=2.25]    (528.97,163.98) -- (528.97,224.4) ;

\draw [color={rgb, 255:red, 74; green, 144; blue, 226 }  ,draw opacity=1 ][line width=2.25]    (460.5,225) -- (528.97,224.4) ;

\draw [color={rgb, 255:red, 74; green, 144; blue, 226 }  ,draw opacity=1 ][line width=2.25]    (491,138) -- (461.12,163.98) ;

\draw [color={rgb, 255:red, 208; green, 2; blue, 27 }  ,draw opacity=1 ][line width=2.25]    (528.97,163.98) -- (552,138) ;

\draw [color={rgb, 255:red, 208; green, 2; blue, 27 }  ,draw opacity=1 ][line width=2.25]    (491,138) -- (552,138) ;

\draw [color={rgb, 255:red, 74; green, 144; blue, 226 }  ,draw opacity=1 ][line width=2.25]    (552,138) -- (552,203.85) ;

\draw [color={rgb, 255:red, 74; green, 144; blue, 226 }  ,draw opacity=1 ][line width=2.25]    (552,203.85) -- (528.97,224.4) ;

\draw [color={rgb, 255:red, 74; green, 144; blue, 226 }  ,draw opacity=1 ][line width=2.25]    (491,138) -- (491,203.85) ;

\draw [color={rgb, 255:red, 208; green, 2; blue, 27 }  ,draw opacity=1 ][line width=2.25]    (491,203.85) -- (460.5,225) ;

\draw [color={rgb, 255:red, 208; green, 2; blue, 27 }  ,draw opacity=1 ][line width=2.25]    (491,203.85) -- (552,203.85) ;

\draw [color={rgb, 255:red, 208; green, 2; blue, 27 }  ,draw opacity=1 ][line width=2.25]    (461.12,163.98) -- (460.5,225) ;

\draw [color={rgb, 255:red, 74; green, 144; blue, 226 }  ,draw opacity=1 ][line width=2.25]    (461.12,163.98) -- (528.97,163.98) ;

\draw [rotate around= { 46.67: (402.44, 127.92)
    }] [color={rgb, 255:red, 144; green, 19; blue, 254 }  ,draw opacity=1 ][line width=1.5]  (393.52,127.92) -- (411.36,127.92)(402.44,118.34) -- (402.44,137.5) ;
\draw  [color={rgb, 255:red, 208; green, 2; blue, 27 }  ,draw opacity=0.44 ][line width=1.5]  (56,133) .. controls (56,86.06) and (93.94,48) .. (140.75,48) .. controls (187.56,48) and (225.5,86.06) .. (225.5,133) .. controls (225.5,179.94) and (187.56,218) .. (140.75,218) .. controls (93.94,218) and (56,179.94) .. (56,133) -- cycle ; \draw  [color={rgb, 255:red, 208; green, 2; blue, 27 }  ,draw opacity=0.44 ][line width=1.5]  (80.82,72.9) -- (200.68,193.1) ; \draw  [color={rgb, 255:red, 208; green, 2; blue, 27 }  ,draw opacity=0.44 ][line width=1.5]  (200.68,72.9) -- (80.82,193.1) ;
\draw [rotate around= { 46.67: (112.44, 167.92)
    }] [color={rgb, 255:red, 144; green, 19; blue, 254 }  ,draw opacity=1 ][line width=1.5]  (103.52,167.92) -- (121.36,167.92)(112.44,158.34) -- (112.44,177.5) ;
\draw [rotate around= { 46.67: (149.44, 158.92)
    }] [color={rgb, 255:red, 144; green, 19; blue, 254 }  ,draw opacity=1 ][line width=1.5]  (140.52,158.92) -- (158.36,158.92)(149.44,149.34) -- (149.44,168.5) ;
\draw [rotate around= { 46.67: (153.44, 92.92)
    }] [color={rgb, 255:red, 144; green, 19; blue, 254 }  ,draw opacity=1 ][line width=1.5]  (144.52,92.92) -- (162.36,92.92)(153.44,83.34) -- (153.44,102.5) ;
\draw [rotate around= { 46.67: (110.44, 105.92)
    }] [color={rgb, 255:red, 144; green, 19; blue, 254 }  ,draw opacity=1 ][line width=1.5]  (101.52,105.92) -- (119.36,105.92)(110.44,96.34) -- (110.44,115.5) ;
\draw [rotate around= { 46.67: (130.44, 178.92)
    }] [color={rgb, 255:red, 144; green, 19; blue, 254 }  ,draw opacity=1 ][line width=1.5]  (121.52,178.92) -- (139.36,178.92)(130.44,169.34) -- (130.44,188.5) ;
\draw [rotate around= { 46.67: (176.44, 168.92)
    }] [color={rgb, 255:red, 144; green, 19; blue, 254 }  ,draw opacity=1 ][line width=1.5]  (167.52,168.92) -- (185.36,168.92)(176.44,159.34) -- (176.44,178.5) ;
\draw [rotate around= { 46.67: (138.44, 118.92)
    }] [color={rgb, 255:red, 144; green, 19; blue, 254 }  ,draw opacity=1 ][line width=1.5]  (129.52,118.92) -- (147.36,118.92)(138.44,109.34) -- (138.44,128.5) ;
\draw [rotate around= { 46.67: (342.44, 138.92)
    }] [color={rgb, 255:red, 144; green, 19; blue, 254 }  ,draw opacity=1 ][line width=1.5]  (333.52,138.92) -- (351.36,138.92)(342.44,129.34) -- (342.44,148.5) ;
\draw [rotate around= { 46.67: (358.44, 117.92)
    }] [color={rgb, 255:red, 144; green, 19; blue, 254 }  ,draw opacity=1 ][line width=1.5]  (349.52,117.92) -- (367.36,117.92)(358.44,108.34) -- (358.44,127.5) ;
\draw [rotate around= { 46.67: (347.44, 177.92)
    }] [color={rgb, 255:red, 144; green, 19; blue, 254 }  ,draw opacity=1 ][line width=1.5]  (338.52,177.92) -- (356.36,177.92)(347.44,168.34) -- (347.44,187.5) ;
\draw [rotate around= { 46.67: (173.44, 107.92)
    }] [color={rgb, 255:red, 144; green, 19; blue, 254 }  ,draw opacity=1 ][line width=1.5]  (164.52,107.92) -- (182.36,107.92)(173.44,98.34) -- (173.44,117.5) ;

\draw (590,62) node  [align=left] {FD (1)};
\draw (590,184) node  [align=left] {FD (1)};
\draw (269,140) node  [align=left] {OR};
\end{tikzpicture}
 \caption{Ways to put opposite edges on opposite faces for 2 edges per face, 6 edges per class.}
    \label{fig:my_label}
\end{figure}
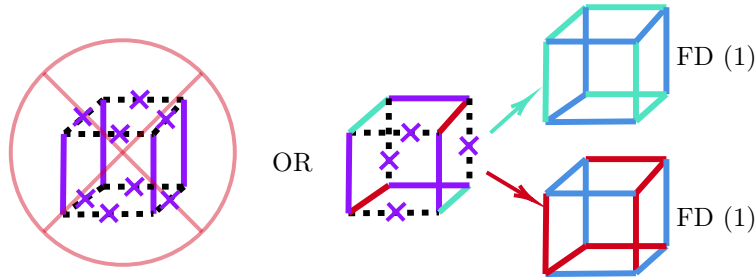

We only have two options in this case because due to symmetry, the $0$ and $\pi$ rotations look the same, as do the $\frac{\pi}{2}$ and $\frac{3\pi}{2}$. On the left of Figure 9: for the 0 rotation case, once we have placed the first four edges, there are no other places to put the remaining two edges without having more than 2 edges on a face. On the right: we have two choices of where to put the remaining two edges, and each gives rise to FD (1)

\textbf{Adjacent edges, Adjacent faces}: We need only consider the cases $\dfrac{\pi}{2}$ rotation and no rotation; otherwise, we will end up with three edges on the top face.

\begin{figure}[H]
    \centering

\tikzset{every picture/.style={line width=0.75pt}} 

\begin{tikzpicture}[x=0.75pt,y=0.75pt,yscale=-1,xscale=1]

\draw [color={rgb, 255:red, 0; green, 0; blue, 0 }  ,draw opacity=1 ][line width=2.25]  [dash pattern={on 2.53pt off 3.02pt}]  (82,119.84) -- (82,167.9) ;

\draw [color={rgb, 255:red, 144; green, 19; blue, 254 }  ,draw opacity=1 ][line width=2.25]    (28,119.84) -- (82,119.84) ;

\draw [color={rgb, 255:red, 0; green, 0; blue, 0 }  ,draw opacity=1 ][line width=2.25]  [dash pattern={on 2.53pt off 3.02pt}]  (27.5,168.38) -- (82,167.9) ;

\draw [color={rgb, 255:red, 144; green, 19; blue, 254 }  ,draw opacity=1 ][line width=2.25]    (51.78,99.18) -- (28,119.84) ;

\draw [color={rgb, 255:red, 0; green, 0; blue, 0 }  ,draw opacity=1 ][line width=2.25]  [dash pattern={on 2.53pt off 3.02pt}]  (82,119.84) -- (100.33,99.18) ;

\draw [color={rgb, 255:red, 0; green, 0; blue, 0 }  ,draw opacity=1 ][line width=2.25]  [dash pattern={on 2.53pt off 3.02pt}]  (51.78,99.18) -- (100.33,99.18) ;

\draw [color={rgb, 255:red, 0; green, 0; blue, 0 }  ,draw opacity=1 ][line width=2.25]  [dash pattern={on 2.53pt off 3.02pt}]  (100.33,99.18) -- (100.33,151.56) ;

\draw [color={rgb, 255:red, 0; green, 0; blue, 0 }  ,draw opacity=1 ][line width=2.25]  [dash pattern={on 2.53pt off 3.02pt}]  (100.33,151.56) -- (82,167.9) ;

\draw [color={rgb, 255:red, 0; green, 0; blue, 0 }  ,draw opacity=1 ][line width=2.25]  [dash pattern={on 2.53pt off 3.02pt}]  (51.78,99.18) -- (51.78,151.56) ;

\draw [color={rgb, 255:red, 0; green, 0; blue, 0 }  ,draw opacity=1 ][line width=2.25]  [dash pattern={on 2.53pt off 3.02pt}]  (51.78,151.56) -- (27.5,168.38) ;

\draw [color={rgb, 255:red, 0; green, 0; blue, 0 }  ,draw opacity=1 ][line width=2.25]  [dash pattern={on 2.53pt off 3.02pt}]  (51.78,151.56) -- (100.33,151.56) ;

\draw [color={rgb, 255:red, 144; green, 19; blue, 254 }  ,draw opacity=1 ][line width=2.25]    (28,119.84) -- (27.5,168.38) ;

\draw [color={rgb, 255:red, 0; green, 0; blue, 0 }  ,draw opacity=1 ][line width=2.25]  [dash pattern={on 2.53pt off 3.02pt}]  (82,119.84) -- (82,167.9) ;

\draw [color={rgb, 255:red, 144; green, 19; blue, 254 }  ,draw opacity=1 ][line width=2.25]    (28,119.84) -- (82,119.84) ;

\draw [color={rgb, 255:red, 0; green, 0; blue, 0 }  ,draw opacity=1 ][line width=2.25]  [dash pattern={on 2.53pt off 3.02pt}]  (27.5,168.38) -- (82,167.9) ;

\draw [color={rgb, 255:red, 0; green, 0; blue, 0 }  ,draw opacity=1 ][line width=2.25]  [dash pattern={on 2.53pt off 3.02pt}]  (82,119.84) -- (100.33,99.18) ;

\draw [color={rgb, 255:red, 0; green, 0; blue, 0 }  ,draw opacity=1 ][line width=2.25]  [dash pattern={on 2.53pt off 3.02pt}]  (51.78,99.18) -- (100.33,99.18) ;

\draw [color={rgb, 255:red, 144; green, 19; blue, 254 }  ,draw opacity=1 ][line width=2.25]  [dash pattern={on 2.53pt off 3.02pt}]  (100.33,99.18) -- (100.33,151.56) ;

\draw [color={rgb, 255:red, 144; green, 19; blue, 254 }  ,draw opacity=1 ][line width=2.25]  [dash pattern={on 2.53pt off 3.02pt}]  (100.33,151.56) -- (82,167.9) ;

\draw [color={rgb, 255:red, 0; green, 0; blue, 0 }  ,draw opacity=1 ][line width=2.25]  [dash pattern={on 2.53pt off 3.02pt}]  (51.78,99.18) -- (51.78,151.56) ;

\draw [color={rgb, 255:red, 0; green, 0; blue, 0 }  ,draw opacity=1 ][line width=2.25]  [dash pattern={on 2.53pt off 3.02pt}]  (51.78,151.56) -- (27.5,168.38) ;

\draw [color={rgb, 255:red, 144; green, 19; blue, 254 }  ,draw opacity=1 ][line width=2.25]  [dash pattern={on 2.53pt off 3.02pt}]  (51.78,151.56) -- (100.33,151.56) ;

\draw [color={rgb, 255:red, 144; green, 19; blue, 254 }  ,draw opacity=1 ][line width=2.25]    (28,119.84) -- (27.5,168.38) ;

\draw [rotate around= { 46.67: (206.56, 132.52)
    }] [color={rgb, 255:red, 144; green, 19; blue, 254 }  ,draw opacity=1 ][line width=1.5]  (199.46,132.52) -- (213.66,132.52)(206.56,124.9) -- (206.56,140.14) ;
\draw [rotate around= { 46.67: (185.06, 166.72)
    }] [color={rgb, 255:red, 144; green, 19; blue, 254 }  ,draw opacity=1 ][line width=1.5]  (177.96,166.72) -- (192.17,166.72)(185.06,159.1) -- (185.06,174.35) ;
\draw [rotate around= { 46.67: (90.34, 110.25)
    }] [color={rgb, 255:red, 144; green, 19; blue, 254 }  ,draw opacity=1 ][line width=1.5]  (83.24,110.25) -- (97.44,110.25)(90.34,102.62) -- (90.34,117.87) ;
\draw [rotate around= { 46.67: (71.24, 99.11)
    }] [color={rgb, 255:red, 144; green, 19; blue, 254 }  ,draw opacity=1 ][line width=1.5]  (64.14,99.11) -- (78.34,99.11)(71.24,91.49) -- (71.24,106.73) ;
\draw [rotate around= { 46.67: (51.34, 129.34)
    }] [color={rgb, 255:red, 144; green, 19; blue, 254 }  ,draw opacity=1 ][line width=1.5]  (44.24,129.34) -- (58.44,129.34)(51.34,121.71) -- (51.34,136.96) ;
\draw [rotate around= { 46.67: (41.78, 157.97)
    }] [color={rgb, 255:red, 144; green, 19; blue, 254 }  ,draw opacity=1 ][line width=1.5]  (34.68,157.97) -- (48.88,157.97)(41.78,150.35) -- (41.78,165.6) ;
\draw [rotate around= { 46.67: (59.3, 168.32)
    }] [color={rgb, 255:red, 144; green, 19; blue, 254 }  ,draw opacity=1 ][line width=1.5]  (52.2,168.32) -- (66.4,168.32)(59.3,160.69) -- (59.3,175.94) ;
\draw [color={rgb, 255:red, 0; green, 0; blue, 0 }  ,draw opacity=1 ][line width=2.25]  [dash pattern={on 2.53pt off 3.02pt}]  (206.18,118.25) -- (206.18,166.31) ;

\draw [color={rgb, 255:red, 144; green, 19; blue, 254 }  ,draw opacity=1 ][line width=2.25]    (152.17,118.25) -- (206.18,118.25) ;

\draw [color={rgb, 255:red, 0; green, 0; blue, 0 }  ,draw opacity=1 ][line width=2.25]  [dash pattern={on 2.53pt off 3.02pt}]  (151.68,166.79) -- (206.18,166.31) ;

\draw [color={rgb, 255:red, 0; green, 0; blue, 0 }  ,draw opacity=1 ][line width=2.25]  [dash pattern={on 2.53pt off 3.02pt}]  (175.96,97.58) -- (152.17,118.25) ;

\draw [color={rgb, 255:red, 144; green, 19; blue, 254 }  ,draw opacity=1 ][line width=2.25]    (206.18,118.25) -- (224.51,97.58) ;

\draw [color={rgb, 255:red, 0; green, 0; blue, 0 }  ,draw opacity=1 ][line width=2.25]  [dash pattern={on 2.53pt off 3.02pt}]  (175.96,97.58) -- (224.51,97.58) ;

\draw [color={rgb, 255:red, 0; green, 0; blue, 0 }  ,draw opacity=1 ][line width=2.25]  [dash pattern={on 2.53pt off 3.02pt}]  (224.51,97.58) -- (224.51,149.97) ;

\draw [color={rgb, 255:red, 0; green, 0; blue, 0 }  ,draw opacity=1 ][line width=2.25]  [dash pattern={on 2.53pt off 3.02pt}]  (224.51,149.97) -- (206.18,166.31) ;

\draw [color={rgb, 255:red, 0; green, 0; blue, 0 }  ,draw opacity=1 ][line width=2.25]  [dash pattern={on 2.53pt off 3.02pt}]  (175.96,97.58) -- (175.96,149.97) ;

\draw [color={rgb, 255:red, 0; green, 0; blue, 0 }  ,draw opacity=1 ][line width=2.25]  [dash pattern={on 2.53pt off 3.02pt}]  (175.96,149.97) -- (151.68,166.79) ;

\draw [color={rgb, 255:red, 0; green, 0; blue, 0 }  ,draw opacity=1 ][line width=2.25]  [dash pattern={on 2.53pt off 3.02pt}]  (175.96,149.97) -- (224.51,149.97) ;

\draw [color={rgb, 255:red, 144; green, 19; blue, 254 }  ,draw opacity=1 ][line width=2.25]    (152.17,118.25) -- (151.68,166.79) ;

\draw [rotate around= { 46.67: (300.89, 217.61)
    }] [color={rgb, 255:red, 144; green, 19; blue, 254 }  ,draw opacity=1 ][line width=1.5]  (293.79,217.61) -- (307.99,217.61)(300.89,209.99) -- (300.89,225.23) ;
\draw [rotate around= { 46.67: (283.38, 157.95)
    }] [color={rgb, 255:red, 144; green, 19; blue, 254 }  ,draw opacity=1 ][line width=1.5]  (276.28,157.95) -- (290.48,157.95)(283.38,150.33) -- (283.38,165.58) ;
\draw [rotate around= { 46.67: (299.69, 74.45)
    }] [color={rgb, 255:red, 126; green, 211; blue, 33 }  ,draw opacity=1 ][line width=1.5]  (292.59,74.45) -- (306.79,74.45)(299.69,66.83) -- (299.69,82.07) ;
\draw [rotate around= { 46.67: (324.77, 187.38)
    }] [color={rgb, 255:red, 144; green, 19; blue, 254 }  ,draw opacity=1 ][line width=1.5]  (317.67,187.38) -- (331.87,187.38)(324.77,179.76) -- (324.77,195.01) ;
\draw [rotate around= { 46.67: (162.78, 107.07)
    }] [color={rgb, 255:red, 144; green, 19; blue, 254 }  ,draw opacity=1 ][line width=1.5]  (155.68,107.07) -- (169.88,107.07)(162.78,99.44) -- (162.78,114.69) ;
\draw [rotate around= { 46.67: (197.01, 96.72)
    }] [color={rgb, 255:red, 144; green, 19; blue, 254 }  ,draw opacity=1 ][line width=1.5]  (189.9,96.72) -- (204.11,96.72)(197.01,89.1) -- (197.01,104.35) ;
\draw [rotate around= { 46.67: (330.73, 76.04)
    }] [color={rgb, 255:red, 144; green, 19; blue, 254 }  ,draw opacity=1 ][line width=1.5]  (323.63,76.04) -- (337.83,76.04)(330.73,68.42) -- (330.73,83.67) ;
\draw [rotate around= { 46.67: (306.06, 110.25)
    }] [color={rgb, 255:red, 144; green, 19; blue, 254 }  ,draw opacity=1 ][line width=1.5]  (298.96,110.25) -- (313.16,110.25)(306.06,102.62) -- (306.06,117.87) ;
\draw [color={rgb, 255:red, 0; green, 0; blue, 0 }  ,draw opacity=1 ][line width=2.25]  [dash pattern={on 2.53pt off 3.02pt}]  (330.36,61.77) -- (330.36,109.83) ;

\draw [color={rgb, 255:red, 144; green, 19; blue, 254 }  ,draw opacity=1 ][line width=2.25]    (276.35,61.77) -- (330.36,61.77) ;

\draw [color={rgb, 255:red, 0; green, 0; blue, 0 }  ,draw opacity=1 ][line width=2.25]  [dash pattern={on 2.53pt off 3.02pt}]  (275.85,110.31) -- (330.36,109.83) ;

\draw [color={rgb, 255:red, 0; green, 0; blue, 0 }  ,draw opacity=1 ][line width=2.25]  [dash pattern={on 2.53pt off 3.02pt}]  (300.13,41.11) -- (276.35,61.77) ;

\draw [color={rgb, 255:red, 144; green, 19; blue, 254 }  ,draw opacity=1 ][line width=2.25]    (330.36,61.77) -- (348.69,41.11) ;

\draw [color={rgb, 255:red, 0; green, 0; blue, 0 }  ,draw opacity=1 ][line width=2.25]  [dash pattern={on 2.53pt off 3.02pt}]  (300.13,41.11) -- (348.69,41.11) ;

\draw [color={rgb, 255:red, 139; green, 87; blue, 42 }  ,draw opacity=1 ][line width=2.25]    (348.69,41.11) -- (348.69,93.49) ;

\draw [color={rgb, 255:red, 0; green, 0; blue, 0 }  ,draw opacity=1 ][line width=2.25]  [dash pattern={on 2.53pt off 3.02pt}]  (348.69,93.49) -- (330.36,109.83) ;

\draw [color={rgb, 255:red, 0; green, 0; blue, 0 }  ,draw opacity=1 ][line width=2.25]  [dash pattern={on 2.53pt off 3.02pt}]  (300.13,41.11) -- (300.13,93.49) ;

\draw [color={rgb, 255:red, 126; green, 211; blue, 33 }  ,draw opacity=1 ][line width=2.25]    (300.13,93.49) -- (275.85,110.31) ;

\draw [color={rgb, 255:red, 139; green, 87; blue, 42 }  ,draw opacity=1 ][line width=2.25]    (300.13,93.49) -- (348.69,93.49) ;

\draw [color={rgb, 255:red, 144; green, 19; blue, 254 }  ,draw opacity=1 ][line width=2.25]    (276.35,61.77) -- (275.85,110.31) ;

\draw [rotate around= { 46.67: (286.95, 50.59)
    }] [color={rgb, 255:red, 144; green, 19; blue, 254 }  ,draw opacity=1 ][line width=1.5]  (279.85,50.59) -- (294.05,50.59)(286.95,42.96) -- (286.95,58.21) ;
\draw [rotate around= { 46.67: (321.18, 40.25)
    }] [color={rgb, 255:red, 144; green, 19; blue, 254 }  ,draw opacity=1 ][line width=1.5]  (314.08,40.25) -- (328.28,40.25)(321.18,32.62) -- (321.18,47.87) ;
\draw [rotate around= { 32.19: (339.49, 102.29)
    }] [color={rgb, 255:red, 208; green, 2; blue, 27 }  ,draw opacity=1 ][line width=1.5]  (332.39,102.29) -- (346.59,102.29)(339.49,94.67) -- (339.49,109.92) ;
\draw [color={rgb, 255:red, 0; green, 0; blue, 0 }  ,draw opacity=1 ][line width=2.25]  [dash pattern={on 2.53pt off 3.02pt}]  (325.19,169.14) -- (325.19,217.2) ;

\draw [color={rgb, 255:red, 144; green, 19; blue, 254 }  ,draw opacity=1 ][line width=2.25]    (271.19,169.14) -- (325.19,169.14) ;

\draw [color={rgb, 255:red, 0; green, 0; blue, 0 }  ,draw opacity=1 ][line width=2.25]  [dash pattern={on 2.53pt off 3.02pt}]  (270.69,217.68) -- (325.19,217.2) ;

\draw [color={rgb, 255:red, 0; green, 0; blue, 0 }  ,draw opacity=1 ][line width=2.25]  [dash pattern={on 2.53pt off 3.02pt}]  (294.97,148.47) -- (271.19,169.14) ;

\draw [color={rgb, 255:red, 144; green, 19; blue, 254 }  ,draw opacity=1 ][line width=2.25]    (325.19,169.14) -- (343.52,148.47) ;

\draw [color={rgb, 255:red, 0; green, 0; blue, 0 }  ,draw opacity=1 ][line width=2.25]  [dash pattern={on 2.53pt off 3.02pt}]  (294.97,148.47) -- (343.52,148.47) ;

\draw [color={rgb, 255:red, 0; green, 0; blue, 0 }  ,draw opacity=1 ][line width=2.25]  [dash pattern={on 2.53pt off 3.02pt}]  (343.52,148.47) -- (343.52,200.86) ;

\draw [color={rgb, 255:red, 126; green, 211; blue, 33 }  ,draw opacity=1 ][line width=2.25]    (343.52,200.86) -- (325.19,217.2) ;

\draw [color={rgb, 255:red, 139; green, 87; blue, 42 }  ,draw opacity=1 ][line width=2.25]    (294.97,148.47) -- (294.97,200.86) ;

\draw [color={rgb, 255:red, 0; green, 0; blue, 0 }  ,draw opacity=1 ][line width=2.25]  [dash pattern={on 2.53pt off 3.02pt}]  (294.97,200.86) -- (270.69,217.68) ;

\draw [color={rgb, 255:red, 139; green, 87; blue, 42 }  ,draw opacity=1 ][line width=2.25]    (294.97,200.86) -- (343.52,200.86) ;

\draw [color={rgb, 255:red, 144; green, 19; blue, 254 }  ,draw opacity=1 ][line width=2.25]    (271.19,169.14) -- (270.69,217.68) ;

\draw [rotate around= { 46.67: (317.61, 147.61)
    }] [color={rgb, 255:red, 144; green, 19; blue, 254 }  ,draw opacity=1 ][line width=1.5]  (310.51,147.61) -- (324.71,147.61)(317.61,139.99) -- (317.61,155.23) ;
\draw [rotate around= { 46.67: (343.08, 174.66)
    }] [color={rgb, 255:red, 126; green, 211; blue, 33 }  ,draw opacity=1 ][line width=1.5]  (335.98,174.66) -- (350.18,174.66)(343.08,167.03) -- (343.08,182.28) ;
\draw [rotate around= { 46.67: (285.77, 207.27)
    }] [color={rgb, 255:red, 208; green, 2; blue, 27 }  ,draw opacity=1 ][line width=1.5]  (278.67,207.27) -- (292.87,207.27)(285.77,199.65) -- (285.77,214.89) ;
\draw    (234,130) -- (267.95,102.27) ;
\draw [shift={(269.5,101)}, rotate = 500.75] [color={rgb, 255:red, 0; green, 0; blue, 0 }  ][line width=0.75]    (10.93,-3.29) .. controls (6.95,-1.4) and (3.31,-0.3) .. (0,0) .. controls (3.31,0.3) and (6.95,1.4) .. (10.93,3.29)   ;

\draw    (234,130) -- (267.93,156.76) ;
\draw [shift={(269.5,158)}, rotate = 218.26] [color={rgb, 255:red, 0; green, 0; blue, 0 }  ][line width=0.75]    (10.93,-3.29) .. controls (6.95,-1.4) and (3.31,-0.3) .. (0,0) .. controls (3.31,0.3) and (6.95,1.4) .. (10.93,3.29)   ;

\draw    (368.5,191) -- (428.5,191) ;
\draw [shift={(430.5,191)}, rotate = 180] [color={rgb, 255:red, 0; green, 0; blue, 0 }  ][line width=0.75]    (10.93,-3.29) .. controls (6.95,-1.4) and (3.31,-0.3) .. (0,0) .. controls (3.31,0.3) and (6.95,1.4) .. (10.93,3.29)   ;

\draw    (370,70) -- (427.5,70) ;
\draw [shift={(429.5,70)}, rotate = 180] [color={rgb, 255:red, 0; green, 0; blue, 0 }  ][line width=0.75]    (10.93,-3.29) .. controls (6.95,-1.4) and (3.31,-0.3) .. (0,0) .. controls (3.31,0.3) and (6.95,1.4) .. (10.93,3.29)   ;

\draw [color={rgb, 255:red, 74; green, 144; blue, 226 }  ,draw opacity=1 ][line width=2.25]    (508.18,56.25) -- (508.18,104.31) ;

\draw [color={rgb, 255:red, 144; green, 19; blue, 254 }  ,draw opacity=1 ][line width=2.25]    (454.17,56.25) -- (508.18,56.25) ;

\draw [color={rgb, 255:red, 74; green, 144; blue, 226 }  ,draw opacity=1 ][line width=2.25]    (453.68,104.79) -- (508.18,104.31) ;

\draw [color={rgb, 255:red, 74; green, 144; blue, 226 }  ,draw opacity=1 ][line width=2.25]    (477.96,35.58) -- (454.17,56.25) ;

\draw [color={rgb, 255:red, 144; green, 19; blue, 254 }  ,draw opacity=1 ][line width=2.25]    (508.18,56.25) -- (526.51,35.58) ;

\draw [color={rgb, 255:red, 74; green, 144; blue, 226 }  ,draw opacity=1 ][line width=2.25]    (477.96,35.58) -- (526.51,35.58) ;

\draw [color={rgb, 255:red, 144; green, 19; blue, 254 }  ,draw opacity=1 ][line width=2.25]    (526.51,35.58) -- (526.51,87.97) ;

\draw [color={rgb, 255:red, 74; green, 144; blue, 226 }  ,draw opacity=1 ][line width=2.25]    (526.51,87.97) -- (508.18,104.31) ;

\draw [color={rgb, 255:red, 74; green, 144; blue, 226 }  ,draw opacity=1 ][line width=2.25]    (477.96,35.58) -- (477.96,87.97) ;

\draw [color={rgb, 255:red, 144; green, 19; blue, 254 }  ,draw opacity=1 ][line width=2.25]    (477.96,87.97) -- (453.68,104.79) ;

\draw [color={rgb, 255:red, 144; green, 19; blue, 254 }  ,draw opacity=1 ][line width=2.25]    (477.96,87.97) -- (526.51,87.97) ;

\draw [color={rgb, 255:red, 144; green, 19; blue, 254 }  ,draw opacity=1 ][line width=2.25]    (454.17,56.25) -- (453.68,104.79) ;

\draw [color={rgb, 255:red, 74; green, 144; blue, 226 }  ,draw opacity=1 ][line width=2.25]    (506.18,181.25) -- (506.18,229.31) ;

\draw [color={rgb, 255:red, 144; green, 19; blue, 254 }  ,draw opacity=1 ][line width=2.25]    (452.17,181.25) -- (506.18,181.25) ;

\draw [color={rgb, 255:red, 74; green, 144; blue, 226 }  ,draw opacity=1 ][line width=2.25]    (451.68,229.79) -- (506.18,229.31) ;

\draw [color={rgb, 255:red, 74; green, 144; blue, 226 }  ,draw opacity=1 ][line width=2.25]    (475.96,160.58) -- (452.17,181.25) ;

\draw [color={rgb, 255:red, 144; green, 19; blue, 254 }  ,draw opacity=1 ][line width=2.25]    (506.18,181.25) -- (524.51,160.58) ;

\draw [color={rgb, 255:red, 74; green, 144; blue, 226 }  ,draw opacity=1 ][line width=2.25]    (475.96,160.58) -- (524.51,160.58) ;

\draw [color={rgb, 255:red, 74; green, 144; blue, 226 }  ,draw opacity=1 ][line width=2.25]    (524.51,160.58) -- (524.51,212.97) ;

\draw [color={rgb, 255:red, 144; green, 19; blue, 254 }  ,draw opacity=1 ][line width=2.25]    (524.51,212.97) -- (506.18,229.31) ;

\draw [color={rgb, 255:red, 144; green, 19; blue, 254 }  ,draw opacity=1 ][line width=2.25]    (475.96,160.58) -- (475.96,212.97) ;

\draw [color={rgb, 255:red, 74; green, 144; blue, 226 }  ,draw opacity=1 ][line width=2.25]    (475.96,212.97) -- (451.68,229.79) ;

\draw [color={rgb, 255:red, 144; green, 19; blue, 254 }  ,draw opacity=1 ][line width=2.25]    (475.96,212.97) -- (524.51,212.97) ;

\draw [color={rgb, 255:red, 144; green, 19; blue, 254 }  ,draw opacity=1 ][line width=2.25]    (452.17,181.25) -- (451.68,229.79) ;

\draw (127,136) node  [align=left] {OR};
\draw (565,65) node  [align=left] {FD (2)};
\draw (562,194) node  [align=left] {FD (1)};
\draw (68,195) node  [align=left] {No rotation};
\draw (196,61) node  [align=left] {$\displaystyle \dfrac{\pi }{2}$ rotation};
\draw (63,64) node  [align=left] {FD (2)};

\end{tikzpicture}

    \caption{Putting adjacent edges on adjacent faces.}
    \label{fig:my_label}
\end{figure}
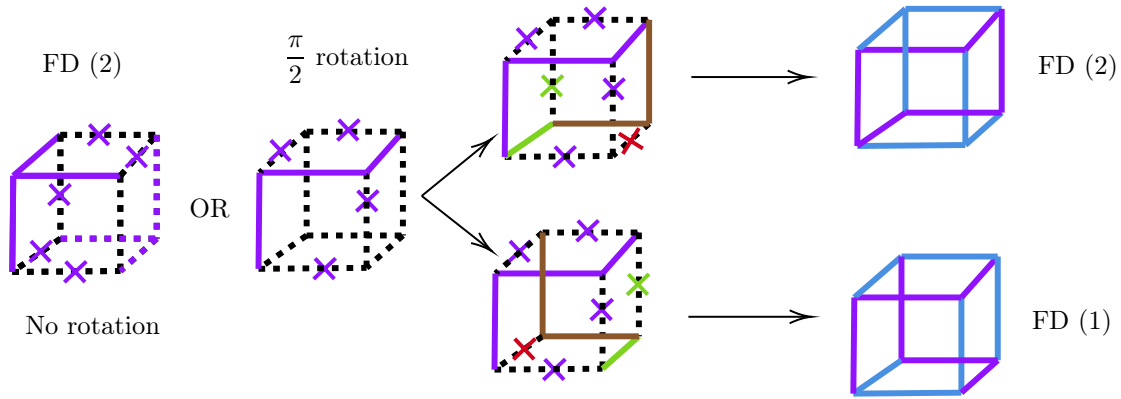

In the "no rotation" case, we end up with FD (2); and in the $\frac{\pi}{2}$ case, we can either end up with FD (1) or FD (2), as shown in Figure 10..

\textbf{Opposite edges, Adjacent Faces}: We need only consider the cases of no rotation and $\dfrac{\pi}{2}$ because the other two cases look the same and therefore have the same problems.

\begin{figure}[H]
    \centering
   
\tikzset{every picture/.style={line width=0.75pt}} 

\begin{tikzpicture}[x=0.6pt,y=0.6pt,yscale=-1,xscale=1]

\draw [color={rgb, 255:red, 144; green, 19; blue, 254 }  ,draw opacity=1 ][line width=2.25]    (192.97,103.98) -- (192.97,164.4) ;

\draw [color={rgb, 255:red, 144; green, 19; blue, 254 }  ,draw opacity=1 ][line width=2.25]    (125.12,103.98) -- (192.97,103.98) ;

\draw [color={rgb, 255:red, 0; green, 0; blue, 0 }  ,draw opacity=1 ][line width=2.25]  [dash pattern={on 2.53pt off 3.02pt}]  (124.5,165) -- (192.97,164.4) ;

\draw [color={rgb, 255:red, 0; green, 0; blue, 0 }  ,draw opacity=1 ][line width=2.25]  [dash pattern={on 2.53pt off 3.02pt}]  (155,78) -- (125.12,103.98) ;

\draw [color={rgb, 255:red, 0; green, 0; blue, 0 }  ,draw opacity=1 ][line width=2.25]  [dash pattern={on 2.53pt off 3.02pt}]  (192.97,103.98) -- (216,78) ;

\draw [color={rgb, 255:red, 144; green, 19; blue, 254 }  ,draw opacity=1 ][line width=2.25]    (155,78) -- (216,78) ;

\draw [color={rgb, 255:red, 0; green, 0; blue, 0 }  ,draw opacity=1 ][line width=2.25]  [dash pattern={on 2.53pt off 3.02pt}]  (216,78) -- (216,143.85) ;

\draw [color={rgb, 255:red, 0; green, 0; blue, 0 }  ,draw opacity=1 ][line width=2.25]  [dash pattern={on 2.53pt off 3.02pt}]  (216,143.85) -- (192.97,164.4) ;

\draw [color={rgb, 255:red, 0; green, 0; blue, 0 }  ,draw opacity=1 ][line width=2.25]  [dash pattern={on 2.53pt off 3.02pt}]  (155,78) -- (155,143.85) ;

\draw [color={rgb, 255:red, 0; green, 0; blue, 0 }  ,draw opacity=1 ][line width=2.25]  [dash pattern={on 2.53pt off 3.02pt}]  (155,143.85) -- (124.5,165) ;

\draw [color={rgb, 255:red, 0; green, 0; blue, 0 }  ,draw opacity=1 ][line width=2.25]  [dash pattern={on 2.53pt off 3.02pt}]  (155,143.85) -- (216,143.85) ;

\draw [color={rgb, 255:red, 144; green, 19; blue, 254 }  ,draw opacity=1 ][line width=2.25]    (125.12,103.98) -- (124.5,165) ;

\draw  [color={rgb, 255:red, 208; green, 2; blue, 27 }  ,draw opacity=0.48 ][line width=2.25]  (83,123) .. controls (83,75.5) and (121.62,37) .. (169.25,37) .. controls (216.88,37) and (255.5,75.5) .. (255.5,123) .. controls (255.5,170.5) and (216.88,209) .. (169.25,209) .. controls (121.62,209) and (83,170.5) .. (83,123) -- cycle ; \draw  [color={rgb, 255:red, 208; green, 2; blue, 27 }  ,draw opacity=0.48 ][line width=2.25]  (108.26,62.19) -- (230.24,183.81) ; \draw  [color={rgb, 255:red, 208; green, 2; blue, 27 }  ,draw opacity=0.48 ][line width=2.25]  (230.24,62.19) -- (108.26,183.81) ;
\draw [color={rgb, 255:red, 144; green, 19; blue, 254 }  ,draw opacity=1 ][line width=2.25]    (445.97,105.98) -- (445.97,166.4) ;

\draw [color={rgb, 255:red, 0; green, 0; blue, 0 }  ,draw opacity=1 ][line width=2.25]  [dash pattern={on 2.53pt off 3.02pt}]  (378.12,105.98) -- (445.97,105.98) ;

\draw [color={rgb, 255:red, 0; green, 0; blue, 0 }  ,draw opacity=1 ][line width=2.25]  [dash pattern={on 2.53pt off 3.02pt}]  (377.5,167) -- (445.97,166.4) ;

\draw [color={rgb, 255:red, 144; green, 19; blue, 254 }  ,draw opacity=1 ][line width=2.25]    (408,80) -- (378.12,105.98) ;

\draw [color={rgb, 255:red, 144; green, 19; blue, 254 }  ,draw opacity=1 ][line width=2.25]    (445.97,105.98) -- (469,80) ;

\draw [color={rgb, 255:red, 0; green, 0; blue, 0 }  ,draw opacity=1 ][line width=2.25]  [dash pattern={on 2.53pt off 3.02pt}]  (408,80) -- (469,80) ;

\draw [color={rgb, 255:red, 0; green, 0; blue, 0 }  ,draw opacity=1 ][line width=2.25]  [dash pattern={on 2.53pt off 3.02pt}]  (469,80) -- (469,145.85) ;

\draw [color={rgb, 255:red, 0; green, 0; blue, 0 }  ,draw opacity=1 ][line width=2.25]  [dash pattern={on 2.53pt off 3.02pt}]  (469,145.85) -- (445.97,166.4) ;

\draw [color={rgb, 255:red, 0; green, 0; blue, 0 }  ,draw opacity=1 ][line width=2.25]  [dash pattern={on 2.53pt off 3.02pt}]  (408,80) -- (408,145.85) ;

\draw [color={rgb, 255:red, 0; green, 0; blue, 0 }  ,draw opacity=1 ][line width=2.25]  [dash pattern={on 2.53pt off 3.02pt}]  (408,145.85) -- (377.5,167) ;

\draw [color={rgb, 255:red, 0; green, 0; blue, 0 }  ,draw opacity=1 ][line width=2.25]  [dash pattern={on 2.53pt off 3.02pt}]  (408,145.85) -- (469,145.85) ;

\draw [color={rgb, 255:red, 144; green, 19; blue, 254 }  ,draw opacity=1 ][line width=2.25]    (378.12,105.98) -- (377.5,167) ;

\draw [rotate around= { 46.67: (416.3, 167.32)
    }] [color={rgb, 255:red, 144; green, 19; blue, 254 }  ,draw opacity=1 ][line width=1.5]  (409.2,167.32) -- (423.4,167.32)(416.3,159.69) -- (416.3,174.94) ;
\draw [rotate around= { 46.67: (458.3, 157.32)
    }] [color={rgb, 255:red, 144; green, 19; blue, 254 }  ,draw opacity=1 ][line width=1.5]  (451.2,157.32) -- (465.4,157.32)(458.3,149.69) -- (458.3,164.94) ;
\draw [rotate around= { 46.67: (468.3, 117.32)
    }] [color={rgb, 255:red, 144; green, 19; blue, 254 }  ,draw opacity=1 ][line width=1.5]  (461.2,117.32) -- (475.4,117.32)(468.3,109.69) -- (468.3,124.94) ;
\draw [rotate around= { 46.67: (431.3, 79.32)
    }] [color={rgb, 255:red, 144; green, 19; blue, 254 }  ,draw opacity=1 ][line width=1.5]  (424.2,79.32) -- (438.4,79.32)(431.3,71.69) -- (431.3,86.94) ;
\draw [rotate around= { 46.67: (423.3, 107.32)
    }] [color={rgb, 255:red, 144; green, 19; blue, 254 }  ,draw opacity=1 ][line width=1.5]  (416.2,107.32) -- (430.4,107.32)(423.3,99.69) -- (423.3,114.94) ;
\draw [rotate around= { 46.67: (407.3, 118.32)
    }] [color={rgb, 255:red, 144; green, 19; blue, 254 }  ,draw opacity=1 ][line width=1.5]  (400.2,118.32) -- (414.4,118.32)(407.3,110.69) -- (407.3,125.94) ;
\draw [rotate around= { 46.67: (395.3, 153.32)
    }] [color={rgb, 255:red, 144; green, 19; blue, 254 }  ,draw opacity=1 ][line width=1.5]  (388.2,153.32) -- (402.4,153.32)(395.3,145.69) -- (395.3,160.94) ;
\draw    (435.13, 145.85) circle [x radius= 33.88, y radius= 9.68]  ;
\draw  [color={rgb, 255:red, 208; green, 2; blue, 27 }  ,draw opacity=0.48 ][line width=2.25]  (337,120) .. controls (337,72.5) and (375.62,34) .. (423.25,34) .. controls (470.88,34) and (509.5,72.5) .. (509.5,120) .. controls (509.5,167.5) and (470.88,206) .. (423.25,206) .. controls (375.62,206) and (337,167.5) .. (337,120) -- cycle ; \draw  [color={rgb, 255:red, 208; green, 2; blue, 27 }  ,draw opacity=0.48 ][line width=2.25]  (362.26,59.19) -- (484.24,180.81) ; \draw  [color={rgb, 255:red, 208; green, 2; blue, 27 }  ,draw opacity=0.48 ][line width=2.25]  (484.24,59.19) -- (362.26,180.81) ;

\draw (175,230) node  [align=left] {Violates 2 edges/face};
\draw (300,116) node  [align=left] {OR};
\draw (427,235) node  [align=left] {Only one edge remains,\\but we need two more.};
\draw (174,20) node  [align=left] {$\displaystyle \pi $ rotation};
\draw (428,17) node  [align=left] {No rotation};

\end{tikzpicture}

    \caption{Putting two opposite edges on two adjacent faces.}
    \label{fig:my_label}
\end{figure}
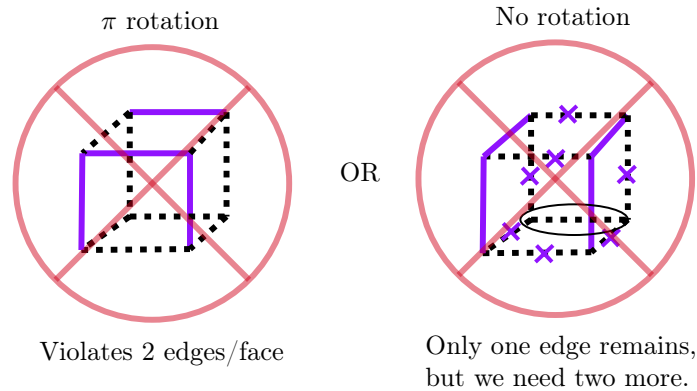

Here, we cannot even end up with a viable fundamental domain.

We conclude that FD(1) and FD(2) are the only possible fundamental domains on the cube with two edges per face and torsion-free groups in the 6-6 breakdown.

\textbf{Case 2}: Three edges per face

We will actually show that there are no fundamental domains on the cube with torsion free groups that have three edges per face. We begin with a face with three edges and consider both the opposite and adjacent cases. Up to rotation, there is only one way to put three edges on a face.

\textbf{Opposite faces:}
We will consider only the cases of a $\pi$ rotation an a $\frac{\pi}{2}$ rotation; the 0 rotation leaves the bottom face completely in the other edge class, and the $\frac{3\pi}{2}$ rotation is analogous to the $\frac{\pi}{2}$ rotation.

\begin{figure}[H]
    \centering
\tikzset{every picture/.style={line width=0.75pt}} 

\begin{tikzpicture}[x=0.75pt,y=0.75pt,yscale=-1,xscale=1]

\draw [color={rgb, 255:red, 0; green, 0; blue, 0 }  ,draw opacity=1 ][line width=2.25]  [dash pattern={on 2.53pt off 3.02pt}]  (335.97,37.98) -- (335.97,98.4) ;

\draw [color={rgb, 255:red, 0; green, 0; blue, 0 }  ,draw opacity=1 ][line width=2.25]  [dash pattern={on 2.53pt off 3.02pt}]  (268.12,37.98) -- (335.97,37.98) ;

\draw [color={rgb, 255:red, 0; green, 0; blue, 0 }  ,draw opacity=1 ][line width=2.25]  [dash pattern={on 2.53pt off 3.02pt}]  (267.5,99) -- (335.97,98.4) ;

\draw [color={rgb, 255:red, 144; green, 19; blue, 254 }  ,draw opacity=1 ][line width=2.25]    (298,12) -- (268.12,37.98) ;

\draw [color={rgb, 255:red, 0; green, 0; blue, 0 }  ,draw opacity=1 ][line width=2.25]  [dash pattern={on 2.53pt off 3.02pt}]  (335.97,37.98) -- (359,12) ;

\draw [color={rgb, 255:red, 0; green, 0; blue, 0 }  ,draw opacity=1 ][line width=2.25]  [dash pattern={on 2.53pt off 3.02pt}]  (298,12) -- (359,12) ;

\draw [color={rgb, 255:red, 0; green, 0; blue, 0 }  ,draw opacity=1 ][line width=2.25]  [dash pattern={on 2.53pt off 3.02pt}]  (359,12) -- (359,77.85) ;

\draw [color={rgb, 255:red, 0; green, 0; blue, 0 }  ,draw opacity=1 ][line width=2.25]  [dash pattern={on 2.53pt off 3.02pt}]  (359,77.85) -- (335.97,98.4) ;

\draw [color={rgb, 255:red, 144; green, 19; blue, 254 }  ,draw opacity=1 ][line width=2.25]    (298,12) -- (298,77.85) ;

\draw [color={rgb, 255:red, 0; green, 0; blue, 0 }  ,draw opacity=1 ][line width=2.25]  [dash pattern={on 2.53pt off 3.02pt}]  (298,77.85) -- (267.5,99) ;

\draw [color={rgb, 255:red, 0; green, 0; blue, 0 }  ,draw opacity=1 ][line width=2.25]  [dash pattern={on 2.53pt off 3.02pt}]  (298,77.85) -- (359,77.85) ;

\draw [color={rgb, 255:red, 144; green, 19; blue, 254 }  ,draw opacity=1 ][line width=2.25]    (268.12,37.98) -- (267.5,99) ;

\draw [color={rgb, 255:red, 144; green, 19; blue, 254 }  ,draw opacity=1 ][line width=2.25]    (239.97,148.98) -- (239.97,209.4) ;

\draw [color={rgb, 255:red, 74; green, 144; blue, 226 }  ,draw opacity=1 ][line width=2.25]    (172.12,148.98) -- (239.97,148.98) ;

\draw [color={rgb, 255:red, 74; green, 144; blue, 226 }  ,draw opacity=1 ][line width=2.25]    (171.5,210) -- (239.97,209.4) ;

\draw [color={rgb, 255:red, 144; green, 19; blue, 254 }  ,draw opacity=1 ][line width=2.25]    (202,123) -- (172.12,148.98) ;

\draw [color={rgb, 255:red, 74; green, 144; blue, 226 }  ,draw opacity=1 ][line width=2.25]    (239.97,148.98) -- (263,123) ;

\draw [color={rgb, 255:red, 74; green, 144; blue, 226 }  ,draw opacity=1 ][line width=2.25]    (202,123) -- (263,123) ;

\draw [color={rgb, 255:red, 144; green, 19; blue, 254 }  ,draw opacity=1 ][line width=2.25]    (263,123) -- (263,188.85) ;

\draw [color={rgb, 255:red, 144; green, 19; blue, 254 }  ,draw opacity=1 ][line width=2.25]    (263,188.85) -- (239.97,209.4) ;

\draw [color={rgb, 255:red, 144; green, 19; blue, 254 }  ,draw opacity=1 ][line width=2.25]    (202,123) -- (202,188.85) ;

\draw [color={rgb, 255:red, 74; green, 144; blue, 226 }  ,draw opacity=1 ][line width=2.25]    (202,188.85) -- (171.5,210) ;

\draw [color={rgb, 255:red, 74; green, 144; blue, 226 }  ,draw opacity=1 ][line width=2.25]    (202,188.85) -- (263,188.85) ;

\draw [color={rgb, 255:red, 144; green, 19; blue, 254 }  ,draw opacity=1 ][line width=2.25]    (172.12,148.98) -- (171.5,210) ;

\draw [color={rgb, 255:red, 144; green, 19; blue, 254 }  ,draw opacity=1 ][line width=2.25]    (413.97,146.98) -- (413.97,207.4) ;

\draw [color={rgb, 255:red, 74; green, 144; blue, 226 }  ,draw opacity=1 ][line width=2.25]    (346.12,146.98) -- (413.97,146.98) ;

\draw [color={rgb, 255:red, 74; green, 144; blue, 226 }  ,draw opacity=1 ][line width=2.25]    (345.5,208) -- (413.97,207.4) ;

\draw [color={rgb, 255:red, 144; green, 19; blue, 254 }  ,draw opacity=1 ][line width=2.25]    (376,121) -- (346.12,146.98) ;

\draw [color={rgb, 255:red, 144; green, 19; blue, 254 }  ,draw opacity=1 ][line width=2.25]    (413.97,146.98) -- (437,121) ;

\draw [color={rgb, 255:red, 74; green, 144; blue, 226 }  ,draw opacity=1 ][line width=2.25]    (376,121) -- (437,121) ;

\draw [color={rgb, 255:red, 74; green, 144; blue, 226 }  ,draw opacity=1 ][line width=2.25]    (437,121) -- (437,186.85) ;

\draw [color={rgb, 255:red, 144; green, 19; blue, 254 }  ,draw opacity=1 ][line width=2.25]    (437,186.85) -- (413.97,207.4) ;

\draw [color={rgb, 255:red, 144; green, 19; blue, 254 }  ,draw opacity=1 ][line width=2.25]    (376,121) -- (376,186.85) ;

\draw [color={rgb, 255:red, 74; green, 144; blue, 226 }  ,draw opacity=1 ][line width=2.25]    (376,186.85) -- (345.5,208) ;

\draw [color={rgb, 255:red, 74; green, 144; blue, 226 }  ,draw opacity=1 ][line width=2.25]    (376,186.85) -- (437,186.85) ;

\draw [color={rgb, 255:red, 144; green, 19; blue, 254 }  ,draw opacity=1 ][line width=2.25]    (346.12,146.98) -- (345.5,208) ;

\draw    (380.05, 146.98) circle [x radius= 38.5, y radius= 9.74]  ;
\draw    (348.55,101.98) -- (365.16,120.51) ;
\draw [shift={(366.5,122)}, rotate = 228.12] [color={rgb, 255:red, 0; green, 0; blue, 0 }  ][line width=0.75]    (10.93,-3.29) .. controls (6.95,-1.4) and (3.31,-0.3) .. (0,0) .. controls (3.31,0.3) and (6.95,1.4) .. (10.93,3.29)   ;

\draw    (281.47,106.98) -- (270.36,130.2) ;
\draw [shift={(269.5,132)}, rotate = 295.57] [color={rgb, 255:red, 0; green, 0; blue, 0 }  ][line width=0.75]    (10.93,-3.29) .. controls (6.95,-1.4) and (3.31,-0.3) .. (0,0) .. controls (3.31,0.3) and (6.95,1.4) .. (10.93,3.29)   ;

\draw (471,99) node  [align=left] {Requires an elliptic\\over circled edge};
\draw (149,106) node  [align=left] {With allowed face\\ IDs, more than 2 \\edge classes};
\draw (210,227) node  [align=left] {$\displaystyle \pi $ rotation};
\draw (391,232) node  [align=left] {$\displaystyle \dfrac{\pi }{2}$ rotation};

\end{tikzpicture}

    \caption{Opposite faces with three edges each.}
    \label{fig:my_label}
\end{figure}
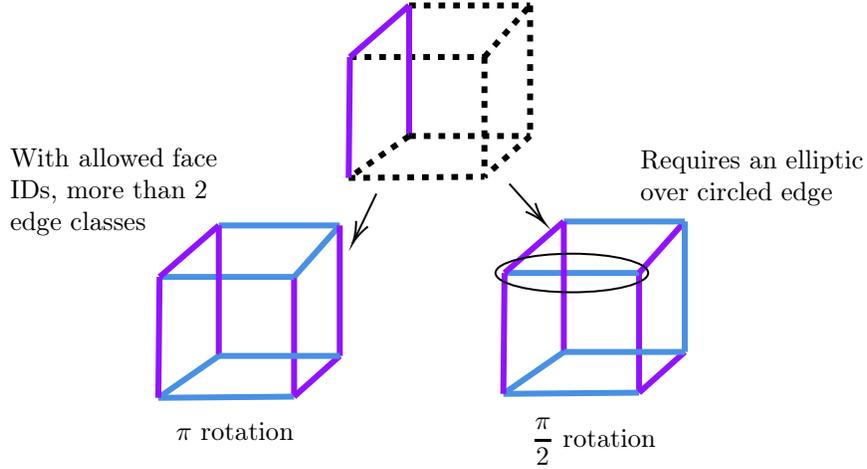

 For the left case, although it appears that there are two edge classes of six, there are actually more distinct edge classes once one carries out the face identifications. The general issue is that we know that top and bottom, left and right, and front and back have to be paired because they are combinatorially equivalent. Additionally, they all have to be paired with either no twist or a twist of $\pi$, and with those ientifications, there is no way that edges on faces separated by a $\dfrac{\pi}{2}$ rotation can be in the same class, yet they need to be in this figure. Thus, the opposite faces case does not produce a viable fundamental domain.
 
\textbf{Adjacent Faces}:

To simplify the argument, we will note that both the 0 and $\pi$ rotation cases require an elliptic generator, and the $\dfrac{\pi}{2}$ and $\dfrac{3\pi}{2}$ are analogous. The $\dfrac{\pi}{2}$ case is shown in Figure 14. The left and top faces have the original five edges, and there is a unique choice for the $6^{\mbox{th}}$ edge in order to allow for face identifications; if it is put anywhere else, there will be a face that is not combinatorially equivalent to any other faces (or an odd number of combinatorially equivalent faces).

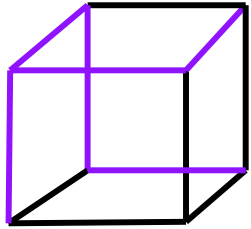
\begin{figure}[H]
    \centering
\tikzset{every picture/.style={line width=0.75pt}} 

\begin{tikzpicture}[x=0.65pt,y=0.65pt,yscale=-1,xscale=1]

\draw [color={rgb, 255:red, 0; green, 0; blue, 0 }  ,draw opacity=1 ][line width=2.25]    (341.77,71.92) -- (341.77,160.12) ;

\draw [color={rgb, 255:red, 144; green, 19; blue, 254 }  ,draw opacity=1 ][line width=2.25]    (239.44,71.92) -- (341.77,71.92) ;

\draw [color={rgb, 255:red, 0; green, 0; blue, 0 }  ,draw opacity=1 ][line width=2.25]    (238.5,161) -- (341.77,160.12) ;

\draw [color={rgb, 255:red, 144; green, 19; blue, 254 }  ,draw opacity=1 ][line width=2.25]    (284.5,34) -- (239.44,71.92) ;

\draw [color={rgb, 255:red, 144; green, 19; blue, 254 }  ,draw opacity=1 ][line width=2.25]    (341.77,71.92) -- (376.5,34) ;

\draw [color={rgb, 255:red, 0; green, 0; blue, 0 }  ,draw opacity=1 ][line width=2.25]    (284.5,34) -- (376.5,34) ;

\draw [color={rgb, 255:red, 0; green, 0; blue, 0 }  ,draw opacity=1 ][line width=2.25]    (376.5,34) -- (376.5,130.13) ;

\draw [color={rgb, 255:red, 0; green, 0; blue, 0 }  ,draw opacity=1 ][line width=2.25]    (376.5,130.13) -- (341.77,160.12) ;

\draw [color={rgb, 255:red, 144; green, 19; blue, 254 }  ,draw opacity=1 ][line width=2.25]    (284.5,34) -- (284.5,130.13) ;

\draw [color={rgb, 255:red, 0; green, 0; blue, 0 }  ,draw opacity=1 ][line width=2.25]    (284.5,130.13) -- (238.5,161) ;

\draw [color={rgb, 255:red, 144; green, 19; blue, 254 }  ,draw opacity=1 ][line width=2.25]    (284.5,130.13) -- (376.5,130.13) ;

\draw [color={rgb, 255:red, 144; green, 19; blue, 254 }  ,draw opacity=1 ][line width=2.25]    (239.44,71.92) -- (238.5,161) ;

\end{tikzpicture}

    \caption{Adjacent faces with three edges each and sharing one edge. This is the only case we need to consider.}
    \label{fig:my_label}
\end{figure}

It is clear that the face identifications have to pair top and left, right and bottom, and front and back because those are the pairs of faces that are combinatorially equivalent. These pairings give rise to the following face identifications:
\begin{itemize}
    \item Top$\longrightarrow$Left with $\dfrac{\pi}{2}$ clockwise twist
    \item Right$\longrightarrow$Bottom with $\frac{\pi}{2}$ clockwise twist
    \item Front$\longrightarrow$Back with $\frac{\pi}{2}$ counterclockwise twist
\end{itemize}
However, when these face identifications are carried out, the purple class is actually broken into two edge classes, each containing three edges, and we do not have a fundamental domain.

We can now conclude that there are no fundamental domains on the cube with non-elliptic generators that have 3 edges per face.

 All that is left for the 6-6 edge breakdown is to justify that two edges per face and three edges per face are the only possible cases. To see this, note that we certainly cannot have 4 edges on a face. If we have three edges on a face, there must be another face with three edges in order to have face identifications, and we considered all of those cases; similarly if we start with two on a face. Thus the only remaining case is that we have one edge on every face, which is not possible on the cube.
 \end{proof}

To complete the proof of Proposition 5.1, we need only consider the 5-7 edge breakdown.

\noindent\textbf{Claim 3}: There is only one fundamental domain on the cube with non-elliptic generators and edge classes of sizes 5 and 7.

First, I will name two "rules" that will help rule out cases for the 5-7 breakdown, shown in Figure 14.
\begin{itemize}
    \item Rule 1: No "4 pillars" in an edge class (Figure 14(a))
   
We only have one more edge to place, so that necessarily leaves a whole face to the other edge class, which we cannot have.

\item Rule 2: Opposite vertices cannot be left open in the 7 edge class; if they are, then 6 of the edges will sum to $4\pi$, but the whole class of 7 has to sum to $5\pi$, which means that the remaining vertex has an exterior dihedral angle of $\pi$, a contradiction (Figure 14(b)).

\end{itemize}
\begin{figure}[H]
    \centering

\tikzset{every picture/.style={line width=0.75pt}} 

\begin{tikzpicture}[x=0.75pt,y=0.75pt,yscale=-1,xscale=1]

\draw [color={rgb, 255:red, 74; green, 144; blue, 226 }  ,draw opacity=1 ][line width=2.25]    (154.83,64) -- (154.83,163.16) ;

\draw [color={rgb, 255:red, 74; green, 144; blue, 226 }  ,draw opacity=1 ][line width=2.25]    (210.02,103.12) -- (210.02,194.09) ;

\draw [color={rgb, 255:red, 0; green, 0; blue, 0 }  ,draw opacity=1 ][line width=2.25]    (111.4,103.12) -- (210.02,103.12) ;

\draw [color={rgb, 255:red, 0; green, 0; blue, 0 }  ,draw opacity=1 ][line width=2.25]    (110.5,195) -- (210.02,194.09) ;

\draw [color={rgb, 255:red, 0; green, 0; blue, 0 }  ,draw opacity=1 ][line width=2.25]    (154.83,64) -- (111.4,103.12) ;

\draw [color={rgb, 255:red, 0; green, 0; blue, 0 }  ,draw opacity=1 ][line width=2.25]    (210.02,103.12) -- (243.5,64) ;

\draw [color={rgb, 255:red, 0; green, 0; blue, 0 }  ,draw opacity=1 ][line width=2.25]    (154.83,64) -- (243.5,64) ;

\draw [color={rgb, 255:red, 74; green, 144; blue, 226 }  ,draw opacity=1 ][line width=2.25]    (243.5,64) -- (243.5,163.16) ;

\draw [color={rgb, 255:red, 0; green, 0; blue, 0 }  ,draw opacity=1 ][line width=2.25]    (243.5,163.16) -- (210.02,194.09) ;

\draw [color={rgb, 255:red, 0; green, 0; blue, 0 }  ,draw opacity=1 ][line width=2.25]    (154.83,163.16) -- (110.5,195) ;

\draw [color={rgb, 255:red, 0; green, 0; blue, 0 }  ,draw opacity=1 ][line width=2.25]    (154.83,163.16) -- (243.5,163.16) ;

\draw [color={rgb, 255:red, 74; green, 144; blue, 226 }  ,draw opacity=1 ][line width=2.25]    (111.4,103.12) -- (110.5,195) ;

\draw [color={rgb, 255:red, 74; green, 144; blue, 226 }  ,draw opacity=1 ][line width=2.25]    (365.83,63) -- (365.83,162.16) ;

\draw [color={rgb, 255:red, 74; green, 144; blue, 226 }  ,draw opacity=1 ][line width=2.25]    (322.4,102.12) -- (421.02,102.12) ;

\draw [color={rgb, 255:red, 0; green, 0; blue, 0 }  ,draw opacity=1 ][line width=2.25]    (321.5,194) -- (421.02,193.09) ;

\draw [color={rgb, 255:red, 0; green, 0; blue, 0 }  ,draw opacity=1 ][line width=2.25]    (365.83,63) -- (322.4,102.12) ;

\draw [color={rgb, 255:red, 74; green, 144; blue, 226 }  ,draw opacity=1 ][line width=2.25]    (421.02,102.12) -- (454.5,63) ;

\draw [color={rgb, 255:red, 0; green, 0; blue, 0 }  ,draw opacity=1 ][line width=2.25]    (365.83,63) -- (454.5,63) ;

\draw [color={rgb, 255:red, 0; green, 0; blue, 0 }  ,draw opacity=1 ][line width=2.25]    (454.5,63) -- (454.5,162.16) ;

\draw [color={rgb, 255:red, 0; green, 0; blue, 0 }  ,draw opacity=1 ][line width=2.25]    (454.5,162.16) -- (421.02,193.09) ;

\draw [color={rgb, 255:red, 74; green, 144; blue, 226 }  ,draw opacity=1 ][line width=2.25]    (365.83,162.16) -- (321.5,194) ;

\draw [color={rgb, 255:red, 74; green, 144; blue, 226 }  ,draw opacity=1 ][line width=2.25]    (365.83,162.16) -- (454.5,162.16) ;

\draw [color={rgb, 255:red, 0; green, 0; blue, 0 }  ,draw opacity=1 ][line width=2.25]    (322.4,102.12) -- (321.5,194) ;

\draw [color={rgb, 255:red, 74; green, 144; blue, 226 }  ,draw opacity=1 ][line width=2.25]    (421.02,102.12) -- (421.02,193.09) ;

\draw (171,229) node  [align=left] {(a)};
\draw (372,229) node  [align=left] {(b)};

\end{tikzpicture}

    \caption{(a) Four pillars (b) Opposite vertices in the same class}
    \label{fig:my_label}
\end{figure}
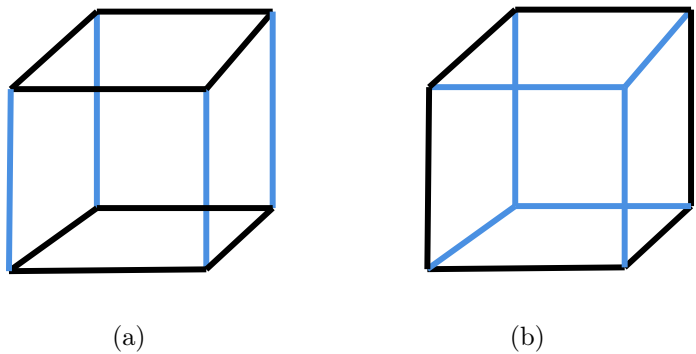
There are two cases to consider: an edge must be shared between adjacent faces in the 5-class since there are an odd number of edges in the class. That edge can either be between faces with three edges each, two edges each, or one edge each.

\textbf{One edge:}

Suppose that we have one edge shared between adjacent faces and no other edges on that face, as in Figure 15.

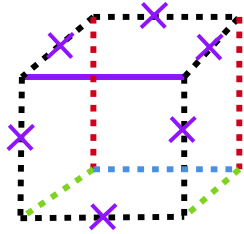
\begin{figure}[H]
    \centering
\tikzset{every picture/.style={line width=0.75pt}} 

\begin{tikzpicture}[x=0.6pt,y=0.6pt,yscale=-1,xscale=1]

\draw [color={rgb, 255:red, 0; green, 0; blue, 0 }  ,draw opacity=1 ][line width=2.25]  [dash pattern={on 2.53pt off 3.02pt}]  (361.77,91.92) -- (361.77,180.12) ;

\draw [color={rgb, 255:red, 144; green, 19; blue, 254 }  ,draw opacity=1 ][line width=2.25]    (259.44,91.92) -- (361.77,91.92) ;

\draw [color={rgb, 255:red, 0; green, 0; blue, 0 }  ,draw opacity=1 ][line width=2.25]  [dash pattern={on 2.53pt off 3.02pt}]  (258.5,181) -- (361.77,180.12) ;

\draw [color={rgb, 255:red, 0; green, 0; blue, 0 }  ,draw opacity=1 ][line width=2.25]  [dash pattern={on 2.53pt off 3.02pt}]  (304.5,54) -- (259.44,91.92) ;

\draw [color={rgb, 255:red, 0; green, 0; blue, 0 }  ,draw opacity=1 ][line width=2.25]  [dash pattern={on 2.53pt off 3.02pt}]  (361.77,91.92) -- (396.5,54) ;

\draw [color={rgb, 255:red, 0; green, 0; blue, 0 }  ,draw opacity=1 ][line width=2.25]  [dash pattern={on 2.53pt off 3.02pt}]  (304.5,54) -- (396.5,54) ;

\draw [color={rgb, 255:red, 208; green, 2; blue, 27 }  ,draw opacity=1 ][line width=2.25]  [dash pattern={on 2.53pt off 3.02pt}]  (396.5,54) -- (396.5,150.13) ;

\draw [color={rgb, 255:red, 126; green, 211; blue, 33 }  ,draw opacity=1 ][line width=2.25]  [dash pattern={on 2.53pt off 3.02pt}]  (396.5,150.13) -- (361.77,180.12) ;

\draw [color={rgb, 255:red, 208; green, 2; blue, 27 }  ,draw opacity=1 ][line width=2.25]  [dash pattern={on 2.53pt off 3.02pt}]  (304.5,54) -- (304.5,150.13) ;

\draw [color={rgb, 255:red, 126; green, 211; blue, 33 }  ,draw opacity=1 ][line width=2.25]  [dash pattern={on 2.53pt off 3.02pt}]  (304.5,150.13) -- (258.5,181) ;

\draw [color={rgb, 255:red, 74; green, 144; blue, 226 }  ,draw opacity=1 ][line width=2.25]  [dash pattern={on 2.53pt off 3.02pt}]  (304.5,150.13) -- (396.5,150.13) ;

\draw [color={rgb, 255:red, 0; green, 0; blue, 0 }  ,draw opacity=1 ][line width=2.25]  [dash pattern={on 2.53pt off 3.02pt}]  (259.44,91.92) -- (258.5,181) ;

\draw [rotate around= { 46.67: (283.68, 72.07)
    }] [color={rgb, 255:red, 144; green, 19; blue, 254 }  ,draw opacity=1 ][line width=1.5]  (273.21,72.07) -- (294.14,72.07)(283.68,60.88) -- (283.68,83.25) ;
\draw [rotate around= { 46.67: (360.68, 125.07)
    }] [color={rgb, 255:red, 144; green, 19; blue, 254 }  ,draw opacity=1 ][line width=1.5]  (350.21,125.07) -- (371.14,125.07)(360.68,113.88) -- (360.68,136.25) ;
\draw [rotate around= { 46.67: (310.68, 180.07)
    }] [color={rgb, 255:red, 144; green, 19; blue, 254 }  ,draw opacity=1 ][line width=1.5]  (300.21,180.07) -- (321.14,180.07)(310.68,168.88) -- (310.68,191.25) ;
\draw [rotate around= { 46.67: (258.68, 130.07)
    }] [color={rgb, 255:red, 144; green, 19; blue, 254 }  ,draw opacity=1 ][line width=1.5]  (248.21,130.07) -- (269.14,130.07)(258.68,118.88) -- (258.68,141.25) ;
\draw [rotate around= { 46.67: (342.68, 53.07)
    }] [color={rgb, 255:red, 144; green, 19; blue, 254 }  ,draw opacity=1 ][line width=1.5]  (332.21,53.07) -- (353.14,53.07)(342.68,41.88) -- (342.68,64.25) ;
\draw [rotate around= { 46.67: (377.68, 73.07)
    }] [color={rgb, 255:red, 144; green, 19; blue, 254 }  ,draw opacity=1 ][line width=1.5]  (367.21,73.07) -- (388.14,73.07)(377.68,61.88) -- (377.68,84.25) ;

\end{tikzpicture}

    \caption{The purple "x"s denote edges that would result in having more than one edge on one of the adjacent faces. None of the remaining edge choices are viable options.}
    \label{fig:my_label}
\end{figure}

We have to pick four of the remaining five available edges to have an edge class of five edges. If we omit the blue edge, we need an elliptic to identify back to bottom. If we omit either green edge, then the back face has three edges in the purple class and no pair. Similarly, if we omit either red edge, then the bottom face has three edges in the purple class and no pair. Thus, no fundamental domains can arise in this case.

\textbf{2 edges}:

Next, we consider the case where the adjacent faces sharing an edge have two edges each. We have to break into two subcases here according to how those two edges are arranged on the faces.

\textbf{Subcase 1}:

The "vertex" case--the edges are arranged around a vertex, as in Figure 16.

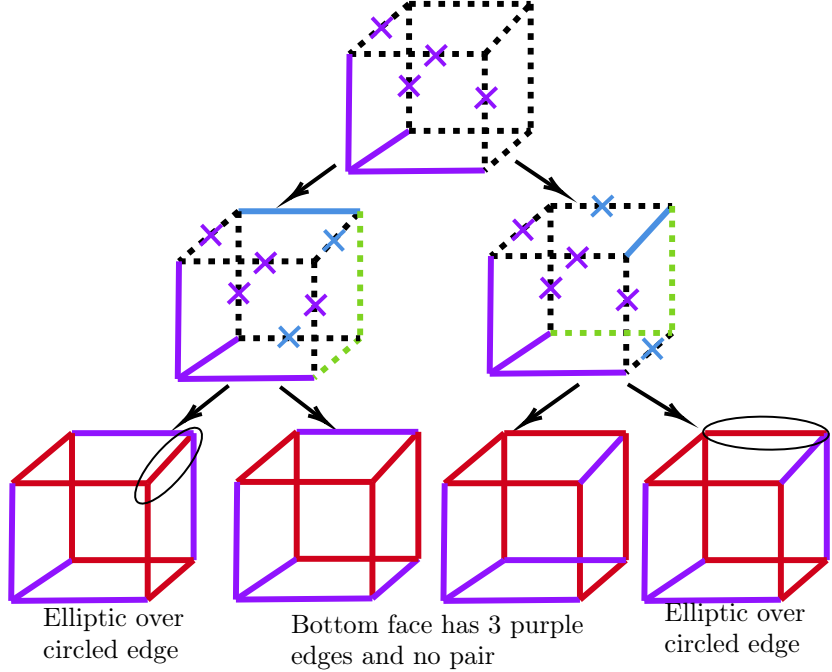
\begin{figure}[H]
    \centering

\tikzset{every picture/.style={line width=0.75pt}} 

\begin{tikzpicture}[x=0.5pt,y=0.5pt,yscale=-1,xscale=1]

\draw [color={rgb, 255:red, 0; green, 0; blue, 0 }  ,draw opacity=1 ][line width=2.25]  [dash pattern={on 2.53pt off 3.02pt}]  (373.77,48.92) -- (373.77,137.12) ;

\draw [color={rgb, 255:red, 0; green, 0; blue, 0 }  ,draw opacity=1 ][line width=2.25]  [dash pattern={on 2.53pt off 3.02pt}]  (271.44,48.92) -- (373.77,48.92) ;

\draw [color={rgb, 255:red, 144; green, 19; blue, 254 }  ,draw opacity=1 ][line width=2.25]    (270.5,138) -- (373.77,137.12) ;

\draw [color={rgb, 255:red, 0; green, 0; blue, 0 }  ,draw opacity=1 ][line width=2.25]  [dash pattern={on 2.53pt off 3.02pt}]  (316.5,11) -- (271.44,48.92) ;

\draw [color={rgb, 255:red, 0; green, 0; blue, 0 }  ,draw opacity=1 ][line width=2.25]  [dash pattern={on 2.53pt off 3.02pt}]  (373.77,48.92) -- (408.5,11) ;

\draw [color={rgb, 255:red, 0; green, 0; blue, 0 }  ,draw opacity=1 ][line width=2.25]  [dash pattern={on 2.53pt off 3.02pt}]  (316.5,11) -- (408.5,11) ;

\draw [color={rgb, 255:red, 0; green, 0; blue, 0 }  ,draw opacity=1 ][line width=2.25]  [dash pattern={on 2.53pt off 3.02pt}]  (408.5,11) -- (408.5,107.13) ;

\draw [color={rgb, 255:red, 0; green, 0; blue, 0 }  ,draw opacity=1 ][line width=2.25]  [dash pattern={on 2.53pt off 3.02pt}]  (408.5,107.13) -- (373.77,137.12) ;

\draw [color={rgb, 255:red, 0; green, 0; blue, 0 }  ,draw opacity=1 ][line width=2.25]  [dash pattern={on 2.53pt off 3.02pt}]  (316.5,11) -- (316.5,107.13) ;

\draw [color={rgb, 255:red, 144; green, 19; blue, 254 }  ,draw opacity=1 ][line width=2.25]    (316.5,107.13) -- (270.5,138) ;

\draw [color={rgb, 255:red, 0; green, 0; blue, 0 }  ,draw opacity=1 ][line width=2.25]  [dash pattern={on 2.53pt off 3.02pt}]  (316.5,107.13) -- (408.5,107.13) ;

\draw [color={rgb, 255:red, 144; green, 19; blue, 254 }  ,draw opacity=1 ][line width=2.25]    (271.44,48.92) -- (270.5,138) ;

\draw [rotate around= { 46.67: (295.68, 29.07)
    }] [color={rgb, 255:red, 144; green, 19; blue, 254 }  ,draw opacity=1 ][line width=1.5]  (285.21,29.07) -- (306.14,29.07)(295.68,17.88) -- (295.68,40.25) ;
\draw [rotate around= { 46.67: (316.68, 73.07)
    }] [color={rgb, 255:red, 144; green, 19; blue, 254 }  ,draw opacity=1 ][line width=1.5]  (306.21,73.07) -- (327.14,73.07)(316.68,61.88) -- (316.68,84.25) ;
\draw [rotate around= { 46.67: (374.68, 82.07)
    }] [color={rgb, 255:red, 144; green, 19; blue, 254 }  ,draw opacity=1 ][line width=1.5]  (364.21,82.07) -- (385.14,82.07)(374.68,70.88) -- (374.68,93.25) ;
\draw [rotate around= { 46.67: (336.68, 50.07)
    }] [color={rgb, 255:red, 144; green, 19; blue, 254 }  ,draw opacity=1 ][line width=1.5]  (326.21,50.07) -- (347.14,50.07)(336.68,38.88) -- (336.68,61.25) ;
\draw [line width=1.5]    (261,133) -- (228.93,156.24) ;
\draw [shift={(226.5,158)}, rotate = 324.07] [color={rgb, 255:red, 0; green, 0; blue, 0 }  ][line width=1.5]    (14.21,-4.28) .. controls (9.04,-1.82) and (4.3,-0.39) .. (0,0) .. controls (4.3,0.39) and (9.04,1.82) .. (14.21,4.28)   ;

\draw [line width=1.5]    (397,130) -- (430.01,152.32) ;
\draw [shift={(432.5,154)}, rotate = 214.06] [color={rgb, 255:red, 0; green, 0; blue, 0 }  ][line width=1.5]    (14.21,-4.28) .. controls (9.04,-1.82) and (4.3,-0.39) .. (0,0) .. controls (4.3,0.39) and (9.04,1.82) .. (14.21,4.28)   ;

\draw [color={rgb, 255:red, 0; green, 0; blue, 0 }  ,draw opacity=1 ][line width=2.25]  [dash pattern={on 2.53pt off 3.02pt}]  (244.77,205.92) -- (244.77,294.12) ;

\draw [color={rgb, 255:red, 0; green, 0; blue, 0 }  ,draw opacity=1 ][line width=2.25]  [dash pattern={on 2.53pt off 3.02pt}]  (142.44,205.92) -- (244.77,205.92) ;

\draw [color={rgb, 255:red, 144; green, 19; blue, 254 }  ,draw opacity=1 ][line width=2.25]    (141.5,295) -- (244.77,294.12) ;

\draw [color={rgb, 255:red, 0; green, 0; blue, 0 }  ,draw opacity=1 ][line width=2.25]  [dash pattern={on 2.53pt off 3.02pt}]  (187.5,168) -- (142.44,205.92) ;

\draw [color={rgb, 255:red, 0; green, 0; blue, 0 }  ,draw opacity=1 ][line width=2.25]  [dash pattern={on 2.53pt off 3.02pt}]  (244.77,205.92) -- (279.5,168) ;

\draw [color={rgb, 255:red, 74; green, 144; blue, 226 }  ,draw opacity=1 ][line width=2.25]    (187.5,168) -- (279.5,168) ;

\draw [color={rgb, 255:red, 126; green, 211; blue, 33 }  ,draw opacity=1 ][line width=2.25]  [dash pattern={on 2.53pt off 3.02pt}]  (279.5,168) -- (279.5,264.13) ;

\draw [color={rgb, 255:red, 126; green, 211; blue, 33 }  ,draw opacity=1 ][line width=2.25]  [dash pattern={on 2.53pt off 3.02pt}]  (279.5,264.13) -- (244.77,294.12) ;

\draw [color={rgb, 255:red, 0; green, 0; blue, 0 }  ,draw opacity=1 ][line width=2.25]  [dash pattern={on 2.53pt off 3.02pt}]  (187.5,168) -- (187.5,264.13) ;

\draw [color={rgb, 255:red, 144; green, 19; blue, 254 }  ,draw opacity=1 ][line width=2.25]    (187.5,264.13) -- (141.5,295) ;

\draw [color={rgb, 255:red, 0; green, 0; blue, 0 }  ,draw opacity=1 ][line width=2.25]  [dash pattern={on 2.53pt off 3.02pt}]  (187.5,264.13) -- (279.5,264.13) ;

\draw [color={rgb, 255:red, 144; green, 19; blue, 254 }  ,draw opacity=1 ][line width=2.25]    (142.44,205.92) -- (141.5,295) ;

\draw [rotate around= { 46.67: (166.68, 186.07)
    }] [color={rgb, 255:red, 144; green, 19; blue, 254 }  ,draw opacity=1 ][line width=1.5]  (156.21,186.07) -- (177.14,186.07)(166.68,174.88) -- (166.68,197.25) ;
\draw [rotate around= { 46.67: (187.68, 230.07)
    }] [color={rgb, 255:red, 144; green, 19; blue, 254 }  ,draw opacity=1 ][line width=1.5]  (177.21,230.07) -- (198.14,230.07)(187.68,218.88) -- (187.68,241.25) ;
\draw [rotate around= { 46.67: (245.68, 239.07)
    }] [color={rgb, 255:red, 144; green, 19; blue, 254 }  ,draw opacity=1 ][line width=1.5]  (235.21,239.07) -- (256.14,239.07)(245.68,227.88) -- (245.68,250.25) ;
\draw [rotate around= { 46.67: (207.68, 207.07)
    }] [color={rgb, 255:red, 144; green, 19; blue, 254 }  ,draw opacity=1 ][line width=1.5]  (197.21,207.07) -- (218.14,207.07)(207.68,195.88) -- (207.68,218.25) ;
\draw [rotate around= { 46.67: (259.68, 190.07)
    }] [color={rgb, 255:red, 74; green, 144; blue, 226 }  ,draw opacity=1 ][line width=1.5]  (249.21,190.07) -- (270.14,190.07)(259.68,178.88) -- (259.68,201.25) ;
\draw [color={rgb, 255:red, 0; green, 0; blue, 0 }  ,draw opacity=1 ][line width=2.25]  [dash pattern={on 2.53pt off 3.02pt}]  (480.77,201.92) -- (480.77,290.12) ;

\draw [color={rgb, 255:red, 0; green, 0; blue, 0 }  ,draw opacity=1 ][line width=2.25]  [dash pattern={on 2.53pt off 3.02pt}]  (378.44,201.92) -- (480.77,201.92) ;

\draw [color={rgb, 255:red, 144; green, 19; blue, 254 }  ,draw opacity=1 ][line width=2.25]    (377.5,291) -- (480.77,290.12) ;

\draw [color={rgb, 255:red, 0; green, 0; blue, 0 }  ,draw opacity=1 ][line width=2.25]  [dash pattern={on 2.53pt off 3.02pt}]  (423.5,164) -- (378.44,201.92) ;

\draw [color={rgb, 255:red, 74; green, 144; blue, 226 }  ,draw opacity=1 ][line width=2.25]    (480.77,201.92) -- (515.5,164) ;

\draw [color={rgb, 255:red, 0; green, 0; blue, 0 }  ,draw opacity=1 ][line width=2.25]  [dash pattern={on 2.53pt off 3.02pt}]  (423.5,164) -- (515.5,164) ;

\draw [color={rgb, 255:red, 126; green, 211; blue, 33 }  ,draw opacity=1 ][line width=2.25]  [dash pattern={on 2.53pt off 3.02pt}]  (515.5,164) -- (515.5,260.13) ;

\draw [color={rgb, 255:red, 0; green, 0; blue, 0 }  ,draw opacity=1 ][line width=2.25]  [dash pattern={on 2.53pt off 3.02pt}]  (515.5,260.13) -- (480.77,290.12) ;

\draw [color={rgb, 255:red, 0; green, 0; blue, 0 }  ,draw opacity=1 ][line width=2.25]  [dash pattern={on 2.53pt off 3.02pt}]  (423.5,164) -- (423.5,260.13) ;

\draw [color={rgb, 255:red, 144; green, 19; blue, 254 }  ,draw opacity=1 ][line width=2.25]    (423.5,260.13) -- (377.5,291) ;

\draw [color={rgb, 255:red, 126; green, 211; blue, 33 }  ,draw opacity=1 ][line width=2.25]  [dash pattern={on 2.53pt off 3.02pt}]  (423.5,260.13) -- (515.5,260.13) ;

\draw [color={rgb, 255:red, 144; green, 19; blue, 254 }  ,draw opacity=1 ][line width=2.25]    (378.44,201.92) -- (377.5,291) ;

\draw [rotate around= { 46.67: (402.68, 182.07)
    }] [color={rgb, 255:red, 144; green, 19; blue, 254 }  ,draw opacity=1 ][line width=1.5]  (392.21,182.07) -- (413.14,182.07)(402.68,170.88) -- (402.68,193.25) ;
\draw [rotate around= { 46.67: (423.68, 226.07)
    }] [color={rgb, 255:red, 144; green, 19; blue, 254 }  ,draw opacity=1 ][line width=1.5]  (413.21,226.07) -- (434.14,226.07)(423.68,214.88) -- (423.68,237.25) ;
\draw [rotate around= { 46.67: (481.68, 235.07)
    }] [color={rgb, 255:red, 144; green, 19; blue, 254 }  ,draw opacity=1 ][line width=1.5]  (471.21,235.07) -- (492.14,235.07)(481.68,223.88) -- (481.68,246.25) ;
\draw [rotate around= { 46.67: (443.68, 203.07)
    }] [color={rgb, 255:red, 144; green, 19; blue, 254 }  ,draw opacity=1 ][line width=1.5]  (433.21,203.07) -- (454.14,203.07)(443.68,191.88) -- (443.68,214.25) ;
\draw [rotate around= { 46.67: (462.68, 164.07)
    }] [color={rgb, 255:red, 74; green, 144; blue, 226 }  ,draw opacity=1 ][line width=1.5]  (452.21,164.07) -- (473.14,164.07)(462.68,152.88) -- (462.68,175.25) ;
\draw [line width=1.5]    (445,299) -- (405.44,327.26) ;
\draw [shift={(403,329)}, rotate = 324.46000000000004] [color={rgb, 255:red, 0; green, 0; blue, 0 }  ][line width=1.5]    (14.21,-4.28) .. controls (9.04,-1.82) and (4.3,-0.39) .. (0,0) .. controls (4.3,0.39) and (9.04,1.82) .. (14.21,4.28)   ;

\draw [line width=1.5]    (482,299) -- (521.45,323.42) ;
\draw [shift={(524,325)}, rotate = 211.76] [color={rgb, 255:red, 0; green, 0; blue, 0 }  ][line width=1.5]    (14.21,-4.28) .. controls (9.04,-1.82) and (4.3,-0.39) .. (0,0) .. controls (4.3,0.39) and (9.04,1.82) .. (14.21,4.28)   ;

\draw [line width=1.5]    (219,301) -- (251.62,326.17) ;
\draw [shift={(254,328)}, rotate = 217.65] [color={rgb, 255:red, 0; green, 0; blue, 0 }  ][line width=1.5]    (14.21,-4.28) .. controls (9.04,-1.82) and (4.3,-0.39) .. (0,0) .. controls (4.3,0.39) and (9.04,1.82) .. (14.21,4.28)   ;

\draw [line width=1.5]    (180,300) -- (148.35,325.13) ;
\draw [shift={(146,327)}, rotate = 321.55] [color={rgb, 255:red, 0; green, 0; blue, 0 }  ][line width=1.5]    (14.21,-4.28) .. controls (9.04,-1.82) and (4.3,-0.39) .. (0,0) .. controls (4.3,0.39) and (9.04,1.82) .. (14.21,4.28)   ;

\draw [color={rgb, 255:red, 208; green, 2; blue, 27 }  ,draw opacity=1 ][line width=2.25]    (118.77,373.92) -- (118.77,462.12) ;

\draw [color={rgb, 255:red, 208; green, 2; blue, 27 }  ,draw opacity=1 ][line width=2.25]    (16.44,373.92) -- (118.77,373.92) ;

\draw [color={rgb, 255:red, 144; green, 19; blue, 254 }  ,draw opacity=1 ][line width=2.25]    (15.5,463) -- (118.77,462.12) ;

\draw [color={rgb, 255:red, 208; green, 2; blue, 27 }  ,draw opacity=1 ][line width=2.25]    (61.5,336) -- (16.44,373.92) ;

\draw [color={rgb, 255:red, 208; green, 2; blue, 27 }  ,draw opacity=1 ][line width=2.25]    (118.77,373.92) -- (153.5,336) ;

\draw [color={rgb, 255:red, 144; green, 19; blue, 254 }  ,draw opacity=1 ][line width=2.25]    (61.5,336) -- (153.5,336) ;

\draw [color={rgb, 255:red, 144; green, 19; blue, 254 }  ,draw opacity=1 ][line width=2.25]    (153.5,336) -- (153.5,432.13) ;

\draw [color={rgb, 255:red, 208; green, 2; blue, 27 }  ,draw opacity=1 ][line width=2.25]    (153.5,432.13) -- (118.77,462.12) ;

\draw [color={rgb, 255:red, 208; green, 2; blue, 27 }  ,draw opacity=1 ][line width=2.25]    (61.5,336) -- (61.5,432.13) ;

\draw [color={rgb, 255:red, 144; green, 19; blue, 254 }  ,draw opacity=1 ][line width=2.25]    (61.5,432.13) -- (15.5,463) ;

\draw [color={rgb, 255:red, 208; green, 2; blue, 27 }  ,draw opacity=1 ][line width=2.25]    (61.5,432.13) -- (153.5,432.13) ;

\draw [color={rgb, 255:red, 144; green, 19; blue, 254 }  ,draw opacity=1 ][line width=2.25]    (16.44,373.92) -- (15.5,463) ;

\draw [color={rgb, 255:red, 208; green, 2; blue, 27 }  ,draw opacity=1 ][line width=2.25]    (288.77,372.92) -- (288.77,461.12) ;

\draw [color={rgb, 255:red, 208; green, 2; blue, 27 }  ,draw opacity=1 ][line width=2.25]    (186.44,372.92) -- (288.77,372.92) ;

\draw [color={rgb, 255:red, 144; green, 19; blue, 254 }  ,draw opacity=1 ][line width=2.25]    (185.5,462) -- (288.77,461.12) ;

\draw [color={rgb, 255:red, 208; green, 2; blue, 27 }  ,draw opacity=1 ][line width=2.25]    (231.5,335) -- (186.44,372.92) ;

\draw [color={rgb, 255:red, 208; green, 2; blue, 27 }  ,draw opacity=1 ][line width=2.25]    (288.77,372.92) -- (323.5,335) ;

\draw [color={rgb, 255:red, 144; green, 19; blue, 254 }  ,draw opacity=1 ][line width=2.25]    (231.5,335) -- (323.5,335) ;

\draw [color={rgb, 255:red, 208; green, 2; blue, 27 }  ,draw opacity=1 ][line width=2.25]    (323.5,335) -- (323.5,431.13) ;

\draw [color={rgb, 255:red, 144; green, 19; blue, 254 }  ,draw opacity=1 ][line width=2.25]    (323.5,431.13) -- (288.77,461.12) ;

\draw [color={rgb, 255:red, 208; green, 2; blue, 27 }  ,draw opacity=1 ][line width=2.25]    (231.5,335) -- (231.5,431.13) ;

\draw [color={rgb, 255:red, 144; green, 19; blue, 254 }  ,draw opacity=1 ][line width=2.25]    (231.5,431.13) -- (185.5,462) ;

\draw [color={rgb, 255:red, 208; green, 2; blue, 27 }  ,draw opacity=1 ][line width=2.25]    (231.5,431.13) -- (323.5,431.13) ;

\draw [color={rgb, 255:red, 144; green, 19; blue, 254 }  ,draw opacity=1 ][line width=2.25]    (186.44,372.92) -- (185.5,462) ;

\draw [color={rgb, 255:red, 208; green, 2; blue, 27 }  ,draw opacity=1 ][line width=2.25]    (445.77,373.92) -- (445.77,462.12) ;

\draw [color={rgb, 255:red, 208; green, 2; blue, 27 }  ,draw opacity=1 ][line width=2.25]    (343.44,373.92) -- (445.77,373.92) ;

\draw [color={rgb, 255:red, 144; green, 19; blue, 254 }  ,draw opacity=1 ][line width=2.25]    (342.5,463) -- (445.77,462.12) ;

\draw [color={rgb, 255:red, 208; green, 2; blue, 27 }  ,draw opacity=1 ][line width=2.25]    (388.5,336) -- (343.44,373.92) ;

\draw [color={rgb, 255:red, 144; green, 19; blue, 254 }  ,draw opacity=1 ][line width=2.25]    (445.77,373.92) -- (480.5,336) ;

\draw [color={rgb, 255:red, 208; green, 2; blue, 27 }  ,draw opacity=1 ][line width=2.25]    (388.5,336) -- (480.5,336) ;

\draw [color={rgb, 255:red, 208; green, 2; blue, 27 }  ,draw opacity=1 ][line width=2.25]    (480.5,336) -- (480.5,432.13) ;

\draw [color={rgb, 255:red, 208; green, 2; blue, 27 }  ,draw opacity=1 ][line width=2.25]    (480.5,432.13) -- (445.77,462.12) ;

\draw [color={rgb, 255:red, 208; green, 2; blue, 27 }  ,draw opacity=1 ][line width=2.25]    (388.5,336) -- (388.5,432.13) ;

\draw [color={rgb, 255:red, 144; green, 19; blue, 254 }  ,draw opacity=1 ][line width=2.25]    (388.5,432.13) -- (342.5,463) ;

\draw [color={rgb, 255:red, 144; green, 19; blue, 254 }  ,draw opacity=1 ][line width=2.25]    (388.5,432.13) -- (480.5,432.13) ;

\draw [color={rgb, 255:red, 144; green, 19; blue, 254 }  ,draw opacity=1 ][line width=2.25]    (343.44,373.92) -- (342.5,463) ;

\draw [color={rgb, 255:red, 208; green, 2; blue, 27 }  ,draw opacity=1 ][line width=2.25]    (597.77,373.92) -- (597.77,462.12) ;

\draw [color={rgb, 255:red, 208; green, 2; blue, 27 }  ,draw opacity=1 ][line width=2.25]    (495.44,373.92) -- (597.77,373.92) ;

\draw [color={rgb, 255:red, 144; green, 19; blue, 254 }  ,draw opacity=1 ][line width=2.25]    (494.5,463) -- (597.77,462.12) ;

\draw [color={rgb, 255:red, 208; green, 2; blue, 27 }  ,draw opacity=1 ][line width=2.25]    (540.5,336) -- (495.44,373.92) ;

\draw [color={rgb, 255:red, 144; green, 19; blue, 254 }  ,draw opacity=1 ][line width=2.25]    (597.77,373.92) -- (632.5,336) ;

\draw [color={rgb, 255:red, 208; green, 2; blue, 27 }  ,draw opacity=1 ][line width=2.25]    (540.5,336) -- (632.5,336) ;

\draw [color={rgb, 255:red, 144; green, 19; blue, 254 }  ,draw opacity=1 ][line width=2.25]    (632.5,336) -- (632.5,432.13) ;

\draw [color={rgb, 255:red, 208; green, 2; blue, 27 }  ,draw opacity=1 ][line width=2.25]    (632.5,432.13) -- (597.77,462.12) ;

\draw [color={rgb, 255:red, 208; green, 2; blue, 27 }  ,draw opacity=1 ][line width=2.25]    (540.5,336) -- (540.5,432.13) ;

\draw [color={rgb, 255:red, 144; green, 19; blue, 254 }  ,draw opacity=1 ][line width=2.25]    (540.5,432.13) -- (494.5,463) ;

\draw [color={rgb, 255:red, 208; green, 2; blue, 27 }  ,draw opacity=1 ][line width=2.25]    (540.5,432.13) -- (632.5,432.13) ;

\draw [color={rgb, 255:red, 144; green, 19; blue, 254 }  ,draw opacity=1 ][line width=2.25]    (495.44,373.92) -- (494.5,463) ;

\draw [rotate around= { 46.67: (225.68, 263.07)
    }] [color={rgb, 255:red, 74; green, 144; blue, 226 }  ,draw opacity=1 ][line width=1.5]  (215.21,263.07) -- (236.14,263.07)(225.68,251.88) -- (225.68,274.25) ;
\draw [rotate around= { 46.67: (500.68, 272.07)
    }] [color={rgb, 255:red, 74; green, 144; blue, 226 }  ,draw opacity=1 ][line width=1.5]  (490.21,272.07) -- (511.14,272.07)(500.68,260.88) -- (500.68,283.25) ;
\draw [rotate around= { 310.41: (134.55, 360.62)
    }]   (134.55, 360.62) circle [x radius= 36.16, y radius= 12.46]  ;
\draw [rotate around= { 1.12: (586.5, 336)
    }]   (586.5, 336) circle [x radius= 47.49, y radius= 12.46]  ;

\draw (563,485) node  [align=left] {Elliptic over \\circled edge};
\draw (93,489) node  [align=left] {Elliptic over \\circled edge};
\draw (338,494) node  [align=left] {Bottom face has 3 purple \\edges and no pair};

\end{tikzpicture}

    \caption{Three edges arranged around one vertex. We end up with four possible arrangements, each of which presents a problem.}
    \label{fig:my_label}
\end{figure}The first tier of the diagram shows that we have to pick one of the edges on the top face; if not, then the top face has four edges in the other class. Note then that if we pick both edges on the top, we will need to identify the back face and the right face with an elliptic over their shared edge. Thus, we can consider the case where we pick exactly one of the top edges. Once we have picked a top edge, we cannot also pick the parallel bottom edge on the same face because that would leave the adjacent face to have four edges in one class. Thus, we are reduced to only four cases, as shown in Figure 16, and each leads to a problem.

\textbf{Subcase 2}:
Consider now the cube drawn in Figure 17. In picking the other two edges, we have to pick one edge from both the right and the back face or pick their shared edge.

If we pick the shared edge, then we have four choices for the last edge, as shown in Figure 17(b). We can break these into two sets of choices: the two that keep the class connected (green) and the two that do not (red). Choosing either green edge leaves two opposite vertices open, violating Rule 2. Choosing either red edge results in either the top or the bottom face having parallel edges in one class and no face pair. Thus, we cannot pick the shared edge.

If we don't pick the shared edge, we have to pick one edge on the back and the opposite edge on the front. Otherwise, we end up with either the top or the bottom face having three edges. This gives us the two cases in Figure 17(c). Only one is a viable fundamental domain, and it is the 5-7 domain we have already classified.

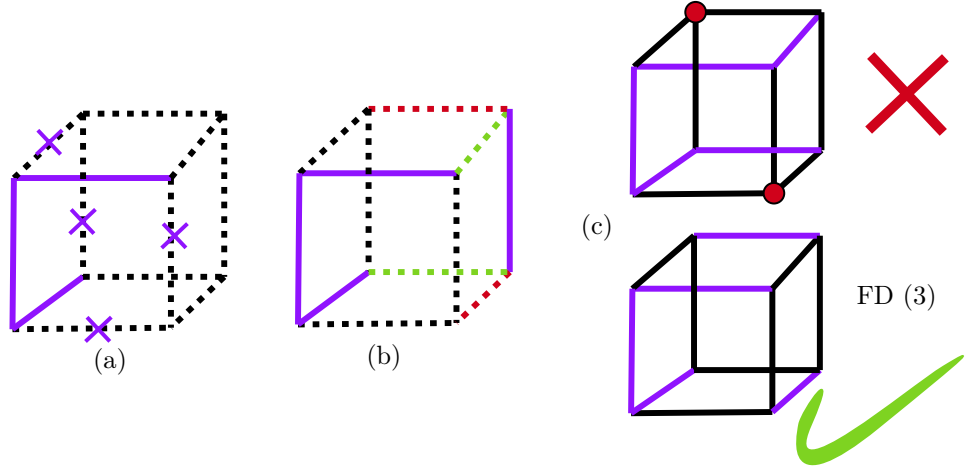
\begin{figure}[H]
    \centering
\tikzset{every picture/.style={line width=0.75pt}} 

\begin{tikzpicture}[x=0.6pt,y=0.6pt,yscale=-1,xscale=1]

\draw [color={rgb, 255:red, 0; green, 0; blue, 0 }  ,draw opacity=1 ][line width=2.25]  [dash pattern={on 2.53pt off 3.02pt}]  (137.4,113.61) -- (137.4,208.06) ;

\draw [color={rgb, 255:red, 144; green, 19; blue, 254 }  ,draw opacity=1 ][line width=2.25]    (38.41,113.61) -- (137.4,113.61) ;

\draw [color={rgb, 255:red, 0; green, 0; blue, 0 }  ,draw opacity=1 ][line width=2.25]  [dash pattern={on 2.53pt off 3.02pt}]  (37.5,209) -- (137.4,208.06) ;

\draw [color={rgb, 255:red, 0; green, 0; blue, 0 }  ,draw opacity=1 ][line width=2.25]  [dash pattern={on 2.53pt off 3.02pt}]  (82,73) -- (38.41,113.61) ;

\draw [color={rgb, 255:red, 0; green, 0; blue, 0 }  ,draw opacity=1 ][line width=2.25]  [dash pattern={on 2.53pt off 3.02pt}]  (137.4,113.61) -- (171,73) ;

\draw [color={rgb, 255:red, 0; green, 0; blue, 0 }  ,draw opacity=1 ][line width=2.25]  [dash pattern={on 2.53pt off 3.02pt}]  (82,73) -- (171,73) ;

\draw [color={rgb, 255:red, 0; green, 0; blue, 0 }  ,draw opacity=1 ][line width=2.25]  [dash pattern={on 2.53pt off 3.02pt}]  (171,73) -- (171,175.94) ;

\draw [color={rgb, 255:red, 0; green, 0; blue, 0 }  ,draw opacity=1 ][line width=2.25]  [dash pattern={on 2.53pt off 3.02pt}]  (171,175.94) -- (137.4,208.06) ;

\draw [color={rgb, 255:red, 0; green, 0; blue, 0 }  ,draw opacity=1 ][line width=2.25]  [dash pattern={on 2.53pt off 3.02pt}]  (82,73) -- (82,175.94) ;

\draw [color={rgb, 255:red, 144; green, 19; blue, 254 }  ,draw opacity=1 ][line width=2.25]    (82,175.94) -- (37.5,209) ;

\draw [color={rgb, 255:red, 0; green, 0; blue, 0 }  ,draw opacity=1 ][line width=2.25]  [dash pattern={on 2.53pt off 3.02pt}]  (82,175.94) -- (171,175.94) ;

\draw [color={rgb, 255:red, 144; green, 19; blue, 254 }  ,draw opacity=1 ][line width=2.25]    (38.41,113.61) -- (37.5,209) ;

\draw [color={rgb, 255:red, 0; green, 0; blue, 0 }  ,draw opacity=1 ][line width=2.25]  [dash pattern={on 2.53pt off 3.02pt}]  (317.4,110.61) -- (317.4,205.06) ;

\draw [color={rgb, 255:red, 144; green, 19; blue, 254 }  ,draw opacity=1 ][line width=2.25]    (218.41,110.61) -- (317.4,110.61) ;

\draw [color={rgb, 255:red, 0; green, 0; blue, 0 }  ,draw opacity=1 ][line width=2.25]  [dash pattern={on 2.53pt off 3.02pt}]  (217.5,206) -- (317.4,205.06) ;

\draw [color={rgb, 255:red, 0; green, 0; blue, 0 }  ,draw opacity=1 ][line width=2.25]  [dash pattern={on 2.53pt off 3.02pt}]  (262,70) -- (218.41,110.61) ;

\draw [color={rgb, 255:red, 126; green, 211; blue, 33 }  ,draw opacity=1 ][line width=2.25]  [dash pattern={on 2.53pt off 3.02pt}]  (317.4,110.61) -- (351,70) ;

\draw [color={rgb, 255:red, 208; green, 2; blue, 27 }  ,draw opacity=1 ][line width=2.25]  [dash pattern={on 2.53pt off 3.02pt}]  (262,70) -- (351,70) ;

\draw [color={rgb, 255:red, 144; green, 19; blue, 254 }  ,draw opacity=1 ][line width=2.25]    (351,70) -- (351,172.94) ;

\draw [color={rgb, 255:red, 208; green, 2; blue, 27 }  ,draw opacity=1 ][line width=2.25]  [dash pattern={on 2.53pt off 3.02pt}]  (351,172.94) -- (317.4,205.06) ;

\draw [color={rgb, 255:red, 0; green, 0; blue, 0 }  ,draw opacity=1 ][line width=2.25]  [dash pattern={on 2.53pt off 3.02pt}]  (262,70) -- (262,172.94) ;

\draw [color={rgb, 255:red, 144; green, 19; blue, 254 }  ,draw opacity=1 ][line width=2.25]    (262,172.94) -- (217.5,206) ;

\draw [color={rgb, 255:red, 126; green, 211; blue, 33 }  ,draw opacity=1 ][line width=2.25]  [dash pattern={on 2.53pt off 3.02pt}]  (262,172.94) -- (351,172.94) ;

\draw [color={rgb, 255:red, 144; green, 19; blue, 254 }  ,draw opacity=1 ][line width=2.25]    (218.41,110.61) -- (217.5,206) ;

\draw [color={rgb, 255:red, 0; green, 0; blue, 0 }  ,draw opacity=1 ][line width=2.25]    (517.55,43.34) -- (517.55,123.2) ;

\draw [color={rgb, 255:red, 144; green, 19; blue, 254 }  ,draw opacity=1 ][line width=2.25]    (429.31,43.34) -- (517.55,43.34) ;

\draw [color={rgb, 255:red, 0; green, 0; blue, 0 }  ,draw opacity=1 ][line width=2.25]    (428.5,124) -- (517.55,123.2) ;

\draw [color={rgb, 255:red, 0; green, 0; blue, 0 }  ,draw opacity=1 ][line width=2.25]    (468.17,9) -- (429.31,43.34) ;

\draw [color={rgb, 255:red, 144; green, 19; blue, 254 }  ,draw opacity=1 ][line width=2.25]    (517.55,43.34) -- (547.5,9) ;

\draw [color={rgb, 255:red, 0; green, 0; blue, 0 }  ,draw opacity=1 ][line width=2.25]    (468.17,9) -- (547.5,9) ;

\draw [color={rgb, 255:red, 0; green, 0; blue, 0 }  ,draw opacity=1 ][line width=2.25]    (547.5,9) -- (547.5,96.05) ;

\draw [color={rgb, 255:red, 0; green, 0; blue, 0 }  ,draw opacity=1 ][line width=2.25]    (547.5,96.05) -- (517.55,123.2) ;

\draw [color={rgb, 255:red, 0; green, 0; blue, 0 }  ,draw opacity=1 ][line width=2.25]    (468.17,9) -- (468.17,96.05) ;

\draw [color={rgb, 255:red, 144; green, 19; blue, 254 }  ,draw opacity=1 ][line width=2.25]    (468.17,96.05) -- (428.5,124) ;

\draw [color={rgb, 255:red, 144; green, 19; blue, 254 }  ,draw opacity=1 ][line width=2.25]    (468.17,96.05) -- (547.5,96.05) ;

\draw [color={rgb, 255:red, 144; green, 19; blue, 254 }  ,draw opacity=1 ][line width=2.25]    (429.31,43.34) -- (428.5,124) ;

\draw [color={rgb, 255:red, 0; green, 0; blue, 0 }  ,draw opacity=1 ][line width=2.25]    (516.55,183.44) -- (516.55,261.22) ;

\draw [color={rgb, 255:red, 144; green, 19; blue, 254 }  ,draw opacity=1 ][line width=2.25]    (428.31,183.44) -- (516.55,183.44) ;

\draw [color={rgb, 255:red, 0; green, 0; blue, 0 }  ,draw opacity=1 ][line width=2.25]    (427.5,262) -- (516.55,261.22) ;

\draw [color={rgb, 255:red, 0; green, 0; blue, 0 }  ,draw opacity=1 ][line width=2.25]    (467.17,150) -- (428.31,183.44) ;

\draw [color={rgb, 255:red, 0; green, 0; blue, 0 }  ,draw opacity=1 ][line width=2.25]    (516.55,183.44) -- (546.5,150) ;

\draw [color={rgb, 255:red, 144; green, 19; blue, 254 }  ,draw opacity=1 ][line width=2.25]    (467.17,150) -- (546.5,150) ;

\draw [color={rgb, 255:red, 0; green, 0; blue, 0 }  ,draw opacity=1 ][line width=2.25]    (546.5,150) -- (546.5,234.78) ;

\draw [color={rgb, 255:red, 144; green, 19; blue, 254 }  ,draw opacity=1 ][line width=2.25]    (546.5,234.78) -- (516.55,261.22) ;

\draw [color={rgb, 255:red, 0; green, 0; blue, 0 }  ,draw opacity=1 ][line width=2.25]    (467.17,150) -- (467.17,234.78) ;

\draw [color={rgb, 255:red, 144; green, 19; blue, 254 }  ,draw opacity=1 ][line width=2.25]    (467.17,234.78) -- (427.5,262) ;

\draw [color={rgb, 255:red, 0; green, 0; blue, 0 }  ,draw opacity=1 ][line width=2.25]    (467.17,234.78) -- (546.5,234.78) ;

\draw [color={rgb, 255:red, 144; green, 19; blue, 254 }  ,draw opacity=1 ][line width=2.25]    (428.31,183.44) -- (427.5,262) ;

\draw [rotate around= { 46.67: (81.68, 141.07)
    }] [color={rgb, 255:red, 144; green, 19; blue, 254 }  ,draw opacity=1 ][line width=1.5]  (71.21,141.07) -- (92.14,141.07)(81.68,129.88) -- (81.68,152.25) ;
\draw [rotate around= { 46.67: (91.68, 209.07)
    }] [color={rgb, 255:red, 144; green, 19; blue, 254 }  ,draw opacity=1 ][line width=1.5]  (81.21,209.07) -- (102.14,209.07)(91.68,197.88) -- (91.68,220.25) ;
\draw [rotate around= { 46.67: (139.68, 150.07)
    }] [color={rgb, 255:red, 144; green, 19; blue, 254 }  ,draw opacity=1 ][line width=1.5]  (129.21,150.07) -- (150.14,150.07)(139.68,138.88) -- (139.68,161.25) ;
\draw [rotate around= { 46.67: (60.68, 91.07)
    }] [color={rgb, 255:red, 144; green, 19; blue, 254 }  ,draw opacity=1 ][line width=1.5]  (50.21,91.07) -- (71.14,91.07)(60.68,79.88) -- (60.68,102.25) ;
\draw  [fill={rgb, 255:red, 208; green, 2; blue, 27 }  ,fill opacity=1 ]  (468.17, 9) circle [x radius= 6.5, y radius= 6.5]  ;
\draw  [fill={rgb, 255:red, 208; green, 2; blue, 27 }  ,fill opacity=1 ]  (517.55, 123.2) circle [x radius= 6.5, y radius= 6.5]  ;
\draw [rotate around= { 46.67: (601.43, 60.57)
    }] [color={rgb, 255:red, 208; green, 2; blue, 27 }  ,draw opacity=1 ][line width=3.75]  (568.09,60.57) -- (634.77,60.57)(601.43,27.31) -- (601.43,93.82) ;
\draw  [color={rgb, 255:red, 126; green, 211; blue, 33 }  ,draw opacity=1 ][fill={rgb, 255:red, 126; green, 211; blue, 33 }  ,fill opacity=1 ] (540.5,270) .. controls (520.5,316) and (602.5,246) .. (624.21,232.75) .. controls (645.93,219.49) and (643.5,222) .. (583.34,270.96) .. controls (523.18,319.91) and (531.96,277.57) .. (533.5,269) .. controls (535.04,260.43) and (560.5,224) .. (540.5,270) -- cycle ;

\draw (99,229) node  [align=left] {(a)};
\draw (271.71,226.59) node  [align=left] {(b)};
\draw (595,189) node  [align=left] {FD (3)};
\draw (406,145) node  [align=left] {(c)};

\end{tikzpicture}

    \caption{(a) The original three edges, two on each adjacent face. The purple 'x's denote forbidden edges.  (b) Once we pick the fourth purple edge, we can either pick the two green edges or the two red edges to finish the class. (c) Picking the green edges is shown in the top cube; picking the red edges is shown in the bottom. The top is not viable, while the bottom gives FD (3).}
    \label{fig:my_label}
\end{figure}

Note that we can skip the case for two edges on the identified faces that are opposite because that case will necessarily require an elliptic.

\textbf{3 edges}:
The last case to check is that there are no new fundamental domain candidates with three edges on adjacent faces. One edge must be shared, and there are three ways that this can happen: either the shared edge is the middle edge of the chain on both faces, on one face, or on neither face. These cases are shown in Figure 18 (a), (b), and (c), respectively. The shared edge is circled in each case.

In Figure 18(a), we must have an elliptic fixing the circled edge. In Figure 18(b), the bottom face is entirely in the purple class and has no pair. Thus, case (c) is the only viable option (there is also a symmetric case that could be considered; however, it has the same problems). In Figure 18(c), there are two faces with three edges, and we claim that the shared edge has to be the middle of the chain for one face and an end of the chain for the other.

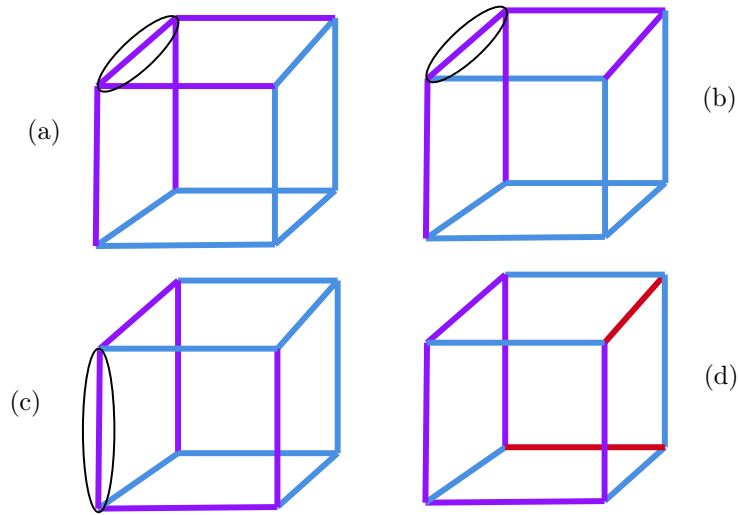
\begin{figure}[H]
    \centering

\tikzset{every picture/.style={line width=0.75pt}} 

\begin{tikzpicture}[x=0.75pt,y=0.75pt,yscale=-1,xscale=1]

\draw [color={rgb, 255:red, 74; green, 144; blue, 226 }  ,draw opacity=1 ][line width=2.25]    (289.86,65.99) -- (289.86,145.94) ;

\draw [color={rgb, 255:red, 144; green, 19; blue, 254 }  ,draw opacity=1 ][line width=2.25]    (200.32,65.99) -- (289.86,65.99) ;

\draw [color={rgb, 255:red, 74; green, 144; blue, 226 }  ,draw opacity=1 ][line width=2.25]    (199.5,146.74) -- (289.86,145.94) ;

\draw [color={rgb, 255:red, 144; green, 19; blue, 254 }  ,draw opacity=1 ][line width=2.25]    (239.75,31.62) -- (200.32,65.99) ;

\draw [color={rgb, 255:red, 74; green, 144; blue, 226 }  ,draw opacity=1 ][line width=2.25]    (289.86,65.99) -- (320.25,31.62) ;

\draw [color={rgb, 255:red, 144; green, 19; blue, 254 }  ,draw opacity=1 ][line width=2.25]    (239.75,31.62) -- (320.25,31.62) ;

\draw [color={rgb, 255:red, 74; green, 144; blue, 226 }  ,draw opacity=1 ][line width=2.25]    (320.25,31.62) -- (320.25,118.76) ;

\draw [color={rgb, 255:red, 74; green, 144; blue, 226 }  ,draw opacity=1 ][line width=2.25]    (320.25,118.76) -- (289.86,145.94) ;

\draw [color={rgb, 255:red, 144; green, 19; blue, 254 }  ,draw opacity=1 ][line width=2.25]    (239.75,31.62) -- (239.75,118.76) ;

\draw [color={rgb, 255:red, 74; green, 144; blue, 226 }  ,draw opacity=1 ][line width=2.25]    (239.75,118.76) -- (199.5,146.74) ;

\draw [color={rgb, 255:red, 74; green, 144; blue, 226 }  ,draw opacity=1 ][line width=2.25]    (239.75,118.76) -- (320.25,118.76) ;

\draw [color={rgb, 255:red, 144; green, 19; blue, 254 }  ,draw opacity=1 ][line width=2.25]    (200.32,65.99) -- (199.5,146.74) ;

\draw [color={rgb, 255:red, 74; green, 144; blue, 226 }  ,draw opacity=1 ][line width=2.25]    (456.36,62.18) -- (456.36,142.13) ;

\draw [color={rgb, 255:red, 74; green, 144; blue, 226 }  ,draw opacity=1 ][line width=2.25]    (366.82,62.18) -- (456.36,62.18) ;

\draw [color={rgb, 255:red, 74; green, 144; blue, 226 }  ,draw opacity=1 ][line width=2.25]    (366,142.93) -- (456.36,142.13) ;

\draw [color={rgb, 255:red, 144; green, 19; blue, 254 }  ,draw opacity=1 ][line width=2.25]    (406.25,27.8) -- (366.82,62.18) ;

\draw [color={rgb, 255:red, 144; green, 19; blue, 254 }  ,draw opacity=1 ][line width=2.25]    (456.36,62.18) -- (486.75,27.8) ;

\draw [color={rgb, 255:red, 144; green, 19; blue, 254 }  ,draw opacity=1 ][line width=2.25]    (406.25,27.8) -- (486.75,27.8) ;

\draw [color={rgb, 255:red, 74; green, 144; blue, 226 }  ,draw opacity=1 ][line width=2.25]    (486.75,27.8) -- (486.75,114.95) ;

\draw [color={rgb, 255:red, 74; green, 144; blue, 226 }  ,draw opacity=1 ][line width=2.25]    (486.75,114.95) -- (456.36,142.13) ;

\draw [color={rgb, 255:red, 144; green, 19; blue, 254 }  ,draw opacity=1 ][line width=2.25]    (406.25,27.8) -- (406.25,114.95) ;

\draw [color={rgb, 255:red, 74; green, 144; blue, 226 }  ,draw opacity=1 ][line width=2.25]    (406.25,114.95) -- (366,142.93) ;

\draw [color={rgb, 255:red, 74; green, 144; blue, 226 }  ,draw opacity=1 ][line width=2.25]    (406.25,114.95) -- (486.75,114.95) ;

\draw [color={rgb, 255:red, 144; green, 19; blue, 254 }  ,draw opacity=1 ][line width=2.25]    (366.82,62.18) -- (366,142.93) ;

\draw [color={rgb, 255:red, 144; green, 19; blue, 254 }  ,draw opacity=1 ][line width=2.25]    (200.75,279.27) -- (291.11,278.47) ;

\draw [color={rgb, 255:red, 144; green, 19; blue, 254 }  ,draw opacity=1 ][line width=2.25]    (241,164.15) -- (201.57,198.53) ;

\draw [color={rgb, 255:red, 74; green, 144; blue, 226 }  ,draw opacity=1 ][line width=2.25]    (291.11,198.53) -- (321.5,164.15) ;

\draw [color={rgb, 255:red, 74; green, 144; blue, 226 }  ,draw opacity=1 ][line width=2.25]    (241,164.15) -- (321.5,164.15) ;

\draw [color={rgb, 255:red, 74; green, 144; blue, 226 }  ,draw opacity=1 ][line width=2.25]    (321.5,164.15) -- (321.5,251.29) ;

\draw [color={rgb, 255:red, 74; green, 144; blue, 226 }  ,draw opacity=1 ][line width=2.25]    (321.5,251.29) -- (291.11,278.47) ;

\draw [color={rgb, 255:red, 144; green, 19; blue, 254 }  ,draw opacity=1 ][line width=2.25]    (241,164.15) -- (241,251.29) ;

\draw [color={rgb, 255:red, 74; green, 144; blue, 226 }  ,draw opacity=1 ][line width=2.25]    (241,251.29) -- (200.75,279.27) ;

\draw [color={rgb, 255:red, 74; green, 144; blue, 226 }  ,draw opacity=1 ][line width=2.25]    (241,251.29) -- (321.5,251.29) ;

\draw [color={rgb, 255:red, 144; green, 19; blue, 254 }  ,draw opacity=1 ][line width=2.25]    (201.57,198.53) -- (200.75,279.27) ;

\draw [color={rgb, 255:red, 144; green, 19; blue, 254 }  ,draw opacity=1 ][line width=2.25]    (291.11,198.53) -- (291.11,278.47) ;

\draw [color={rgb, 255:red, 74; green, 144; blue, 226 }  ,draw opacity=1 ][line width=2.25]    (201.57,198.53) -- (291.11,198.53) ;

\draw [rotate around= { 317.15: (220.88, 49.75)
    }]   (220.88, 49.75) circle [x radius= 26.66, y radius= 8.08]  ;
\draw [rotate around= { 317.15: (386.54, 44.99)
    }]   (386.54, 44.99) circle [x radius= 26.66, y radius= 8.08]  ;
\draw [rotate around= { 270: (200.96, 239.81)
    }]   (200.96, 239.81) circle [x radius= 41.28, y radius= 7.94]  ;
\draw [color={rgb, 255:red, 144; green, 19; blue, 254 }  ,draw opacity=1 ][line width=2.25]    (365.75,276.27) -- (456.11,275.47) ;

\draw [color={rgb, 255:red, 144; green, 19; blue, 254 }  ,draw opacity=1 ][line width=2.25]    (406,161.15) -- (366.57,195.53) ;

\draw [color={rgb, 255:red, 208; green, 2; blue, 27 }  ,draw opacity=1 ][line width=2.25]    (456.11,195.53) -- (486.5,161.15) ;

\draw [color={rgb, 255:red, 74; green, 144; blue, 226 }  ,draw opacity=1 ][line width=2.25]    (406,161.15) -- (486.5,161.15) ;

\draw [color={rgb, 255:red, 74; green, 144; blue, 226 }  ,draw opacity=1 ][line width=2.25]    (486.5,161.15) -- (486.5,248.29) ;

\draw [color={rgb, 255:red, 74; green, 144; blue, 226 }  ,draw opacity=1 ][line width=2.25]    (486.5,248.29) -- (456.11,275.47) ;

\draw [color={rgb, 255:red, 144; green, 19; blue, 254 }  ,draw opacity=1 ][line width=2.25]    (406,161.15) -- (406,248.29) ;

\draw [color={rgb, 255:red, 74; green, 144; blue, 226 }  ,draw opacity=1 ][line width=2.25]    (406,248.29) -- (365.75,276.27) ;

\draw [color={rgb, 255:red, 208; green, 2; blue, 27 }  ,draw opacity=1 ][line width=2.25]    (406,248.29) -- (486.5,248.29) ;

\draw [color={rgb, 255:red, 144; green, 19; blue, 254 }  ,draw opacity=1 ][line width=2.25]    (366.57,195.53) -- (365.75,276.27) ;

\draw [color={rgb, 255:red, 144; green, 19; blue, 254 }  ,draw opacity=1 ][line width=2.25]    (456.11,195.53) -- (456.11,275.47) ;

\draw [color={rgb, 255:red, 74; green, 144; blue, 226 }  ,draw opacity=1 ][line width=2.25]    (366.57,195.53) -- (456.11,195.53) ;

\draw (173.44,89.5) node  [align=left] {(a)};
\draw (514.19,72.78) node  [align=left] {(b)};
\draw (164.06,225.94) node  [align=left] {(c)};
\draw (515,213) node  [align=left] {(d)};

\end{tikzpicture}

    \caption{(a) Shared edge is middle edge for both faces  (b) Shared edge is middle edge for left face only (c) Shared edge is not middle edge for either face  (d) Carrying out face identifications for case (c); we end up with three edge classes.}
    \label{fig:my_label}
\end{figure} Considering the diagram in Figure 18(c), we have one choice of identifications. The left and front faces have to be identified with a $\pi$ twist; sending the right face to the bottom or the back requires an elliptic generator, so we must identify the right face with the top face. That leaves the bottom and the back face to be identified. If we follow these pairings, we obtain that the class of seven edges is actually one class of five edges and another class of two edges, as shown in Figure 18(d).

Thus, we conclude that FD (3) is the only candidates to be fundamental domain on the cube with the edge breakdown 5-7, and more generally, that FD (1)-(3) are the only fundamental domain candidates on the cube with torsion-free groups.

\subsection{Explicitly finding the group for a given set of edge classes}
\begin{figure}[H]
    \centering
 
\tikzset{every picture/.style={line width=0.75pt}} 

\begin{tikzpicture}[x=0.75pt,y=0.75pt,yscale=-1,xscale=1]

\draw [color={rgb, 255:red, 144; green, 19; blue, 254 }  ,draw opacity=1 ][line width=2.25]    (219.5,127) -- (218.5,228) ;

\draw [color={rgb, 255:red, 74; green, 144; blue, 226 }  ,draw opacity=1 ][line width=2.25]    (328.5,127) -- (328.5,227) ;

\draw [color={rgb, 255:red, 74; green, 144; blue, 226 }  ,draw opacity=1 ][line width=2.25]    (219.5,127) -- (328.5,127) ;

\draw [color={rgb, 255:red, 144; green, 19; blue, 254 }  ,draw opacity=1 ][line width=2.25]    (218.5,228) -- (328.5,227) ;

\draw [color={rgb, 255:red, 74; green, 144; blue, 226 }  ,draw opacity=1 ][line width=2.25]    (267.5,84) -- (219.5,127) ;

\draw [color={rgb, 255:red, 144; green, 19; blue, 254 }  ,draw opacity=1 ][line width=2.25]    (328.5,127) -- (365.5,84) ;

\draw [color={rgb, 255:red, 144; green, 19; blue, 254 }  ,draw opacity=1 ][line width=2.25]    (267.5,84) -- (365.5,84) ;

\draw [color={rgb, 255:red, 74; green, 144; blue, 226 }  ,draw opacity=1 ][line width=2.25]    (365.5,84) -- (365.5,193) ;

\draw [color={rgb, 255:red, 144; green, 19; blue, 254 }  ,draw opacity=1 ][line width=2.25]    (365.5,193) -- (328.5,227) ;

\draw [color={rgb, 255:red, 144; green, 19; blue, 254 }  ,draw opacity=1 ][line width=2.25]    (267.5,84) -- (267.5,193) ;

\draw [color={rgb, 255:red, 74; green, 144; blue, 226 }  ,draw opacity=1 ][line width=2.25]    (267.5,193) -- (218.5,228) ;

\draw [color={rgb, 255:red, 74; green, 144; blue, 226 }  ,draw opacity=1 ][line width=2.25]    (267.5,193) -- (365.5,193) ;

\draw (242,198) node  [align=left] {\textcolor[rgb]{0.29,0.56,0.89}{2}};
\draw (276,239) node  [align=left] {\textcolor[rgb]{0.56,0.07,1}{1}};
\draw (359,215) node  [align=left] {\textcolor[rgb]{0.56,0.07,1}{3}};
\draw (308,183) node  [align=left] {\textcolor[rgb]{0.29,0.56,0.89}{4}};
\draw (290,139) node  [align=left] {\textcolor[rgb]{0.29,0.56,0.89}{5}};
\draw (232,102) node  [align=left] {\textcolor[rgb]{0.29,0.56,0.89}{6}};
\draw (352,118) node  [align=left] {\textcolor[rgb]{0.56,0.07,1}{7}};
\draw (316,72) node  [align=left] {\textcolor[rgb]{0.56,0.07,1}{8}};
\draw (205,171) node  [align=left] {\textcolor[rgb]{0.56,0.07,1}{9}};
\draw (341,167) node  [align=left] {\textcolor[rgb]{0.29,0.56,0.89}{10}};
\draw (379,147) node  [align=left] {\textcolor[rgb]{0.29,0.56,0.89}{11}};
\draw (254,151) node  [align=left] {\textcolor[rgb]{0.56,0.07,1}{12}};

\end{tikzpicture}
 \caption{Edge classes for the fundamental domain on the cube. All exterior dihedral angles are $\dfrac{2\pi}{3}$.}
    \label{fig:my_label}

\end{figure}
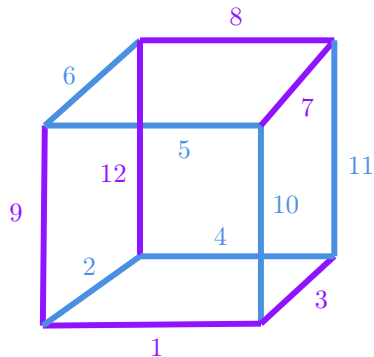

To find the generators, note that we have paired faces based on the face identifications, and thus we can specify that three vertices on the front face, for example, go to three specific vertices on the back face.

Recall that the face identifications for this cube are as follows:

\begin{itemize}
    \item FD (1):

    \subitem A: Front$\longrightarrow$Back with $\dfrac{\pi}{2}$ clockwise twist on back interior face
    \subitem B: Left$\longrightarrow$Right with $\dfrac{\pi}{2}$ clockwise twist on right interior face
    \subitem C: Top$\longrightarrow$Bottom with $\dfrac{\pi}{2}$ clockwise twist on bottom interior face
   
\end{itemize}

A Mobius transformation is completely determined once we know what it does to three points, and thus we can use the cross ratio to find the specific generators. The cross ratio is given by
\begin{equation}
    [z,P_1,P_2,P_3]=\dfrac{(z-P_2)(P_1-P_3)}{(z-P_3)(P_1-P_2)},
\end{equation}where $P_1,P_2,P_3$ go to $0,1,\infty$, respectively.

We first found the vertices of a cube inscribed in the unit cube centered at $(0,0,1)$ in $\mathbb{R}^3$. These are also the vertices of an ideal cube in the Poincare ball; to apply the cross ratio, we want to have the vertices in the upper half space model. To go from the Poincare ball to the upper half space model, the vertices are inverted over the sphere of radius 2 centered at $(0,0,2)$ and then reflected over the xy-plane (see Figure 20).

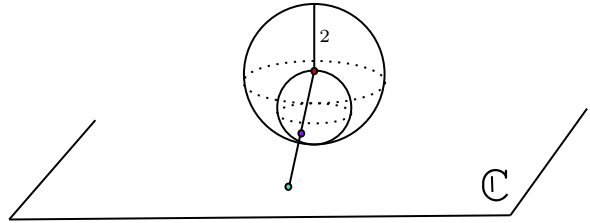
\begin{figure} [H]
    \centering
    \tikzset{every picture/.style={line width=0.75pt}} 

\begin{tikzpicture}[x=0.4pt,y=0.4pt,yscale=-1,xscale=1]

\draw    (343,84) -- (318.5,197) ;

\draw    (343, 87.5) circle [x radius= 66.5, y radius= 66.5]  ;
\draw    (136,131) -- (54.5,225) ;

\draw    (601,120) -- (528.5,221) ;

\draw    (54.5,225) -- (528.5,221) ;

\draw  [dash pattern={on 0.84pt off 2.51pt}]  (343, 124.5) circle [x radius= 35, y radius= 9.5]  ;
\draw    (343,21) -- (343,84) ;

\draw  [dash pattern={on 0.84pt off 2.51pt}]  (343, 95) circle [x radius= 67, y radius= 20]  ;
\draw    (343, 119) circle [x radius= 35, y radius= 35]  ;
\draw  [fill={rgb, 255:red, 144; green, 19; blue, 254 }  ,fill opacity=1 ] (328,143.5) .. controls (328,141.84) and (329.23,140.5) .. (330.75,140.5) .. controls (332.27,140.5) and (333.5,141.84) .. (333.5,143.5) .. controls (333.5,145.16) and (332.27,146.5) .. (330.75,146.5) .. controls (329.23,146.5) and (328,145.16) .. (328,143.5) -- cycle ;
\draw  [fill={rgb, 255:red, 208; green, 2; blue, 27 }  ,fill opacity=1 ] (340.25,84.5) .. controls (340.25,82.84) and (341.48,81.5) .. (343,81.5) .. controls (344.52,81.5) and (345.75,82.84) .. (345.75,84.5) .. controls (345.75,86.16) and (344.52,87.5) .. (343,87.5) .. controls (341.48,87.5) and (340.25,86.16) .. (340.25,84.5) -- cycle ;
\draw  [fill={rgb, 255:red, 80; green, 227; blue, 194 }  ,fill opacity=1 ] (315.75,194) .. controls (315.75,192.34) and (316.98,191) .. (318.5,191) .. controls (320.02,191) and (321.25,192.34) .. (321.25,194) .. controls (321.25,195.66) and (320.02,197) .. (318.5,197) .. controls (316.98,197) and (315.75,195.66) .. (315.75,194) -- cycle ;
\draw    (511,184) -- (511.5,199) ;

\draw (353,51) node  [align=left] {{\scriptsize 2}};
\draw (514,192) node  [align=left] {{\fontfamily{pcr}\selectfont {\huge C}}};

\end{tikzpicture}

    \caption{Inversion over the sphere containing the Poincare ball. A vertex (purple) on the cube inscribed in the ball model is sent to the new vertex (green).}
    \label{fig:my_label}
\end{figure}
Because all of the points lie on a sphere passing through the south pole and the center of the sphere of inversion, all of the vertices land in the $z=0$ plane. In the upper half model, the $z=0$ plane is $\mathbb{C}$, the complex plane. We then used the cross ratio along with the scheme in Figure 21 to determine the Mobius transformation.

\tikzset{every picture/.style={line width=0.75pt}} 
\begin{figure}[H]
    \centering

\begin{tikzpicture}[x=0.75pt,y=0.75pt,yscale=-1,xscale=1]

\draw    (204,120) -- (257.92,161.78) ;
\draw [shift={(259.5,163)}, rotate = 217.77] [color={rgb, 255:red, 0; green, 0; blue, 0 }  ][line width=0.75]    (10.93,-3.29) .. controls (6.95,-1.4) and (3.31,-0.3) .. (0,0) .. controls (3.31,0.3) and (6.95,1.4) .. (10.93,3.29)   ;

\draw    (376,117) -- (322.03,162.71) ;
\draw [shift={(320.5,164)}, rotate = 319.74] [color={rgb, 255:red, 0; green, 0; blue, 0 }  ][line width=0.75]    (10.93,-3.29) .. controls (6.95,-1.4) and (3.31,-0.3) .. (0,0) .. controls (3.31,0.3) and (6.95,1.4) .. (10.93,3.29)   ;

\draw    (205,103) -- (375.5,103) ;
\draw [shift={(377.5,103)}, rotate = 180] [color={rgb, 255:red, 0; green, 0; blue, 0 }  ][line width=0.75]    (10.93,-3.29) .. controls (6.95,-1.4) and (3.31,-0.3) .. (0,0) .. controls (3.31,0.3) and (6.95,1.4) .. (10.93,3.29)   ;

\draw (214,151) node  [align=left] {X};
\draw (159,103) node  [align=left] {$\displaystyle ( l_{1} ,\ l_{2} ,\ l_{3} \ )$};
\draw (291,177.33) node  [align=left] {$\displaystyle ( 0,1,\infty )$};
\draw (422,99) node  [align=left] {$\displaystyle ( r_{1} ,r_{2} ,r_{3})$};
\draw (365,150) node  [align=left] {Y};
\draw (291,53) node  [align=left] {Left face$\displaystyle \rightarrow $Right face};
\draw (296,84) node  [align=left] {$\displaystyle XY^{-1}$};

\end{tikzpicture}
\caption{Diagram demonstrating how to find the Mobius transformation for the left to right face identification.}
    \label{fig:my_label}
\end{figure}
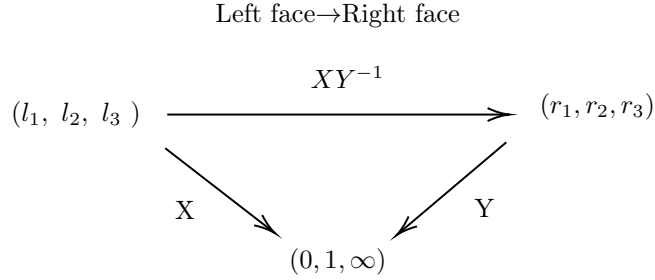

Following this pattern, we found the three generators to be:
\begin{equation}
    A(z)=\dfrac{(i-\sqrt{3})z+4}{z+(i-\sqrt{3})};  B(z)=\dfrac{(1-\sqrt{3}i)z+4}{-z+(1-\sqrt{3}i)}; C(z)=\dfrac{(1-\sqrt{3})(1-i)z}{-(1+\sqrt{3})(1+i)}
\end{equation}

To determine the relators for the group, we start at an edge and follow the face identifications until we return to the original edge. Intuitively, these products should be the identity because we return to the starting location; however, we will need to explicitly check to verify that this is the case (see Section 5.3). Combining the results with the generators above, we obtain the group given by:
\begin{equation}
    G=\langle A,B,C\ |\ AB^{-1}CA^{-1}B^{-1}C^{-1}, ABC^{-1}A^{-1}BC\rangle.
\end{equation} There are three generators, corresponding to the three pairwise face identifications, and two relators corresponding to the two edge classes.

We further prove that for this particular set of face identifications, there is no other angle solution that satisfies both the Rivin equations and the fundamental domain equations.

\textbf{Proposition 5.2}: Given the face identifications above, the set of angles all equal to $\dfrac{2\pi}{3}$ is the \textit{unique} angle solution.

\begin{proof}

 Between the Rivin equations and the fundamental equations, we have a system of 10 equations in 12 variables: 8 equations corresponding to the eight vertices and one equation for each of the two edge classes. Refer to Figure 19 for the labelling of the edges. First we prove that each of the connected components of an edge class (in this particular case) sums to $2\pi$ (of the total $4\pi$). We will show it for one of the edge classes and assert that it is true for the other (it is straightforward to check). Suppose that the blue edge class is not divided evenly between its components and that one is equal to $2\pi+k$ and one is equal to $2\pi-k$, where $k\in\mathbb{Z}$ can be positive or negative. Assume WLOG that the component containing edges 5, 6 and 10 is equal to $2\pi+k$. Starting with the fundamental domain equation for the purple edge class, we have:
\begin{center}
$\pi-x_1+\pi-x_3+\pi-x_9+\pi-x_{12}+\pi-x_8+\pi-x_7=2\pi$\\
$\rightarrow 6\pi-(x_1+x_3+x_9+x_{12}+x_8+x_7)=2\pi$\\
$\rightarrow x_1+x_3+x_9+x_{12}+x_8+x_7=4\pi$\\
\end{center}
Then, by making two substitutions from the Rivin equations, we can obtain an equation containing 2 equations from each of the classes. We will substitute
$x_1+x_3=2\pi-x_{10}$ and $x_{12}+x_8=2\pi-x_6$ in and obtain:
\begin{center}

 $2\pi-x_{10}+x_9+2\pi-x_6+x_7=\pi$\\
$\longrightarrow x_{10}-x_{9}=x_7-x_6$
   
\end{center}

Note that we can also combine the equations for the component $x_5+x_6+x_{10}=2\pi+k$ (WLOG because k can have either sign) with the two vertex equations that overlap in two edges: we have $x_5+x_6+x_9=2\pi$ and $x_5+x_{10}+x_7=2\pi$. Then, by combining each of these equations with the component equation, we get the two equations
\begin{equation}
x_6-x_7=k\  \mbox{and}\  x_{10}-x_9=k.
\end{equation} But combining with the above equations, we get $k=-k$ and thus $k=0$. This tells us that both components of the blue edge class are equal to $2\pi$. Similar arguments can be made for the purple class starting with the equation for the blue edge class and making two substitutions. Thus, we now have 8 equations from the vertices and 4 equations for the 4 connected components, all equal to $2\pi$. Each equation has three vertices being added together, and they are linearly independent because no vertex has three equations in one class and thus they are distinct sums of three vertices. 12 equations in 12 variables has a unique solution, and we know what that solution is.
\end{proof}

\subsection{Determining Whether the Candidates are Fundamental Domains}

In order to determine whether or not these polyhedra and groups are viable fundamental domains, we need to ensure that the relators are equal to the identity. We have only imposed the condition that the angle sum of edges in a class is $2\pi$; it is still possible that the relator is a loxodromic containing the starting edge in the class in its axis.

For FD (1), we used the group calculated in section 5.2 to verify that the relators are indeed the identity; thus, FD(1) is a fundamental domain. We only explicitly verified this claim for the version of FD (1) shown in Figure 19; to see that it is true for the mirror image version, note that the face identifications in that case require a $\dfrac{\pi}{2}$ counterclockwise twist, meaning that group generators are the inverses of the generators for the group that we checked. Let $X=A^{-1}, Y=B^{-1}, Z=C^{-1}$. The relators can be found as \begin{equation}
    XYZX^{-1}YZ^{-1}\ \mbox{and}\ XY^{-1}Z^{-1}X^{-1}Y{-1}Z.
\end{equation} When we plug in the substitutions for $A,B,$ and $C$, we obtain the relators that we just verified were the identity.

For FD (2), we obtain the following group using the same method as in Section 5.2:
\begin{center}
\begin{equation}
 G=\langle P,Q,R\ |\ PR^{-1}R^{-1}PQ^{-1}Q^{-1}, PQR^{-1}P^{-1}Q^{-1}R\rangle,
\end{equation}
\end{center}where $P,Q,$ and $R$ are given by
\begin{equation}
    P(z)=\dfrac{2(1+i)z-4\sqrt{3}(1+i)}{-\sqrt{3}(1+i)z+2(1+i)};
\end{equation}

\begin{equation}
     Q(z)=\dfrac{(2+2\sqrt{3}+i(-2+2\sqrt{3}))z+20+12\sqrt{3}+i(4+4\sqrt{3})}{(\sqrt{3}-1+i(-1-\sqrt{3}))z-2\sqrt{3}-2+i(10+6\sqrt{3})}
\end{equation}
and
\begin{equation}
R(z)=\dfrac{(10\sqrt{3}-18+(6\sqrt{3}-10)i)z-12\sqrt{3}+20+(20\sqrt{3}-36)i}{(\sqrt{3}-3+(-1+\sqrt{3})i)z+2\sqrt{3}-2+(6-2\sqrt{3})i}.
\end{equation}

Plugging in these equations to the relators in Equation 20 yields a multiple of the identity, confirming that FD (2) is indeed a fundamental domain.

We will not explicitly treat FD (3) (and its mirror image) in this paper. The same procedure of finding the generators and relators and confirming that the relators yield the identity will work for this case as well; however, finding the vertices of the cube is more difficult here because the exterior dihedral angles are not all equal to $\dfrac{2\pi}{3}$. We leave this verification for future work.

\section{Combinatorial Restrictions on General Fundamental Domains}

After investigating fundamental domains on the cube specifically, we want to consider fundamental domains on more general polyhedra. One fact that helped us fully classify the fundamental domains on the cube with non-elliptic generators was the fact that there were only certain edge breakdowns that were allowed, restricting our consideration to fewer cases. Extending this idea, we wondered if there might be a more general cutoff for the size of an edge class. We know that an edge class has to have at least 2 edges; but how can we impose conditions on the polyhedron and group to improve that bound?

The following results give a set of conditions for which we know that we can improve the lower bound on edge class size to three edges.

\textbf{Proposition 6.1:} If no pair of generators in the group corresponding to the fundamental domain commute with each other, then the fundamental domain does not have an edge class of size 3.

\begin{proof}
The proof will consist of Claims 6.1 and 6.2, which when combined give the full result.

\textbf{Claim 6.1:} There is an edge class containing three edges $\iff$ $Y^2Z$ is a relator in the corresponding group, where $Y,Z$ are generators of the group.
\begin{proof}

$\implies$ Suppose that there is an edge class containing three edges. Then there must be two adjacent faces in the polyhedron that are identified with each other, and they must share an edge. Otherwise, we will not be able to put face identifications on the polyhedron.

\begin{figure}[H]
    \centering

\tikzset{every picture/.style={line width=0.75pt}} 

\begin{tikzpicture}[x=0.75pt,y=0.75pt,yscale=-1,xscale=1]

\draw [color={rgb, 255:red, 74; green, 144; blue, 226 }  ,draw opacity=1 ][line width=2.25]    (268.37,151.99) -- (251,192.83) ;

\draw    (326.5,145.65) -- (354.44,180.16) ;

\draw    (268.37,151.99) -- (326.5,145.65) ;

\draw    (277.43,233.67) -- (331.79,230.86) ;

\draw    (251,192.83) -- (277.43,233.67) ;

\draw  [fill={rgb, 255:red, 0; green, 0; blue, 0 }  ,fill opacity=1 ] (349.15,191.78) .. controls (349.15,190.8) and (349.83,190.01) .. (350.66,190.01) .. controls (351.5,190.01) and (352.17,190.8) .. (352.17,191.78) .. controls (352.17,192.75) and (351.5,193.54) .. (350.66,193.54) .. controls (349.83,193.54) and (349.15,192.75) .. (349.15,191.78) -- cycle ;
\draw  [fill={rgb, 255:red, 0; green, 0; blue, 0 }  ,fill opacity=1 ] (343.11,205.15) .. controls (343.11,204.18) and (343.79,203.39) .. (344.62,203.39) .. controls (345.46,203.39) and (346.13,204.18) .. (346.13,205.15) .. controls (346.13,206.13) and (345.46,206.91) .. (344.62,206.91) .. controls (343.79,206.91) and (343.11,206.13) .. (343.11,205.15) -- cycle ;
\draw  [fill={rgb, 255:red, 0; green, 0; blue, 0 }  ,fill opacity=1 ] (336.32,219.24) .. controls (336.32,218.27) and (336.99,217.48) .. (337.83,217.48) .. controls (338.66,217.48) and (339.34,218.27) .. (339.34,219.24) .. controls (339.34,220.21) and (338.66,221) .. (337.83,221) .. controls (336.99,221) and (336.32,220.21) .. (336.32,219.24) -- cycle ;
\draw [color={rgb, 255:red, 74; green, 144; blue, 226 }  ,draw opacity=1 ][line width=2.25]    (326.5,145.65) -- (354.44,180.16) ;

\draw [color={rgb, 255:red, 74; green, 144; blue, 226 }  ,draw opacity=1 ][line width=2.25]    (358.97,105.52) -- (412.57,104.81) ;

\draw    (326.5,145.65) -- (358.97,105.52) ;

\draw    (439,145.65) -- (409.55,185.09) ;

\draw    (412.57,104.81) -- (439,145.65) ;

\draw  [fill={rgb, 255:red, 0; green, 0; blue, 0 }  ,fill opacity=1 ] (362.74,181.21) .. controls (362.74,180.24) and (363.42,179.45) .. (364.25,179.45) .. controls (365.09,179.45) and (365.76,180.24) .. (365.76,181.21) .. controls (365.76,182.19) and (365.09,182.97) .. (364.25,182.97) .. controls (363.42,182.97) and (362.74,182.19) .. (362.74,181.21) -- cycle ;
\draw  [fill={rgb, 255:red, 0; green, 0; blue, 0 }  ,fill opacity=1 ] (392.94,184.73) .. controls (392.94,183.76) and (393.62,182.97) .. (394.45,182.97) .. controls (395.29,182.97) and (395.96,183.76) .. (395.96,184.73) .. controls (395.96,185.71) and (395.29,186.49) .. (394.45,186.49) .. controls (393.62,186.49) and (392.94,185.71) .. (392.94,184.73) -- cycle ;
\draw  [fill={rgb, 255:red, 0; green, 0; blue, 0 }  ,fill opacity=1 ] (377.09,182.62) .. controls (377.09,181.65) and (377.76,180.86) .. (378.6,180.86) .. controls (379.43,180.86) and (380.11,181.65) .. (380.11,182.62) .. controls (380.11,183.59) and (379.43,184.38) .. (378.6,184.38) .. controls (377.76,184.38) and (377.09,183.59) .. (377.09,182.62) -- cycle ;
\draw    (335.56,237.9) .. controls (360.72,259.98) and (425.76,217.14) .. (416.98,188.91) ;
\draw [shift={(416.35,187.2)}, rotate = 426.78999999999996] [fill={rgb, 255:red, 0; green, 0; blue, 0 }  ][line width=0.75]  [draw opacity=0] (8.93,-4.29) -- (0,0) -- (8.93,4.29) -- cycle    ;

\draw    (347,103) .. controls (335.1,70.4) and (246.02,111.17) .. (273.66,140.66) ;
\draw [shift={(275,142)}, rotate = 223] [fill={rgb, 255:red, 0; green, 0; blue, 0 }  ][line width=0.75]  [draw opacity=0] (8.93,-4.29) -- (0,0) -- (8.93,4.29) -- cycle    ;

\draw    (196,201) -- (251,192.83) ;

\draw    (234,103) -- (268.37,151.99) ;

\draw    (334,59) -- (358.97,105.52) ;

\draw    (412.57,104.81) -- (443,67) ;

\draw    (207,127) .. controls (195.1,94.4) and (334.02,18.74) .. (368.5,44.74) ;
\draw [shift={(370,46)}, rotate = 223] [fill={rgb, 255:red, 0; green, 0; blue, 0 }  ][line width=0.75]  [draw opacity=0] (8.93,-4.29) -- (0,0) -- (8.93,4.29) -- cycle    ;

\draw  [fill={rgb, 255:red, 0; green, 0; blue, 0 }  ,fill opacity=1 ] (256.11,108.08) .. controls (256.16,107.11) and (256.87,106.36) .. (257.7,106.4) .. controls (258.54,106.43) and (259.18,107.25) .. (259.13,108.22) .. controls (259.09,109.2) and (258.37,109.95) .. (257.54,109.91) .. controls (256.71,109.87) and (256.07,109.06) .. (256.11,108.08) -- cycle ;
\draw  [fill={rgb, 255:red, 0; green, 0; blue, 0 }  ,fill opacity=1 ] (285.11,90.15) .. controls (285.11,89.18) and (285.79,88.39) .. (286.62,88.39) .. controls (287.46,88.39) and (288.13,89.18) .. (288.13,90.15) .. controls (288.13,91.13) and (287.46,91.91) .. (286.62,91.91) .. controls (285.79,91.91) and (285.11,91.13) .. (285.11,90.15) -- cycle ;
\draw  [fill={rgb, 255:red, 0; green, 0; blue, 0 }  ,fill opacity=1 ] (313.15,75.78) .. controls (313.15,74.8) and (313.83,74.01) .. (314.66,74.01) .. controls (315.5,74.01) and (316.17,74.8) .. (316.17,75.78) .. controls (316.17,76.75) and (315.5,77.54) .. (314.66,77.54) .. controls (313.83,77.54) and (313.15,76.75) .. (313.15,75.78) -- cycle ;

\draw (350.42,155.51) node  [align=left] {$\displaystyle \textcolor[rgb]{0.29,0.56,0.89}{e}$};
\draw (269.12,175.23) node  [align=left] {$\displaystyle \textcolor[rgb]{0.29,0.56,0.89}{x_{a}}$};
\draw (386.9,113.97) node  [align=left] {$\displaystyle \textcolor[rgb]{0.29,0.56,0.89}{x_{b}}$};
\draw (300.83,186.49) node   {$A$};
\draw (384.64,138.61) node  [align=left] {$\displaystyle B$};
\draw (383.13,218.18) node  [align=left] {$\displaystyle Y$};
\draw (305.36,114.67) node   {$Y^{-1}$};
\draw (220,153) node   {$C$};
\draw (387,75) node   {$D$};
\draw (264,46) node   {$Z$};

\end{tikzpicture}

    \caption{Faces $A$ an $B$ are identified by generator $Y$ and share edge $e$. We will show that C an D must be identified by some generator $Z$.}
    \label{fig:my_label}
\end{figure}
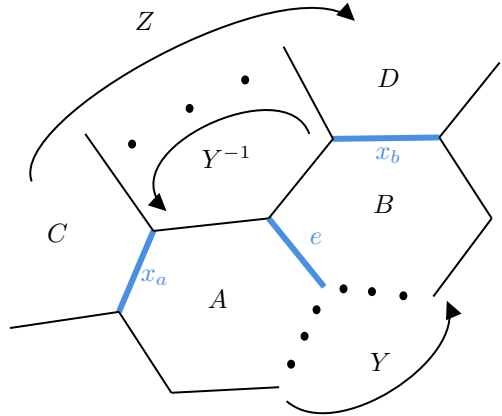

Figure 22 shows the general setup. Faces $A$ and $B$ are adjacent and share edge $e$; the edge class also contains edges $x_a$ and $x_b$. Edge $x_a$ is also contained in face C, and edge $x_b$ is also contained in face D. Face $A$ is sent to face $B$ by isometry $Y$. We will not worry about isometry $Z$ yet.

The relators of the group are constructed by moving from edge to edge in a given class and recording the generators (face identifications) that are used at each step. Each edge is contained in two faces, and we use this fact to "jump" from edge to edge. For example, consider a cube where the front is identified with the back via generator P and the bottom is identified with the top via generator K. Consider the edge shared between the bottom and the back faces. We can get to that edge from the front face via generator X. Now we want to consider the edge as part of the bottom face because if we consider it as part of the back face, we will end up back at the original edge on the front face. In this way, we have already accounted for that edge as part of the back face. We then use generator K to account for the edge as part of the bottom face. This allows us to reach every edge in the class, and eventually we come back to the starting edge.

We can choose any edge in the class as a starting point. Because we have to pick one face to start with, the other face that contains the starting edge will be the endpoint of the last generator in the relator, bringing us back to the original edge. If we have adjacent faces that are identified with each other via generator $Y$ sharing an edge, then if we start with the shared edge, the two faces containing the edge are identified with each other, and so the first and terms of the relator are each either $Y$ or $Y^{-1}$. We argue that they both have to be the same; for if they were not, then we could just consider the case where we find the relator starting at another edge, which gives a cyclic permutation of this relator; there would then be an adjacent $YY^{-1}$, which should not happen. Thus, both the first and last terms are the same generator.

In this specific case, we will start the relator at edge $e$. The first and last term in the relator are either both $Y$ or both $Y^{-1}$ by the argument above. We will consider the case with both being $Y$. Then our relator so far looks like this: $Y$\underline{\ \ \ }$Y$. We want to fill in the blank space. Since our first relator is $Y$, that means that we consider $e$ as part of face $A$ originally. Thus, the missing generator takes edge $x_b$ to edge $x_a$, as the last generator in the relator takes $x_a$ to $e$ via $Y$.

Edge $x_b$ is contained in two faces: $B$ and $D$. However, if we consider $x_b$ as part of face $B$, then we have to use the generator $Y^{-1}$, and we know that $Y^{-1}$ sends $x_b$ to $e$, so we cannot use face B. Thus, we consider $x_b$ as part of face D. Since $x_b$ has to go to $x_a$ and face $A$ is already identified with face $B$, this means that face $D$ has to be identified with face $C$ and has to match $x_b$ and $x_a$. Thus, the relator is of the form $YZY$. As we mentioned earlier, any cyclic permutation of the relator can be used, so we can rearrange this relator to obtain $Y^2Z$, which is what we set out to show.

$\impliedby$ The term $YYZ$ has three generators, which means that there are three edges being identified. Suppose that there were originally more generators that were cancelled out. Then there would have been adjacent terms of the form $XX^{-1}$. We will use an example to explain why this cannot happen: suppose WLOG that the front and back faces of a cube are identified so that edge $f$ on the front face is sent to edge $b$ on the back face. Let $X$ be the generator that sends the front face to the back face. When we are at $f$ in the chain of edges in the class, we act with $X$ to get to $b$. But now we are on the back face and we are supposed to use $X^{-1}$, which sends the back face to the front face; in particular, it sends $b$ back to $f$. Thus, we have ended up at the same edge we started at, meaning that those two edges should be in their own class, a contradiction. So we cannot have cancelled out any terms of the form $XX^{-1}$, and the class must have three edges if there is a relator of the form $Y^2Z$.

\end{proof}

\noindent Now we prove the following result:

\textbf{Claim 6.2}: $Y^2Z$ is a relator $\implies$ $Y,Z$ commute.

\begin{proof}
We begin with $Y^2Z=1$. Thus, $Y=Z^{-1}Y^{-1}$ and $Y=Y^{-1}Z^{-1}$, obtained by solving for each of the $Y$s. Then:
\begin{equation}
    Z^{-1}Y^{-1}=Y^{-1}Z^{-1} \implies  Y^{-1}= ZY^{-1}Z^{-1}\implies Y^{-1}Z=ZY^{-1}
\end{equation}
\begin{equation*}
     \implies Z=YZY^{-1}\implies ZY=YZ
\end{equation*}
\end{proof}
By the contrapositive of Claim 6.2, if Y and Z do not commute, then $Y^2Z$ is not a relator; and by Claim 6.1, if $Y^2Z$ is not a relator, then there is not an edge class containing three edges. This combination gives us the full result that we were looking for.

\end{proof}

One might ask when we can extend Claim 6.2 to an "if and only if" statement. We cannot in general, but the following result gives us a partial extension.

\noindent \textbf{Claim 6.3}: Y,Z commute $\implies$ $Y^2Z=K$ where $K$ commutes with $Y$ ($K\neq Z$).

\begin{proof} Suppose that $Y$ and $Z$ commute. Then $YZ=ZY$ and thus $Y^{-1}Z^{-1}=Z^{-1}Y^{-1}$. Let $Y^2Z=K$. Then $Y=KZ^{-1}Y^{-1}$ and $Y=Y^{-1}KZ^{-1}$.
\begin{equation}
    KZ^{-1}Y^{-1}=Y^{-1}KZ^{-1} \implies K=YKZ^{-1}Y^{-1}Z.
\end{equation}
Then since $Y^{-1}$ and $Z^{-1}$ commute,
\begin{equation}
    K=YKY^{-1}Z^{-1}Z=YKY^{-1}\implies KY=YK,
\end{equation} so K commutes with Y. Note that $K\neq Z$ because if it was, then $Y^2Z=Z\implies Y^2=1$. But $Y$ does not have order 2 because it is squared in the relator but applying $Y$ twice in that case does not take us back to the original edge.
\end{proof}

\noindent\textbf{Corollary 6.1:} If $Z$ is the only element in the group that commutes with $Y$, then $Y^2Z=1$.
\noindent\textbf{Corollary 6.2:} If there is a generator in the group that only commutes with one other element that is also a generator, then there is an edge class with three edges.

Now that we have proven Proposition 6.1, we can prove the following result:

\textbf{Proposition 6.2}: If the group corresponding to a fundamental domain does not have any generators that commute, then the polyhedron must satisfy the equation $\overline{E}\leq 2\overline{V}$.

\begin{proof}
Suppose that the group does not have any generators that commute. Then by Proposition 6.1, there is not an edge class of size three. Thus, all edge classes have at least four elements. That means that the number of edges in the abstract polyhedron divided by the number of edge classes in the quotient must be at least 4. Thus:

\begin{equation}
    \dfrac{\overline{E}}{E}\geq 4 \implies \dfrac{\overline{E}}{\frac{\overline{E}-\overline{V}}{2}}\geq 4 \implies \dfrac{\overline{E}}{\overline{E}-\overline{V}}\geq 2 \implies \overline{E}\geq 2\overline{E}-2\overline{V} \implies \overline{E}\leq 2\overline{V}
\end{equation}
\end{proof}

\textbf{Corollary 6.3}: Any fundamental domain on the icosahedron without elliptic generators must have at least two generators that commute.

\begin{proof}
The icosahedron has 30 edges and 9 edge classes, and $\dfrac{30}{9}<4$. Thus, it violates the equation in Proposition 6.2, and by the contrapositive of that statement, it has generators that commute.
\end{proof}

We can further generalize the methods used above to edge classes of odd size.

\textbf{Proposition 6.3}: There is a squared term in a relator for the group $\iff$ A pair of adjacent faces are identified in the polyhedron.

\begin{proof}
$\implies$ Assume that there is a squared term in one of the relators for the group. That is, the relator is of the form $X_1...X_i^2...X_k$, where the $X_j$ are generators for the group. Let $X_i$ be the generator that takes face $A$ to face $B$. Suppose that we get edge $e$ from $X_{i-1}$. Then the image face of $X_{i-1}$ is face A. We act $X_i$ on face A and obtain face B. Edge $e$ is now edge $X_i(e)$ on face B. We do not want the generator that has face B as its domain because that is $X_i^{-1}$ and will undo what we just did. So the next term in the relator comes from considering $X_i(e)$ as part of its other face. But the next generator is $X_i$ again, which has face A as its domain face. Thus, edge $e$ is shared between faces A and B, which are identified by $X_i$.

$\impliedby$ Assume that a pair of adjacent faces are identified. Call them faces A and B and let them share edge $e$. Isometry X sends A to B. We will proceed as we did in the proof of Claim 6.1. WLOG, we can construct the relator for the class of edge $e$ by starting at $e$. Edge $e$ is part of faces A and B; we will first consider $e$ as part of face A (again WLOG), so the first term in the relator is X. That means that we want to end with the isometry that takes another edge back to $e$, where $e$ is part of face B. In other words, we need the isometry that has B as its codomain face. But that isometry is $X: A\rightarrow B$. Thus, the relator is of the form $X...X$, and if we take a cyclic permutation, the generator $X$ will be squared.
\end{proof}

\noindent\textbf{Corollary 6.4:} If no relators have squared terms, then there are no edge classes of odd size and moreover, both $\overline{E}$ and $\overline{V}$ have to be even.

\begin{proof}
If no relators have squared terms, then by Proposition 6.3, there are no adjacent faces identified. In order to have an edge class of odd size, we need to identify adjacent faces so that there can be a shared edge. Otherwise, every edge has a pair somewhere else and there are an even number of edges. Thus, there are no edge classes of odd size. To see that $\overline{E}$ and $\overline{V}$ are both even, note first that in order for the equation $E=\dfrac{\overline{E}-\overline{V}}{2}$ to be well defined, $\overline{E}$ and $\overline{V}$ have to have the same parity. If $\overline{E}$ is odd, then there must be at least one edge class of odd size. Thus, in this case both $\overline{E}$ and $\overline{V}$ must be even.
\end{proof}

These results provide us with methods of checking whether or not a given group and/or polyhedron is suited to be a fundamental domain.

\section{Conclusions and Next Steps}

The study of fundamental domains is a key component of current research in hyperbolic 3-manifolds. One approach is to determine different fundamental domains, and in this paper we have explored combinatorial restrictions on fundamental polyhedra and conditions on the associated groups. These constraints offer a way of determining whether given polyhedra and/or groups are candidates to be fundamental domains. We only considered torsion-free groups; a potential future direction would be to consider groups with elliptic elements. In Section 3.2.1, we described the role of elliptic elements in general groups; an approach that could be taken would be to go back through the proofs outlined in Section 5 and figure out what is happening with the cases where there are elliptic generators. We could also attempt to make a modification to the fundamental domain equations used in this paper. For an elliptic generator fixing edge $i$, we need that

\begin{equation}
    x_i=\dfrac{2\pi}{k}, k\in\mathbb{Z}
\end{equation} so that some power of the elliptic is the identity. We could explore how this equation combines with the Rivin equations and the fundamental domain equations used in this paper in a manner similar to our methods in Section 4.
Another possible extension would be to classify fundamental domains on the octahedron, analogous to what we did for the cube in Section 5. It would be interesting to analyze the similarities between the cube and the octahedron, given that they are polyhedral duals of one another.

We also proved a general result about the number of edge classes required for a given polyhedron based on the Euler characteristic equation, assuming that the associated group is torsion-free. Using that result, we know that the cube has to have 2 edge classes, and we classify all of the fundamental domains on the cube. In verifying that the candidates that we found are fundamental domains, we describe the process by which we can explicitly find the group associated with a set of edge classes. Finally, we prove a variety of results about combinatorial restrictions on general fundamental domains; specifically, we investigated how different properties of the relators in a group lead to restrictions on the arrangement of the edge classes on the polyhedron. All of these results give us insight into what polyhedra and groups can be fundamental domains. They also illuminate the connection between the group and polyhedron associated to a fundamental domain. Armed with these results, we move towards a method of determining whether or not certain polyhedra can be fundamental domains, bringing us a step closer to understanding hyperbolic 3-manifolds.

\section{Acknowledgements}
A huge thank you to Franco Vargas Pallete, my advisor on this project, for answering my questions and helping me develop my interests in this area. Thank you also to University of California, Berkeley for the REU funding and experience during the summer of 2018, and to the other students that participated in the program for their support and collaboration.

\section{References}
[1] Albert Marden. (2007). \textit{Outer circles: An introduction to hyperbolic 3-manifolds}. Cambridge, UK: Cambridge University Press.

[2] Rivin, I. (1996). A characterization of ideal polyhedra in hyperbolic 3-space. \textit{Annals of Mathematics}, 143, 51-70.

[3]  E. M. Andre'ev, On convex polyhedra of finite volume in Lobachevskii space, Math. USSR, Sbornik, 12(1970), 255-259.
\end{document}